\documentclass{article}
\usepackage[english]{babel}
\usepackage[utf8]{inputenc}

\usepackage{bm}
\usepackage{amsmath}
\usepackage{graphicx}
\usepackage{color}
\usepackage{subcaption}
\usepackage{amsfonts}
\usepackage{float}
\usepackage{placeins}
\usepackage{algorithm}
\usepackage{algpseudocode}
\usepackage{amsmath}
\usepackage{hyperref}
\usepackage{lineno}
\RequirePackage{booktabs}
\usepackage{verbatim}

\graphicspath{ {./figs/} }

\title{Multicontinuum Homogenization for Coupled Flow and Transport Equations}
\author{Dmitry Ammosov\footnote{Computational Technologies \& Artificial Intelligence, North-Eastern Federal University, Yakutsk 677980, Russia}, W.T. Leung\footnote{Department of Mathematics, City University of Hong Kong, Hong Kong}, Buzheng Shan \footnote{Department of Mathematics, Texas A\&M University, College Station, TX 77843, USA}, Jian Huang\footnote{School of Mathematics and Computational Science, Xiangtan University, National Center for Applied Mathematics in Hunan and Hunan Key Laboratory for Computation and Simulation in Science and Engineering, Xiangtan 411105, Hunan, China}}

\begin{document}

\maketitle

\begin{abstract}
    In this paper, we present the derivation of a multicontinuum model for the coupled flow and transport equations by applying multicontinuum homogenization. 
    We perform the multicontinuum expansion for both flow and transport solutions and formulate novel coupled constraint cell problems to capture the multiscale property, where oversampled regions are utilized to avoid boundary effects. Assuming the smoothness of macroscopic variables, we obtain a multicontinuum system composed of macroscopic elliptic equations and convection-diffusion-reaction equations with homogenized effective properties.
    Finally, we present numerical results for various coefficient fields and boundary conditions to validate our proposed algorithm.

    \textbf{Keywords:} multiscale, multicontinuum, homogenization, coupled flow and transport, miscible displacement.
\end{abstract}

\section{Introduction}

The coupled system of flow and transport equations is fundamental in modeling phenomena across various scientific and engineering fields, including hydrology, petroleum engineering, and environmental science. 
For instance, the miscible displacement of one incompressible fluid by another, which occurs in processes ranging from enhanced oil and gas recovery to the remediation of contaminated groundwater, can be governed by these equations \cite{stalkup1983miscible}. 
However, a ubiquitous challenge in these scenarios is the inherent heterogeneities and high contrast in material properties. These complexities not only complicate simulation efforts but also significantly impact the accuracy of predictive models. Consequently, the adoption of upscaling techniques is necessary to accurately capture detailed multiscale properties, enabling time-saving and precise computations.

Many different methods have been developed to perform simulations on a coarse grid, one type of which is homogenization methods \cite{papanicolau1978asymptotic,wu2002analysis,hornung2012homogenization}. 
In homogenization methods, one defines upscaled coefficients for multiscale properties through local solutions of fine-grid problems and utilizes them to formulate the macroscopic flow and transport equations on the coarse grid. The homogenized solution for the transport equation can be computed after the macroscopic flow equation has been solved. Such homogenization methods have been applied in petroleum engineering \cite{durlofsky1997nonuniform,efendiev2000modeling}.

Another approach is multiscale methods, including the Multiscale Finite Element Method \cite{hou1997multiscale, hou1999convergence, efendiev2009multiscale}, the Multiscale Finite Volume Method \cite{jenny2003multi, jenny2005adaptive, jenny2006adaptive}, the Generalized Multiscale Finite Element Method \cite{efendiev2013generalized,chung2015mixed,chung2016adaptive}, and others.
Multiscale methods construct local basis functions via fine-grid local problems and use them to approximate the solutions. One could either apply multiscale methods to the flow equation and solve the transport equation on the fine grid, or apply these methods to both of them, depending on the media properties \cite{jenny2005adaptive, chung2018multiscale}. 
Mixed multiscale methods have also been developed to solve the coarse transport equation \cite{chen2003mixed, aarnes2008multiscale}.

This paper presents the multicontinuum homogenization method following earlier works on the elliptic equation \cite{efendiev2023multicontinuum} using the ideas of the above two types of methods.
For deriving a macroscopic equation, multicontinuum homogenization posits a multicontinuum expansion in each coarse block (or Representative Volume Element (RVE) in the case of scale separation) about macroscopic quantities and multiscale basis functions. The macroscopic quantities are introduced to represent local averages of solution in each continuum and assumed to be smooth over coarse regions. Besides, the multiscale basis functions are obtained from some oversampled locally constraint cell problems solved on a fine grid, which are inspired by \cite{chung2018constraint,chung2018non} and can be designed to promise localization.

In our study, we will apply multicontinuum homogenization to both the flow and the transport equations in the coupled system. However, challenges are that the transport equation contains heterogeneous and high-contrast convection and temporal terms, and its convection coefficient further depends on the solution of the flow equation. To address these problems, we carefully design new coupled constraint cell problems, where the one for transport reflects the multiscale property for each term. 
With multiscale basis functions obtained from the coupled cell problems, we perform the multicontinuum expansions in the variational forms of the flow and transport equations. We will derive a multicontinuum system composed of macroscopic elliptic equations and convection-diffusion-reaction equations, where the effective properties have been homogenized in each coarse block to capture the fine-scale details. We remark here that the properties for convection-diffusion-reaction equations will depend on the flow equation. Therefore, machine learning techniques could be used to predict effective properties, and we could, as in many model reduction methods, divide the computation into offline and online stages to save online computation. We also note that the fine-scale information can be recovered from macroscopic solutions by the multicontinuum expansions. We apply it to coupled flow and transport equations with different coefficient fields subject to different boundary conditions to validate our proposed method. We change the coarse mesh size and compare the macroscopic results with the averaged reference solutions. The relative $L^2$ errors are very small, and our method performs the simulations accurately.

Our main contributions consist of the following. 
Firstly, we formulate appropriate coupled constraint cell problems to deal with the multiscale property of the coupled equations. 
Secondly, we derive a novel multicontinuum coupled model about macroscopic solutions on the coarse grid.
Thirdly, we present numerical results for various media fields and boundary conditions to show the validity.

The paper is organized as follows. In Section \ref{sec:preliminaries}, we present preliminaries and review some earlier works on multiscale methods. Section \ref{sec:multicontinuum} is contributed to the description of multicontinuum homogenization and the derivation of the multicontinuum model. We present numerical results in Section \ref{sec:numerical} and finally draw conclusions in Section \ref{sec:conclusion}.

\section{Preliminaries}
\label{sec:preliminaries}

In this work, we consider the following linearized and one-way coupled system of flow and transport equations for pressure $p$ and concentration $c$ in a bounded domain $\Omega \subset \mathbb{R}^d$
\begin{equation}\label{eq:fine_model}
\begin{split}
&-\nabla \cdot (\kappa \nabla p) = g,\\
&\phi \frac{\partial c}{\partial t} + u \cdot \nabla c - \nabla \cdot (D \nabla c) = h(c),
\end{split}
\end{equation}
where $u = - \kappa \nabla p$, $\kappa$ is a permeability coefficient, $\phi$ is a porosity, $D$ is a diffusion coefficient, and $g$ and $h$ are source terms. We assume $\kappa$ and $D$ are heterogeneous and of high contrast.
Before delving into multicontinuum homogenization, we briefly review earlier efforts on some methods for multiscale problems, including (numerical) homogenization, MsFEM, GMsFEM.

The main idea of the homogenization method is to derive macroscopic equations, with the multiscale coefficients in original equations replaced by homogenized effective properties. While the classical homogenization method relies on periodicity and posits a formal asymptotic expansion of solution, numerical homogenization defines an upscaled property for each coarse block, conversing averages in a physical sense for specific boundary conditions, without the periodicity assumption. To homogenize the flow equation in (\ref{eq:fine_model}), for instance, we formulate local problems for each coarse block $K$
\begin{equation}
    -\nabla \cdot (\kappa \nabla \phi_l) = 0,
\end{equation}
subject to some appropriate boundary conditions, like $\phi_l = x_l$ on $\partial K$. The elements of the upscaled coefficient $\kappa^{* }$ are defined by averaging the fluxes 
\begin{equation}
    \int_K \kappa_{i j}^{* } \frac{\partial}{\partial x_j} \phi_l^*=\int_K \kappa_{i j}(x) \frac{\partial}{\partial x_j} \phi_l,
\end{equation}
and this leads to 
\begin{equation}
    \kappa_{\cdot,l}^{* } = \frac{1}{|K|} \int_K \kappa (x) \nabla \phi_l.
\end{equation}
Then, the transport equation can be homogenized with a similar idea to get upscaled diffusion $D^{*}$ while the upscaled porosity $\phi^{*}$ would be its volume averaging, and consequently, we derive the following coupled macroscopic equations
\begin{equation}
\begin{split}
&-\nabla \cdot (\kappa^{*} \nabla p^{*}) = g,\\
&\phi^{*} \frac{\partial c^{*}}{\partial t} + u^{*} \cdot \nabla c^{*} - \nabla \cdot (D^{*} \nabla c^{*}) = h,
\end{split}
\end{equation}
where $u^{*} = - \kappa^{*} \nabla p^{*}$, $\kappa^{* }$ and $D^{*}$ are anisotropic upscaled coefficients. The macroscopic solutions can be obtained by solving the coupled system above on the coarse grid.

The Multiscale Finite Element Method (MsFEM), rather than deriving a macroscopic equation, establishes local multiscale basis functions via local problems formulated on coarse blocks to approximate the solution. Specifically, with multiscale basis functions $\phi^{\omega_i,p}$ for pressure and $\phi^{\omega_i,c}$ for concentration at each coarse node $i$ supported in coarse neighborhood $\omega_i$, we seek approximate solutions $p_H$ and $c_H$ that lie in $V_\text{ms}^p = \text{span}_i \{\phi^{\omega_i,p}\}$ and $V_\text{ms}^c = \text{span}_i \{\phi^{\omega_i,c}\}$ respectively, and
\begin{equation}
\label{eq:MsFEM approx}
    p_H = \sum_i p_i \phi^{\omega_i,p},\quad c_H = \sum_i c_i \phi^{\omega_i,c}.
\end{equation}
Now, we illustrate how to construct local problems. Since the transport equation in (\ref{eq:fine_model}) lacks temporal heterogeneity, $\phi^{\omega_i,p}$ and $\phi^{\omega_i,c}$ can be chosen to satisfy
\begin{equation}
\begin{split}
    & -\nabla \cdot (\kappa \nabla \phi^{\omega_i,p}) = 0,\\
    & (-\kappa \sum_i p_i \nabla \phi^{\omega_i,p}) \cdot \nabla \phi^{\omega_i,c} - \nabla \cdot (D \nabla \phi^{\omega_i,c}) = 0
\end{split}
\end{equation}
in each coarse block K such that $\phi^{\omega_i,p} = \phi^{\omega_i,c} = \phi_0^{\omega_i}$ on $\partial K$, where $\phi_0^{\omega_i}$ represents the nodal basis at node $i$ of the coarse finite element space. Indeed, the mesh grids for pressure and concentration are not necessarily the same, and the  boundary conditions for local problems can be chosen differently. Note also that local problems can be defined in representative volume element (RVE) in the case of scale separation, and fine-scale information could be recovered, which is unachievable by most homogenization methods.

The Generalized Multiscale Finite Element Method (GMsFEM), evolving from MsFEM, adds additional local degrees of freedom as necessary and thus can systematically enrich coarse spaces and capture multiscale features more precisely. The first step is to construct the snapshot space $V_\text {snap }^{(i),p}$ and $V_\text {snap }^{(i),c}$ for each generic coarse region $\omega_i$, which can be composed of, for example, harmonic extensions of fine-grid functions on $\partial \omega^i$, that is,
\begin{equation}
\begin{split}
    & -\nabla \cdot (\kappa \nabla \psi_l^{\omega_i,p}) = 0,\\
    & (-\kappa \nabla p_H) \cdot \nabla \psi_l^{\omega_i,c} - \nabla \cdot (D \nabla \psi_l^{\omega_i,c}) = 0
\end{split}
\end{equation}
in $\omega_i$ with $\psi_l^{\omega_i,p} = \psi_l^{\omega_i,c} = \delta _l^h$ on $\partial \omega^i$. Here $\delta _l^h$ is a piecewise linear function on $\partial \omega^i$ satisfying $\delta _l^h(x_k) = \delta_{l,k}$ for any fine-grid boundary node $x_k$. Subsequently, we only keep the dominant component subspace of each local snapshot space and form a multiscale space. A local spectral decomposition for each snapshot space $V_\text {snap }^{(i)}$ is constructed for model reduction purposes. This involves finding $\lambda_j^{\omega_i} \in \mathbb{R}$ and $\phi_j^{\omega_i} \in V_\text {snap }^{(i)}$ such that
\begin{equation}
    a_{\omega_i}(\phi_j^{\omega_i},w) = \lambda_j^{\omega_i} s_{\omega_i}(\phi_j^{\omega_i},w),\quad \forall w \in  V_\text {snap }^{(i)}.
\end{equation}
\begin{comment}
where $a_{\omega_i}(v,w) = \int_{\omega_i} \kappa \nabla v \cdot \nabla w$ and $s_{\omega_i}(v,w) = \int_{\omega_i} \tilde{\kappa} vw$ for the diffusion equation for example.
\end{comment}
We select the eigenvectors corresponding to the smallest eigenvalues, denoted as $\phi_l^{\omega_i,p}$ and $\phi_m^{\omega_i,c}$, and multiply them by multiscale partition of unity $\chi_i$ to form the multiscale basis, i.e., $V_\text{ms}^p = \text{span}_{i,l}\{ \chi_i \phi_l^{\omega_i,p} \}$ and $V_\text{ms}^c = \text{span}_{i,m}\{ \chi_i \phi_m^{\omega_i,c} \}$. We remark that the degrees of freedom can be added adaptively according to a-posteriori error estimates. For more details on the methods above, one may refer to \cite{hornung2012homogenization, hou1997multiscale, efendiev2013generalized, chung2023multiscale}.

\section{Multicontinuum homogenization and multicontinuum model}
\label{sec:multicontinuum}

In this section, we will present the main ideas and concepts of multicontinuum homogenization and apply it to the coupled flow and transport equations to derive a multicontinuum model.

\subsection{Main ideas}

The multicontinuum homogenization method, similar to classical homogenization, also aims for macroscopic equations. It introduces macroscopic quantities to represent the local averages of solution and posits the general multicontinuum expansion in each local region, like RVE,
\begin{equation}
\label{eq:mc_expansion}
    u = \phi_i U_i + \phi_i^m \nabla_m U_i + \phi_i^{mn} \nabla_{mn}^2 U_i + \cdots,
\end{equation}
where $U_i$ is the macroscopic quantity for continuum $i$ with physical meaning, $\phi_i$, $\phi_i^m$, $\phi_i^{mn}$, $\ldots$ are auxiliary basis functions obtained via some carefully constructed cell problems, which minimize local energy under certain constraints. The summation over repeated indices is taken. Here, the choice of cell problems is crucial for ensuring localization of the local basis. In the following, we will ignore the terms after the second term and only use the two-term expansion for simplicity. We apply the above expansion to both the solution and test function and substitute them into the variational form, where RVE concepts could be utilized. Using the smoothness of macroscopic quantities, we can take them out of the local integrals and, after some calculation, arrive at a macroscopic equation about $U_i$, which can recover the fine-scale information by expansion (\ref{eq:mc_expansion}).
Note that the smoothness of macroscopic quantities is reasonable and crucial for multicontinuum homogenization. We also remark that the continua are not necessarily defined apriori but could be achieved by local spectral problems.

Next, we will describe the multicontinuum homogenization method for (\ref{eq:fine_model}) in detail. Since the pressure field's equation does not depend on the concentration, we can perform homogenization in a split way.

\subsection{Flow homogenization}

Let us consider the following variational formulation of the flow equation
\begin{equation}\label{eq:p_var_model}
\int_\Omega \kappa \nabla p \cdot \nabla q = \int_\Omega g q, \quad \forall q \in H_0^1(\Omega).
\end{equation}

\paragraph{Cell problems.}

We assume that our computational domain $\Omega$ is partitioned into coarse-grid elements $\omega$'s, whose sizes are larger than the scale of heterogeneities. Also, we assume there is a Representative Volume Element (RVE) $R_\omega$ inside each coarse element $\omega$. These RVEs have the following properties:
\begin{itemize}
\item $R_\omega$ can represent the entire $\omega$ in terms of its heterogeneities;
\item $R_\omega$ contains N continua (or components). For each continuum of $R_\omega$, we introduce the following characteristic function $\psi_j$ ($j = 1, ..., N$)
\end{itemize}
\[
\psi_j = \begin{cases}
  1, & \text{in continuum} ~j,\\
  0, & \text{otherwise}.
\end{cases}
\]

\begin{figure}
    \centering
    \includegraphics[width=8cm]{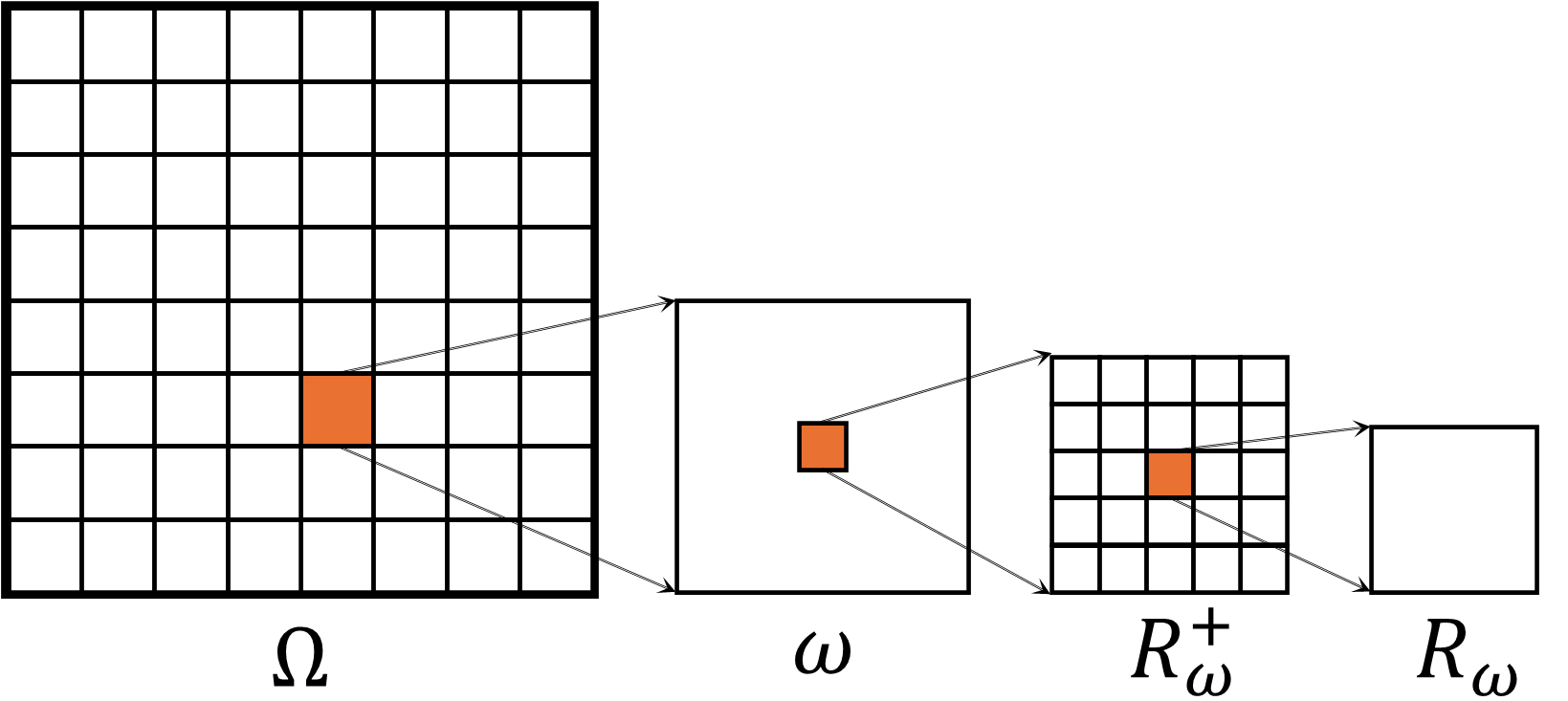}
    \caption{Illustration of computational domain $\Omega$, coarse block $\omega$, RVE $R_\omega$ and oversampled RVE $R_\omega^+$}
    \label{fig:Illustration}
\end{figure}

For performing homogenization, we introduce two cell problems in the oversampled RVE $R_\omega^+$ constructed around $R_\omega$ and consisted of several RVEs $R_\omega^k$, where $R_\omega^{k_0} = R_\omega$. The relationship between the computational domain $\Omega$, coarse block $\omega$, RVE $R_\omega$ and oversampled RVE $R_\omega^+$ has been illustrated in Figure \ref{fig:Illustration}. We would like to minimize the local energy in ${R_\omega^+}$ under different constraints, where Lagrange multipliers are applied. The first cell problem \eqref{eq:p_cell_grad_aver} considers gradient effects and imposes constraints to represent the linear functions in the average behavior of each continuum
\begin{equation}
\label{eq:p_cell_grad_aver}
\begin{split}
&\int_{R_\omega^+}\kappa \nabla \varphi^{m,p}_{i}\cdot \nabla q -  \sum_{j,k} {\beta_{ij}^{mk}\over \int_{R_\omega^k}\psi_j^k} \int_{R_\omega^k}\psi_j^k q =0,\\
&\int_{R_\omega^k}  \varphi^{m,p}_{i} \psi_j^k = \delta_{ij} \int_{R_\omega^k} (x_m-\tilde{x}_{m})\psi_j^k ,\\
&\int_{R_\omega^{k_0}} (x_m-\tilde{x}_{m})\psi_j^{k_0} =0, \\
\end{split}
\end{equation}
where $\tilde{x}_{m}$ is a constant.
The second cell problem considers different averages in each continuum and represents the constants in the average behavior
\begin{equation}
\label{eq:p_cell_aver_aver}
\begin{split}
&\int_{R_\omega^+}\kappa \nabla \varphi_{i}^{p} \cdot \nabla q -
\sum_{j,k} {\beta_{ij}^k\over \int_{R_\omega^k}\psi_j^k} \int_{R_\omega^k}\psi_j^k  q = 0,\\
&\int_{R_\omega^k}  \varphi_{i}^{p} \psi_j^k = \delta_{ij} \int_{R_\omega^k} \psi_j^k.
\end{split}
\end{equation}
We remark that the solution in $R_{\omega}$ should be independent of the oversampling size as a result of the construction of our cell problems.

\paragraph{Coarse-grid model.}

Next, we introduce the following multicontinuum expansion of $p$ in each $R_\omega$
\begin{equation}
\label{eq:p_mc_expansion_p}
\begin{split}
p \approx  \varphi_i^p P_i  + \varphi_i^{m,p} \nabla_m P_i,
\end{split}
\end{equation}
where $P_i$ is a smooth function in a macroscopic sense representing the homogenized solution for the $i$-th continuum. In our case, we have $P_i(x_{\omega}^*) = \dfrac{\int_{R_{\omega}} p \psi_i }{\int_{R_{\omega}} \psi_i}$ for some $x_{\omega}^* \in R_{\omega}$.

Using the definition of RVE, we obtain the following approximation of the variational formulation
\begin{equation}
\label{eq:p_RVE_assump}
\begin{split}
&\int_\Omega g q = \int_\Omega \kappa \nabla p \cdot \nabla q = \sum_\omega\int_\omega \kappa \nabla p \cdot \nabla q \approx \sum_\omega {|\omega|\over |R_\omega|} \int_{R_\omega} \kappa \nabla p \cdot \nabla q.
\end{split}
\end{equation}
Next, if we replace $p$ with its multicontinuum expansion \eqref{eq:p_mc_expansion_p}, we get
\begin{equation}\label{eq:p_mc_approx}
\begin{split}
\int_{R_\omega} \kappa \nabla p \cdot \nabla q \approx
\int_{R_\omega} \kappa \nabla (\varphi_i^p P_i) \cdot\nabla q  +
\int_{R_\omega} \kappa \nabla (\varphi_i^{m,p} \nabla_m P_i) \cdot\nabla q. 
\end{split}
\end{equation}
Additionally, we define the following multicontinuum expansion of $q$
\begin{equation}\label{eq:p_mc_expansion_v}
q \approx \varphi_s^p Q_s + \varphi_s^{k,p} \nabla_k Q_s.
\end{equation}
According to our assumption on the smoothness of the macroscopic variables, the variations of $P_i$ and $\nabla_m P_i$ are relatively small when compared to the variations of $\varphi_i^p$ and $\varphi_i^{m,p}$. Similarly, this applies to $Q_s$ and $\nabla_k Q_s$.

Then, we can approximate the first term of \eqref{eq:p_mc_approx} as follows
\begin{equation}
\begin{split}
&\int_{R_\omega} \kappa \nabla (\varphi_i^p P_i) \cdot \nabla q \approx
P_i(x_{\omega}) \int_{R_\omega} \kappa \nabla \varphi_i^p \cdot \nabla q \approx\\
&P_i(x_{\omega}) Q_s(x_\omega)\int_{R_\omega} \kappa \nabla \varphi_i^p \cdot \nabla \varphi_s^p +
P_i(x_{\omega}) \nabla_m Q_s(x_\omega)\int_{R_\omega} \kappa \nabla \varphi_i^p \cdot \nabla \varphi_s^{m,p} =\\
&P_i(x_{\omega})  \beta_{is}^* Q_s(x_\omega) + \beta_{is}^{m*} P_i(x_{\omega}) \nabla  Q_s(x_\omega),
\end{split}
\end{equation}
where
\begin{equation}\label{eq:p_beta_defintions}
\beta_{ik}^{m*}=\int_{R_\omega} \kappa \nabla \varphi_i^p \cdot \nabla \varphi_k^{m,p}, \ \ 
\beta_{ik}^{*}=\int_{R_\omega} \kappa \nabla \varphi_i^p \cdot \nabla \varphi_k^p.
\end{equation}
The second term of \eqref{eq:p_mc_approx} is approximated in a similar way
\begin{equation}
\begin{split}
&\int_{R_\omega} \kappa \nabla (\varphi_i^{m,p} \nabla_m P_i)  \cdot \nabla q\approx
\nabla_m P_i (x_\omega)\int_{R_\omega} \kappa \nabla \varphi_i^{m,p}\cdot \nabla q \approx\\
&\nabla_m P_i (x_\omega)  \nabla_k Q_s (x_\omega) \int_{R_\omega} \kappa \nabla \varphi_i^{m,p} \cdot \nabla \varphi_s^{k,p} +
\nabla_m P_i (x_\omega)  Q_s (x_\omega) \int_{R_\omega} \kappa \nabla \varphi_i^{m,p} \cdot \nabla \varphi_s^p =\\
&\nabla_m P_i (x_\omega)  \nabla_k Q_s (x_\omega) \alpha_{is}^{km} + \nabla_m P_i (x_\omega) Q_s (x_\omega) \beta_{is}^{m*},
\end{split}
\end{equation}
where
\begin{equation}
    \alpha_{is}^{km}= \int_{R_\omega} \kappa \nabla \varphi_i^{m,p} \cdot \nabla \varphi_s^{k,p}.
\end{equation}
Therefore, by utilizing continuous approximations for $P_i$ and $Q_j$, we have
\begin{equation}\label{eq:p_mc_final}
\begin{split}
\int_{R_\omega} \kappa \nabla p  \cdot \nabla q \approx
P_i  \beta_{ij}^{n*} \nabla_n Q_j + P_i  \beta_{ij}^* Q_j +
\nabla_m P_i\alpha_{ij}^{mn} \nabla_n Q_j + \nabla_m P_i \beta_{ij}^{m*} Q_j.
\end{split}
\end{equation}

According to \cite{efendiev2023multicontinuum}, we have the following estimates for the solutions of the local cell problems and their gradients
\begin{equation}
\label{eq:p_mc_scalings}
\begin{split}
&\|\varphi_i^p\|=O(1),\ \|\nabla \varphi_i^p\|=O({1\over\epsilon}),\\
&\|\varphi_i^{m,p}\|=O(\epsilon),\ \|\nabla \varphi_i^{m,p}\|=O(1),\\
\end{split}
\end{equation}
where $\epsilon$ is the size of RVE.
Then, we can obtain the following scalings
\begin{equation*}
\beta_{ij}^* =O(\cfrac{|R_\omega|}{\epsilon^{2}}),\
\beta_{ij}^{m*}=O(\cfrac{|R_\omega|}{\epsilon}),\ 
\alpha_{ij}^{mn}=O(|R_\omega|).\
\end{equation*}
We introduce $\widehat{\beta_{ij}}$, $\widehat{\beta_{ij}^{m*}}$, and $\widehat{\alpha_{ij}^{mn}}$ as follows
\begin{equation}\label{eq:p_rescaled_coefficients}
\begin{split}
\widehat{\beta_{ij}}  =\cfrac{\epsilon^{2}}{|R_\omega|}\beta_{ij}^*,\ 
\widehat{\beta_{ij}^{m*}}  =\cfrac{\epsilon}{|R_\omega|}\beta_{ij}^{m*},\ 
\widehat{\alpha_{ij}^{mn}} = \frac{1}{|R_\omega|}\alpha_{ij}^{mn}.\
\end{split}
\end{equation}
Using these scalings, we obtain
\begin{equation}
\label{eq:p_mc}
\begin{split}
\int_\Omega \kappa \nabla p \cdot \nabla q \approx
\int_\Omega \widehat{\alpha_{ij}^{mn}} \nabla_m P_i \nabla_n Q_j +
{1\over\epsilon}\int_\Omega \widehat{\beta_{ij}^m} \nabla_m P_i Q_j +
{1\over\epsilon}\int_\Omega  \widehat{\beta_{ij}^m} P_i \nabla_m Q_j +
{1\over\epsilon^2}   \int_\Omega \widehat{\beta_{ij}} P_i Q_j.
\end{split}
\end{equation}

By applying the integration by parts formula, one can find out that that the sum of the second and third terms is negligible. Finally, we obtain the following multicontinuum flow equations
\begin{equation}\label{eq:p_mc_index}
\begin{split}
-\nabla_n (\widehat{\alpha_{ij}^{mn}} \nabla_m P_j) +
{1\over \epsilon^2}\widehat{\beta_{ij}} P_j=g_i.
\end{split}
\end{equation}

\subsection{Transport homogenization}

Let us consider the following variational formulation of the transport equation
\begin{equation}\label{eq:original1}
\int_\Omega \phi \frac{\partial c}{\partial t} v + \int_\Omega (u \cdot \nabla c) v + \int_\Omega D \nabla c \cdot \nabla v = \int_\Omega h(c) v, \quad \forall v \in H_0^1(\Omega).
\end{equation}

\paragraph{Cell problems.} To formulate and solve the cell problems, we need to linearize the convective term using known $P_i$. One can obtain these functions by solving the flow problem prior to that. Another approach is to consider many possible variants of $P_i$ to approximate the effective properties using machine learning techniques.
% We need to solve cell problems and compute effective properties for various $P_i$. After that we can use machine learning techniques to predict effective properties. 

Again, the first cell problem \eqref{eq:c_cell_grad_aver} takes into account gradient effects
\begin{equation}
\label{eq:c_cell_grad_aver}
\begin{split}
&\int_{R_\omega^+} \left[-\kappa \nabla ( \varphi_i^p P_i  + \varphi_i^{m,p} \nabla_m P_i)\right] \cdot \nabla \varphi^{m,c}_{i} v + \int_{R_\omega^+} D \nabla \varphi^{m,c}_{i} \cdot \nabla v -  
\sum_{j,k} {\theta_{ij}^{mk}\over \int_{R_\omega^k}\psi_j^k} \int_{R_\omega^k}\psi_j^k v =0,\\
&\int_{R_\omega^k}  \varphi^{m,c}_{i} \psi_j^k = \delta_{ij} \int_{R_\omega^k} (x_m-\tilde{x}_{m})\psi_j^k ,\\
&\int_{R_\omega^{k_0}} (x_m-\tilde{x}_{m})\psi_j^{k_0} =0. \\
\end{split}
\end{equation}
The second cell problem considers different averages in each continuum
\begin{equation}
\label{eq:c_cell_aver_aver}
\begin{split}
&\int_{R_\omega^+} \left[-\kappa \nabla ( \varphi_i^p P_i  + \varphi_i^{m,p} \nabla_m P_i)\right] \cdot \nabla \varphi^{c}_{i} v + \int_{R_\omega^+} D \nabla \varphi_{i}^{c} \cdot \nabla v -
\sum_{j,k} {\theta_{ij}^k\over \int_{R_\omega^k}\psi_j^k} \int_{R_\omega^k}\psi_j^k  v = 0,\\
&\int_{R_\omega^k}  \varphi_{i}^{c} \psi_j^k = \delta_{ij} \int_{R_\omega^k} \psi_j^k.
\end{split}
\end{equation}

\paragraph{Coarse-grid model.}

We apply the following multicontinuum expansion of $c$ in $R_\omega$ using $\varphi_i^c$ and $\varphi_i^{m,c}$ as follows
\begin{equation}
\label{eq:c_mc_expansion}
\begin{split}
c \approx  \varphi_i^c C_i  + \varphi_i^{m, c} \nabla_m C_i,
\end{split}
\end{equation}
where $C_i$ is a smooth function representing the homogenized solution for the $i$-th continuum.

We obtain the following approximation of the variational formulation
\begin{equation}
\label{eq:c_RVE_assump}
\begin{split}
&\int_\Omega \phi {\partial c \over \partial t} v + \int_\Omega (u \cdot \nabla c) v + \int_\Omega D \nabla c \cdot \nabla v = 
\sum_\omega \int_\omega \phi {\partial c \over \partial t} v + \sum_\omega \int_\omega (u \cdot \nabla c) v + \sum_\omega\int_\omega D \nabla c \cdot \nabla v \approx\\
&\sum_\omega {|\omega|\over |R_\omega|} \int_{R_\omega}
 \phi {\partial c \over \partial t} v +
 \sum_\omega {|\omega|\over |R_\omega|} \int_{R_\omega} (u \cdot \nabla c) v +
\sum_\omega {|\omega|\over |R_\omega|} \int_{R_\omega} D \nabla c \cdot \nabla v = \sum_\omega {|\omega|\over |R_\omega|} \int_{R_\omega} h(c) v
\end{split}
\end{equation}

The homogenization of the diffusion term is similar to the flow case. We can represent it as follows
\begin{equation}\label{eq:c_mc_diffusion}
\begin{split}
\int_{R_\omega} D \nabla c  \cdot \nabla v \approx
C_i  \theta_{ij}^{n*} \nabla_n V_j + C_i  \theta_{ij}^* V_j +
\nabla_m C_i \eta_{ij}^{mn} \nabla_n V_j + \nabla_m C_i \theta_{ij}^{m*} V_j,
\end{split}
\end{equation}
where
\begin{equation}\label{eq:c_diffusion_and_transfer}
\eta_{is}^{km}= \int_{R_\omega} D \nabla \varphi_i^{m,c} \cdot \nabla \varphi_s^{k,c}, \ \
\theta_{ik}^{m*}=\int_{R_\omega} D \nabla \varphi_i^c \cdot \nabla \varphi_k^{m,c}, \ \ 
\theta_{ik}^{*}=\int_{R_\omega} D \nabla \varphi_i^c \cdot \nabla \varphi_k^c.
\end{equation}
We can also approximate the term with time derivative 
\begin{equation}
\begin{split}
&\int_{R_\omega} \phi {\partial c \over \partial t} v \approx
\int_{R_\omega} \phi 
{\partial ( \varphi_i^c C_i  + \varphi_i^{m,c} \nabla_m C_i) \over \partial t}
( \varphi_j^c V_j  + \varphi_j^{n,c} \nabla_n V_j)\approx
\int_{R_\omega} \phi {\partial ( \varphi_i^c C_i) \over \partial t} \varphi_j^c V_j \approx\\
&{\partial C_i(x_\omega) \over \partial t} V_j(x_\omega) \int_{R_\omega} \phi \varphi_i^c \varphi_j^c = 
{\partial C_i(x_\omega) \over \partial t} V_j(x_\omega) \gamma_{ij},
\end{split}
\end{equation}
where
\begin{equation*}
\gamma_{ij} = \int_{R_\omega} \phi \varphi_i^c \varphi_j^c.
\end{equation*}
Here for the second approximate equality, we use the fact that $\varphi_i^{m,c}$ is of the order $\epsilon$ and $\varphi_i^{c}$ is of the order $1$ as will be stated in (\ref{eq:c_scalings}).
Again we make use of continuous approximations for $C_i$ and $V_j$, allowing us to derive the following approximation
\begin{equation}\label{eq:c_mc_time_appr_final}
\int_{R_\omega} \phi {\partial c \over \partial t} v \approx
{\partial c_i \over \partial t} \gamma_{ij} v_j.
\end{equation}

Next, we need to perform homogenization for the convection term. With our previous assumptions and estimates for the solutions of cell problems, we write
\begin{equation}
\begin{split}
&\int_{R_\omega} (u \cdot \nabla c) v \approx 
\int_{R_\omega} \left[ -\kappa \nabla (\varphi_s^p P_s  + \varphi_s^{l,p} \nabla_l P_s) \cdot \nabla (\varphi_i^c C_i  + \varphi_i^{m,c} \nabla_m C_i) \right] \left[ \varphi_j^c V_j  + \varphi_j^{n,c} \nabla_n V_j \right] \approx \\
&\int_{R_\omega} \left[ -\kappa \nabla (\varphi_s^p P_s  + \varphi_s^{l,p} \nabla_l P_s) \cdot \nabla (\varphi_i^c C_i  + \varphi_i^{m,c} \nabla_m C_i) \right] \varphi_j^c V_j =\\
&\int_{R_\omega} (-\kappa) (\nabla \varphi_s^p P_s \cdot \nabla \varphi_i^c C_i) \varphi_j^c V_j +
\int_{R_\omega} (-\kappa) (\nabla \varphi_s^p P_s \cdot \nabla \varphi_i^{m,c} \nabla_m C_i) \varphi_j^c V_j +\\
&\int_{R_\omega} (-\kappa) (\nabla \varphi_s^{l,p} \nabla_l P_s \cdot \nabla \varphi_i^c c_i) \varphi_j^c V_j +
\int_{R_\omega} (-\kappa) (\nabla \varphi_s^{l,p} \nabla_l P_s \cdot \nabla \varphi_i^{m,c} \nabla_m C_i) \varphi_j^c V_j \approx\\
&P_s(x_\omega) C_i(x_\omega) V_j(x_\omega) \int_{R_\omega} (-\kappa) (\nabla \varphi_s^p \cdot \nabla \varphi_i^c) \varphi_j^c +
P_s(x_\omega) \nabla_m C_i(x_\omega) V_j(x_\omega) \int_{R_\omega} (-\kappa) (\nabla \varphi_s^p \cdot \nabla \varphi_i^{m,c}) \varphi_j^c +\\
&\nabla_l P_s(x_\omega) C_i(x_\omega) V_j(x_\omega) \int_{R_\omega} (-\kappa) (\nabla \varphi_s^{l,p} \cdot \nabla \varphi_i^c) \varphi_j^c +
\nabla_l P_s(x_\omega) \nabla_m C_i(x_\omega) V_j(x_\omega) \int_{R_\omega} (-\kappa) (\nabla \varphi_s^{l,p} \cdot \nabla \varphi_i^{m,c}) \varphi_j^c =\\
&P_s(x_\omega) C_i(x_\omega) V_j(x_\omega) \zeta_{sij} +
P_s(x_\omega) \nabla_m C_i(x_\omega) V_j(x_\omega) \chi_{sij}^{m} +\\
&\nabla_l P_s(x_\omega) C_i(x_\omega) V_j(x_\omega) \Upsilon_{sij}^{l} +
\nabla_l P_s(x_\omega) \nabla_m C_i(x_\omega) V_j(x_\omega) \iota_{sij}^{lm},
\end{split}
\end{equation}
where
\begin{equation}
\begin{split}
&\zeta_{sij} = \int_{R_\omega} (-\kappa) (\nabla \varphi_s^p \cdot \nabla \varphi_i^c) \varphi_j^c, \quad
\chi_{sij}^{m} = \int_{R_\omega} (-\kappa) (\nabla \varphi_s^p \cdot \nabla \varphi_i^{m,c}) \varphi_j^c, \\
&\Upsilon_{sij}^{l} = \int_{R_\omega} (-\kappa) (\nabla \varphi_s^{l,p} \cdot \nabla \varphi_i^c) \varphi_j^c, \quad
\iota_{sij}^{lm} = \int_{R_\omega} (-\kappa) (\nabla \varphi_s^{l,p} \cdot \nabla \varphi_i^{m,c}) \varphi_j^c.
\end{split}
\end{equation}
Applying the continuous approximations for $C_i$ and $V_j$ gives us
\begin{equation}
\begin{split}
&\int_{R_\omega} (u \cdot \nabla c) v \approx P_s C_i V_j \zeta_{sij} +
P_s \nabla_m C_i V_j \chi_{sij}^{m} + \nabla_l P_s C_i V_j \Upsilon_{sij}^{l} +
\nabla_l P_s \nabla_m C_i V_j \iota_{sij}^{lm}.
\end{split}
\end{equation}

We can establish the following estimates
\begin{equation}
\label{eq:c_scalings}
\begin{split}
&\|\varphi_i^c\|=O(1),\ \|\nabla \varphi_i^c\|=O({1\over\epsilon}),\\
&\|\varphi_i^{m,c}\|=O(\epsilon),\ \|\nabla \varphi_i^{m,c}\|=O(1).\\
\end{split}
\end{equation}
Subsequently, we derive the following scalings
\begin{equation*}
\begin{split}
&\theta_{ij}^* =O(\cfrac{|R_\omega|}{\epsilon^{2}}),\
\theta_{ij}^{m*}=O(\cfrac{|R_\omega|}{\epsilon}),\ 
\eta_{ij}^{mn}=O(|R_\omega|),\
\gamma_{ij}=O(|R_\omega|),\\
&\zeta_{sij} = O(\cfrac{|R_\omega|}{\epsilon^{2}}),\ 
\chi_{sij}^{m} = O(\cfrac{|R_\omega|}{\epsilon}),\ 
\Upsilon_{sij}^{l} = O(\cfrac{|R_\omega|}{\epsilon}),\ 
\iota_{sij}^{lm} = O(|R_\omega|).
\end{split}
\end{equation*}
We introduce the scaled effective properties
\begin{equation}\label{eq:c_rescaled_coefficients}
\begin{split}
&\widehat{\theta_{ij}}  =\cfrac{\epsilon^{2}}{|R_\omega|}\beta_{ij}^*,\ 
\widehat{\theta_{ij}^{m*}}  =\cfrac{\epsilon}{|R_\omega|}\beta_{ij}^{m*},\ 
\widehat{\eta_{ij}^{mn}} = \frac{1}{|R_\omega|}\alpha_{ij}^{mn},\
\widehat{\gamma_{ij}} = \frac{1}{|R_\omega|}\gamma_{ij},\\
&\widehat{\zeta_{sij}} = \cfrac{\epsilon^{2}}{|R_\omega|} \zeta_{sij},\
\widehat{\chi_{sij}^{m}} = \cfrac{\epsilon}{|R_\omega|} \chi_{sij}^{m},\ 
\widehat{\Upsilon_{sij}^{l}} = \cfrac{\epsilon}{|R_\omega|}\Upsilon_{sij}^{l},\ 
\widehat{\iota_{sij}^{lm}} = \frac{1}{|R_\omega|} \iota_{sij}^{lm}.
\end{split}
\end{equation}
With these scalings, we have
\begin{equation}
\label{eq:c_mc}
\begin{split}
&\int_\Omega \phi {\partial c \over \partial t} v + \int_\Omega (u \cdot \nabla c) v + \int_\Omega D \nabla c \cdot \nabla v \approx\\
&\int_\Omega \widehat{\gamma_{ij}} {\partial C_i \over \partial t} V_j+
 \int_\Omega \widehat{\eta_{ij}^{mn}} \nabla_m C_i \nabla_n V_j +
{1\over\epsilon}\int_\Omega \widehat{\theta_{ij}^m} \nabla_m C_i V_j +
{1\over\epsilon}\int_\Omega  \widehat{\theta_{ij}^m} C_i \nabla_m V_j +
{1\over\epsilon^2}   \int_\Omega \widehat{\theta_{ij}} C_i V_j +\\
&{1\over\epsilon^2} \int_\Omega P_s C_i V_j \widehat{\zeta_{sij}} +
{1\over\epsilon} \int_\Omega P_s \nabla_m C_i V_j \widehat{\chi_{sij}^{m}} + 
{1\over\epsilon} \int_\Omega \nabla_l P_s C_i V_j \widehat{\Upsilon_{sij}^{l}} +
\int_\Omega \nabla_l P_s \nabla_m C_i V_j \widehat{\iota_{sij}^{lm}}.
\end{split}
\end{equation}

One can see that the sum of the third and fourth terms in \eqref{eq:c_mc} is negligible. It can be shown by using integration by parts. Consequently, we have
\begin{equation}
\label{eq:c_mc_2}
\begin{split}
&\int_\Omega \phi {\partial c \over \partial t} v + \int_\Omega (u \cdot \nabla c) v + \int_\Omega D \nabla c \cdot \nabla v \approx\\
&\int_\Omega \widehat{\gamma_{ij}} {\partial C_i \over \partial t} V_j+
 \int_\Omega \widehat{\eta_{ij}^{mn}} \nabla_m C_i \nabla_n V_j +
{1\over\epsilon^2}   \int_\Omega (\widehat{\theta_{ij}} + P_s \widehat{\zeta_{sij}} + \epsilon \nabla_l P_s \widehat{\Upsilon_{sij}^{l}}) C_i V_j +\\
&{1\over\epsilon} \int_\Omega (P_s \widehat{\chi_{sij}^{m}} + \epsilon \nabla_l P_s \widehat{\iota_{sij}^{lm}}) \nabla_m C_i V_j =\\
&\int_\Omega \widehat{\gamma_{ij}} {\partial C_i \over \partial t} V_j+
 \int_\Omega \widehat{\eta_{ij}^{mn}} \nabla_m C_i \nabla_n V_j +
 {1\over\epsilon} \int_\Omega \widehat{\xi_{ij}^m} \nabla_m C_i V_j +
{1\over\epsilon^2}   \int_\Omega \widehat{\Theta_{ij}} C_i V_j,
\end{split}
\end{equation}
where
\begin{equation*}
\widehat{\xi_{ij}^m} = P_s \widehat{\chi_{sij}^{m}} + \epsilon \nabla_l P_s \widehat{\iota_{sij}^{lm}}, \quad
\widehat{\Theta_{ij}} = \widehat{\theta_{ij}} + P_s \widehat{\zeta_{sij}} + \epsilon \nabla_l P_s \widehat{\Upsilon_{sij}^{l}}.
\end{equation*}
Therefore, we obtain the following multicontinuum equations
\begin{equation}\label{eq:c_mc_index}
\begin{split}
\widehat{\gamma_{ij}}(x, P) {\partial C_j \over \partial t} -
\nabla_n (\widehat{\eta_{ij}^{mn}} (x, P) \nabla_m C_j) +
{1\over\epsilon} \widehat{\xi_{ij}^m}(x, P) \nabla_m C_j +
{1\over \epsilon^2}\widehat{\Theta_{ij}}(x, P) C_j=h_i(x, C, P).
\end{split}
\end{equation}

The effective properties depend on $x$ because they can be different in each coarse block. Moreover, the coefficients of the time term and the source term depend on $P = (P_1, P_2, ..., P_N)$ because the cell problems' solutions $\phi_i^{c}$ and $\phi_i^{m,c}$ depend on $P$.

\subsection{Coupled multicontinuum model}

Finally, we have the following coupled multicontinuum model described by macroscopic elliptic equations and convection-diffusion-reaction equations
\begin{equation}\label{eq:mc_index}
\begin{split}
&-\nabla_n (\widehat{\alpha_{ij}^{mn}} \nabla_m P_j) +
{1\over \epsilon^2}\widehat{\beta_{ij}} P_j=g_i,\\
&\widehat{\gamma_{ij}}(x, P) {\partial C_j \over \partial t} -
\nabla_n (\widehat{\eta_{ij}^{mn}} (x, P) \nabla_m C_j) +
{1\over\epsilon} \widehat{\xi_{ij}^m}(x, P) \nabla_m C_j +
{1\over \epsilon^2}\widehat{\Theta_{ij}}(x, P) C_j=h_i(x, C, P).
\end{split}
\end{equation}
We can obtain the macroscopic solutions by solving (\ref{eq:mc_index}) on the coarse grid. In the next section, we present numerical results to verify the obtained model.

\section{Numerical examples}
\label{sec:numerical}

We will apply our proposed multicontinuum homogenization approach to solve the system of equations ~(\ref{eq:fine_model}) subject to certain boundary conditions, where the spatial domain $\Omega$ is chosen to be a unit square $[0,1] \times [0,1]$. We assume that the permeability coefficient $\kappa$ and the diffusion coefficient $D$ are of high contrast, and we will consider different coefficient fields in the following examples.

We will divide the computational domain $\Omega$ into $M\times M$ square coarse blocks of the same size and define the coarse mesh size to be $H = 1/M$. Each whole coarse block will be taken as an RVE for itself. The oversampling RVE is defined as an extension of itself by $l$ layers of coarse blocks. We denote the regions of low and high values by $\Omega_1$ and $\Omega_2$, respectively, and define the relative $L^2$ errors of solution in $\Omega_1$ and $\Omega_2$ at a specific time $t$ by
\begin{equation}
    e^{(i)}(t)=\sqrt{ \frac{\sum_K|\frac1{|K|}\int_K C_i(x,t) dx -  \frac1{|K\cap \Omega_i|}\int_{K\cap \Omega_i}c(x,t) dx|^2}{\sum_K|\frac1{|K\cap \Omega_i|}\int_{K\cap \Omega_i}c(x,t) dx|^2} },
\end{equation}
where $i=1,2$, and $K$ denotes the RVE, which is taken to be the coarse block.

\subsection{Example 1: Layered field}

\subsubsection{Case 1}

\begin{figure}
    \centering
    \includegraphics[width=8cm]{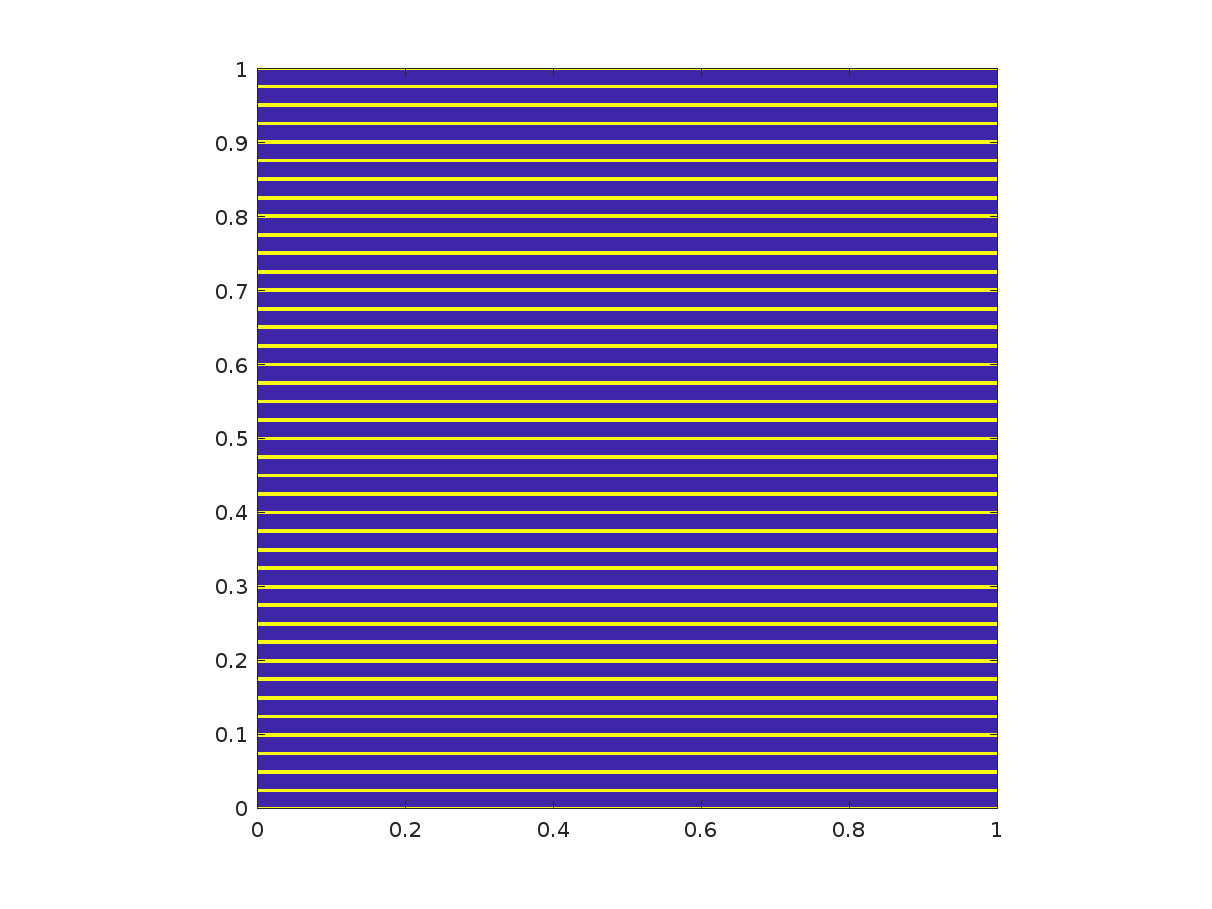}
    \caption{Layered field for Example 1 ($\Omega_1$: blue regions; $\Omega_2$: yellow regions)}
    \label{fig:layer_e40}
\end{figure}

First, we consider the case when the permeability field $\kappa$ and the diffusion field $D$ consist of horizontal layers, as shown in Figure \ref{fig:layer_e40}, and 
\begin{equation}
\label{eq:Example1_coefficient}
    \kappa(x) = D(x) = 
    \begin{cases}
        10^{-4}, \quad & x\in \Omega_1,\\
        1, \quad & x \in \Omega_2.
    \end{cases}
\end{equation}
We assume $\phi = 1$, $g(x)=e^{-40 \left( (x_1-0.5)^2+(x_2-0.5)^2 \right)}$ and $h(x) = 0.1 e^{-40 \left((x_1-0.5)^2+(x_2-0.5)^2 \right) }$ for any $x=(x_1,x_2)\in \Omega$.
We set the following homogeneous Dirichlet boundary conditions
\begin{equation}
    \begin{split}
        p(x,t) = 0, \quad & x\in \partial \Omega, \quad \forall t\in [0,\infty),\\
        c(x,t) = 0, \quad & x\in \partial \Omega, \quad \forall t\in [0,\infty).
    \end{split}
\end{equation}

To discretize the spatial domain, we fix the fine mesh size to be $1/400$ and take the coarse mesh size to be $1/20$ and $1/40$, respectively. We use an implicit Euler scheme for the time discretization and choose the time step $\tau = 0.001$. The number of oversampling layers $l$ is taken to be $ \lceil -2\log(H) \rceil$, that is, $l=6$ for $H=1/20$ and $l= 8$ for $H=1/40$, which is inspired by the analysis in \cite{chung2018constraint, efendiev2023multicontinuum}. 

\begin{figure}
    \centering
    \includegraphics[width=5.5cm]{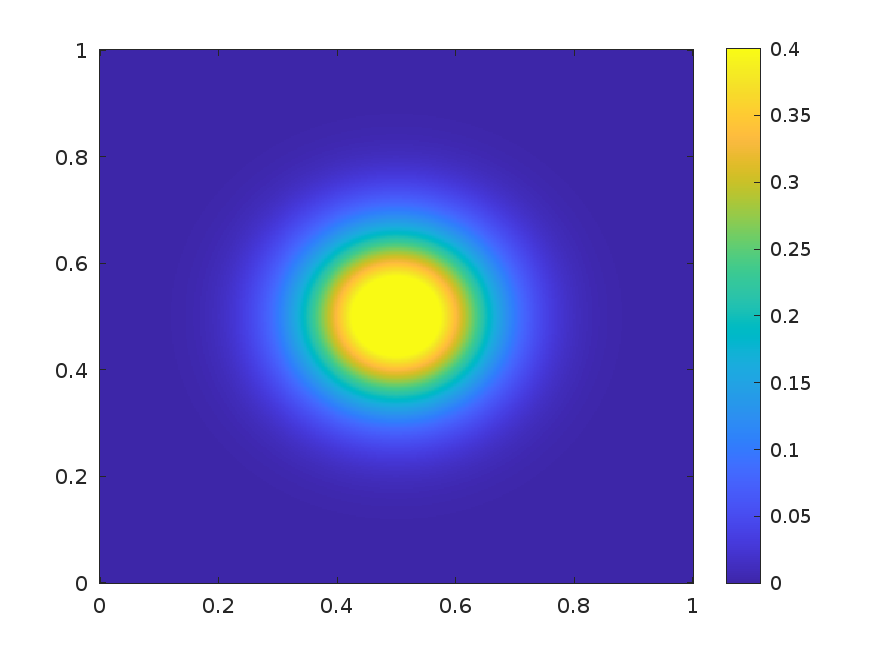}
    \caption{Initial concentration when $t=0$}
    \label{fig:Example1_init}
\end{figure}

\begin{figure}
  \centering
  \includegraphics[width=5.5cm]{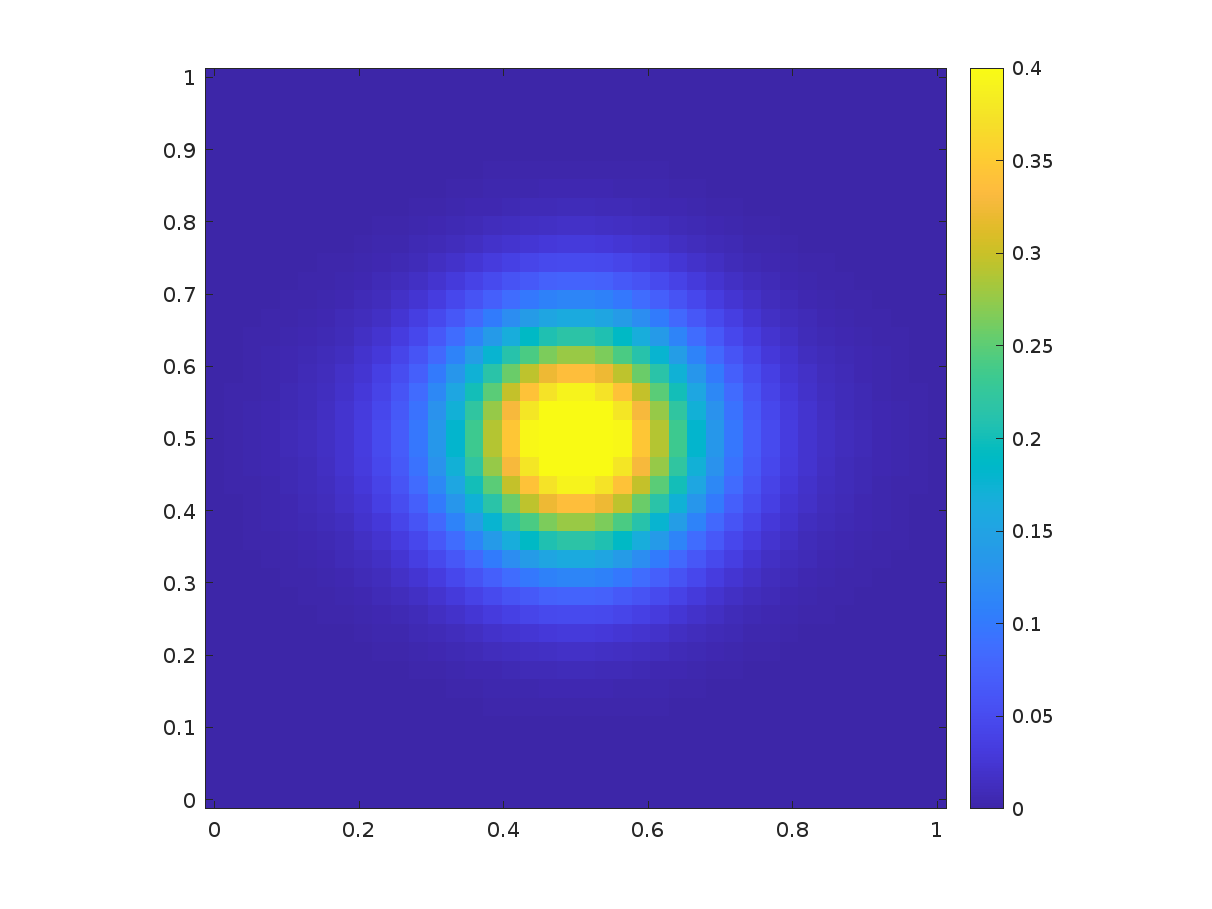}
  \quad
  \includegraphics[width=5.5cm]{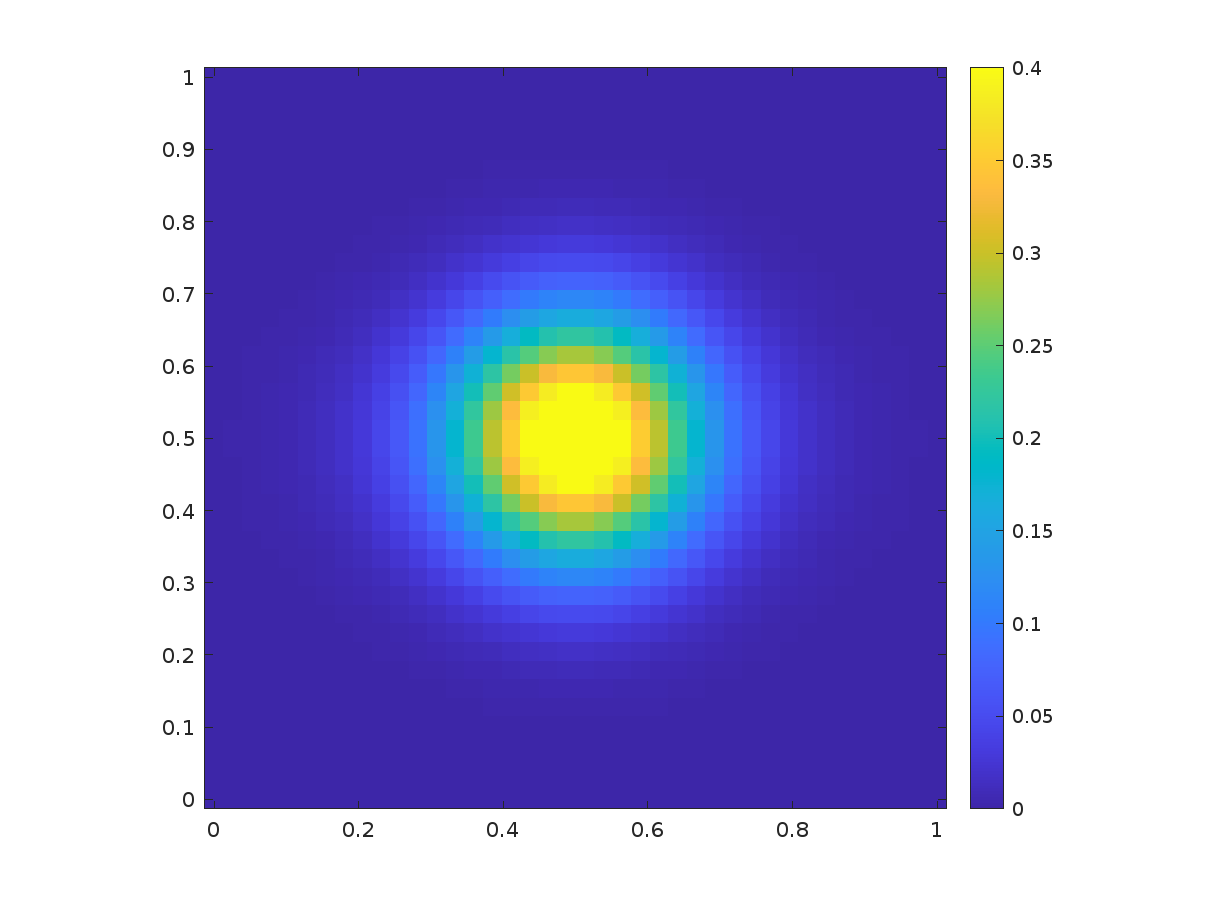}
  \quad
  \includegraphics[width=5.5cm]{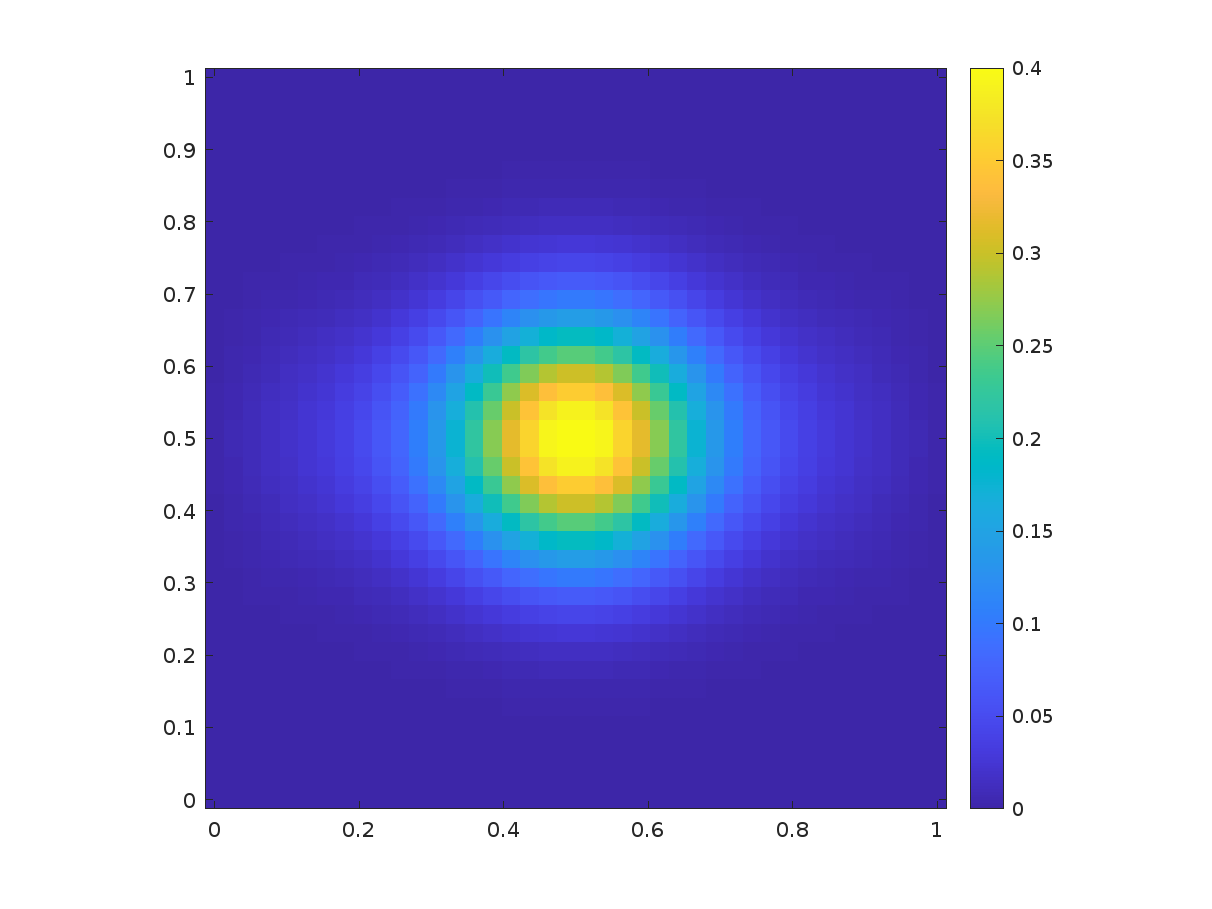}
  \quad
  \includegraphics[width=5.5cm]{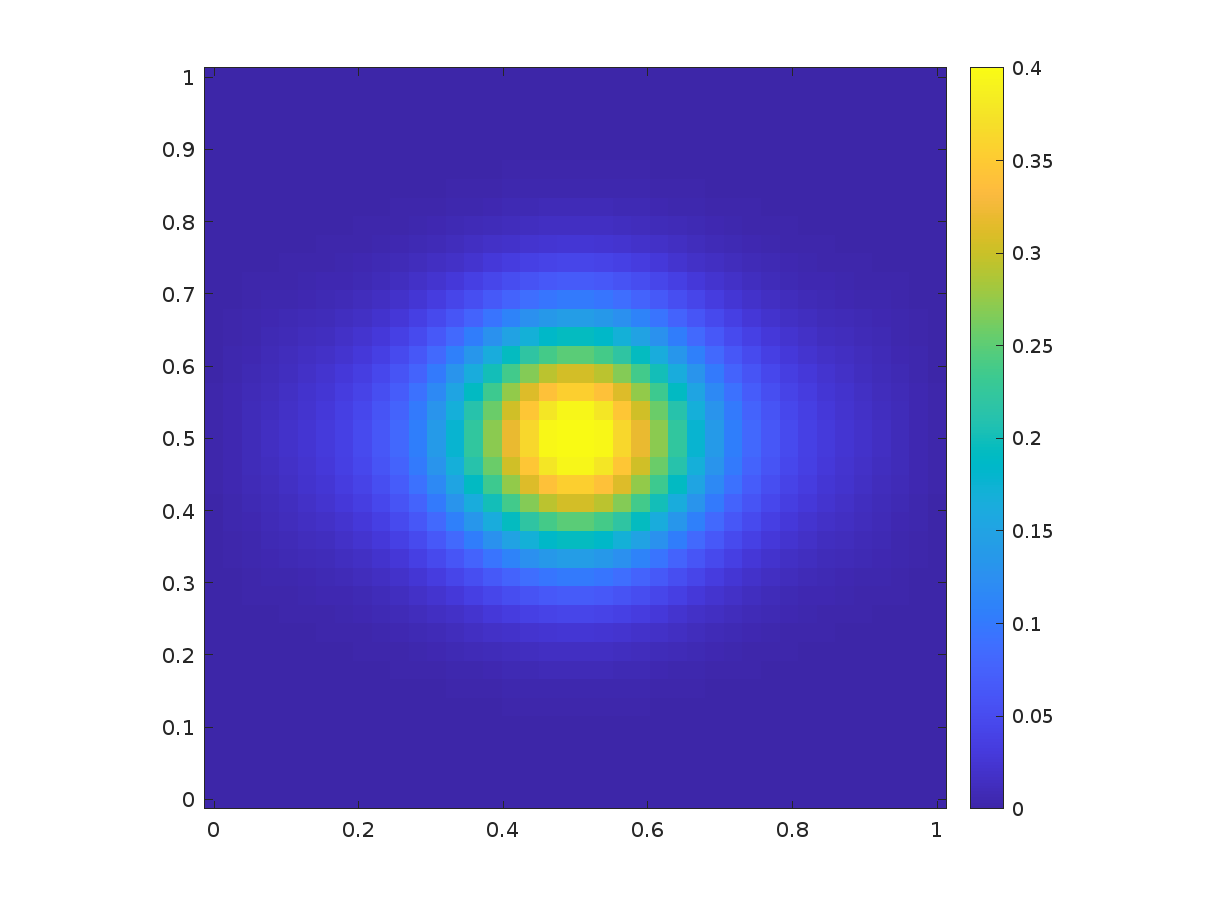}
  \quad
  \includegraphics[width=5.5cm]{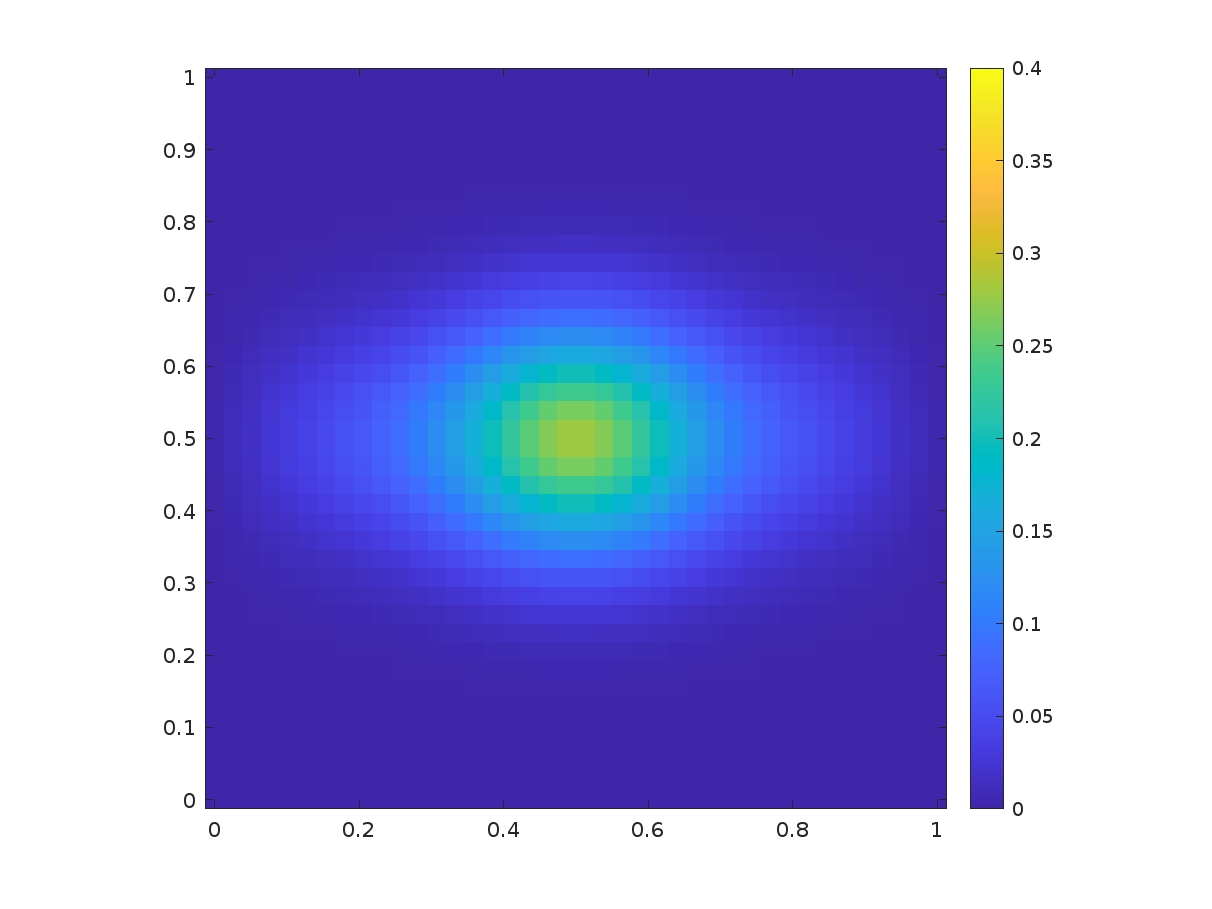}
  \quad
  \includegraphics[width=5.5cm]{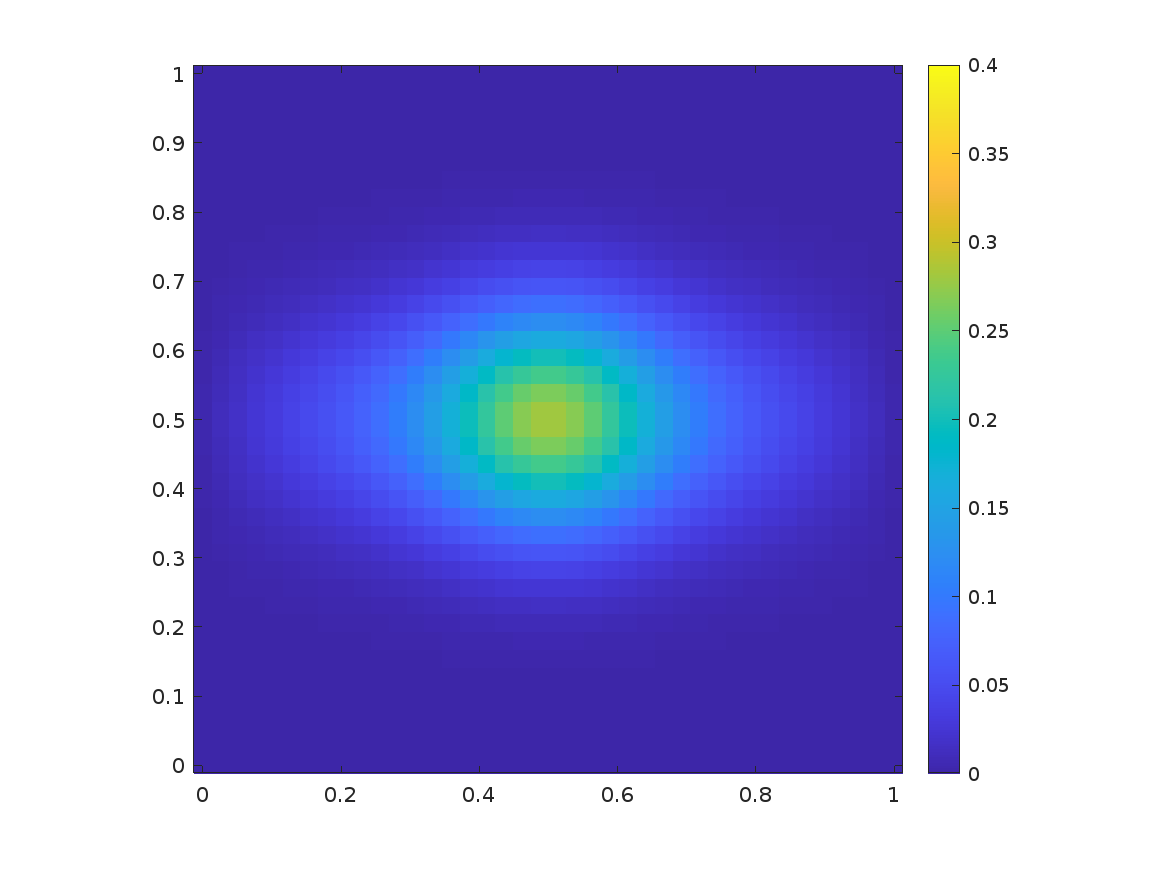}
  \quad
  \includegraphics[width=5.5cm]{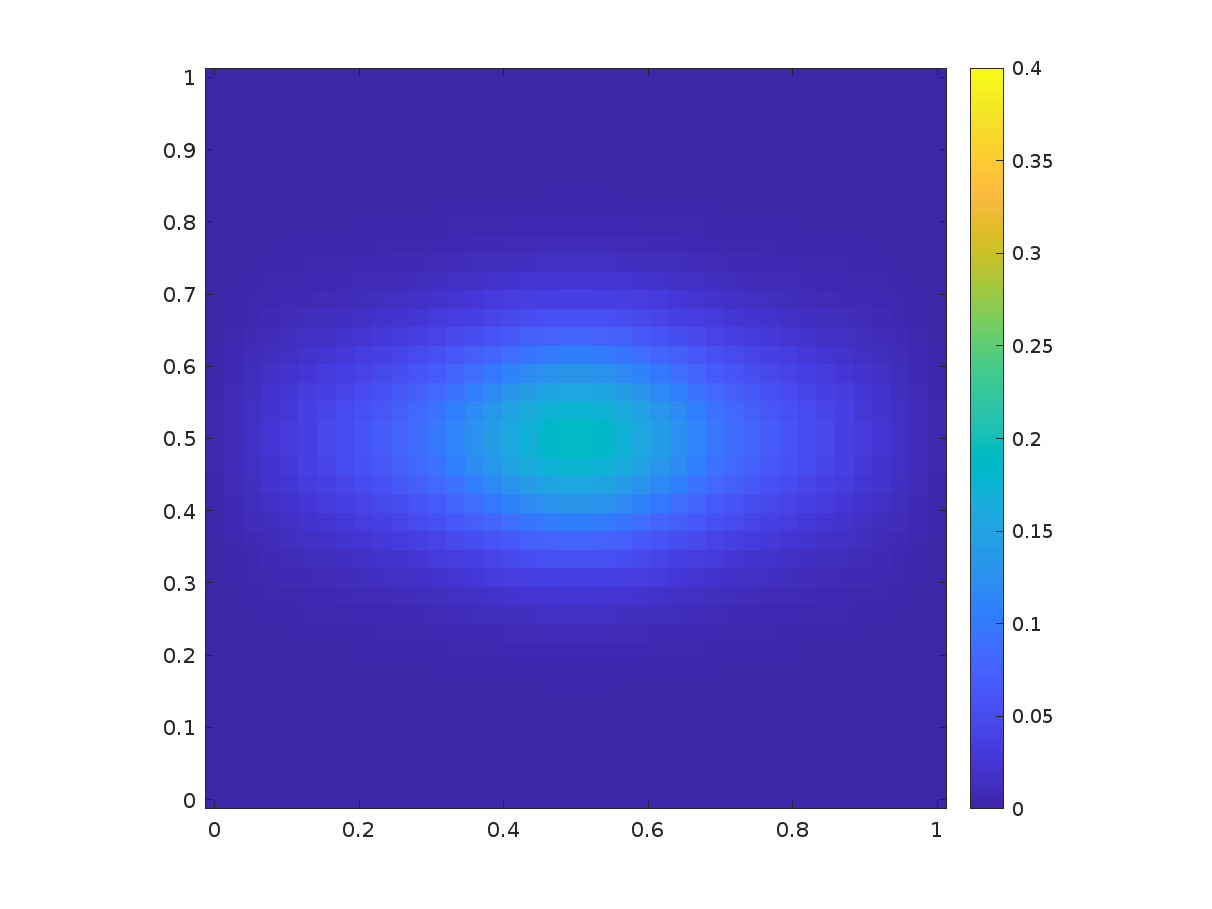}
  \quad
  \includegraphics[width=5.5cm]{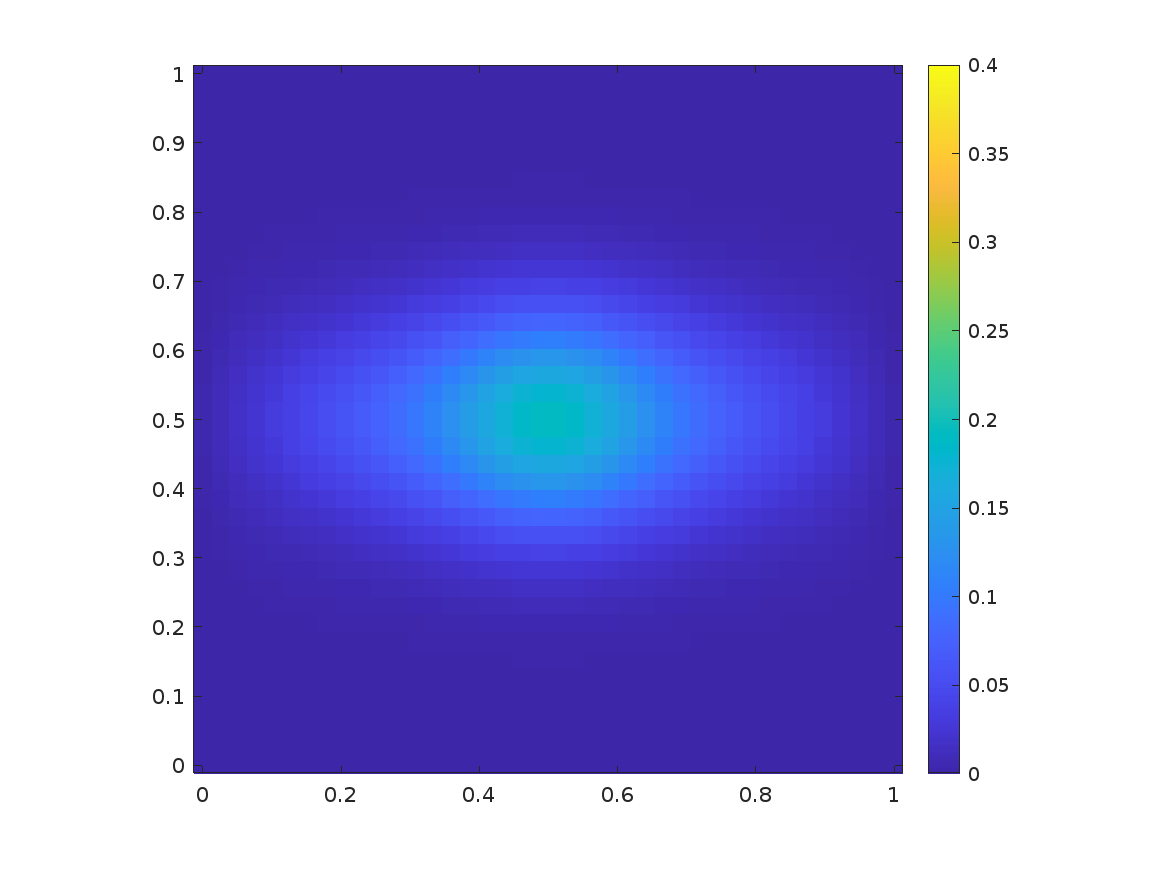}
  \quad
  \includegraphics[width=5.5cm]{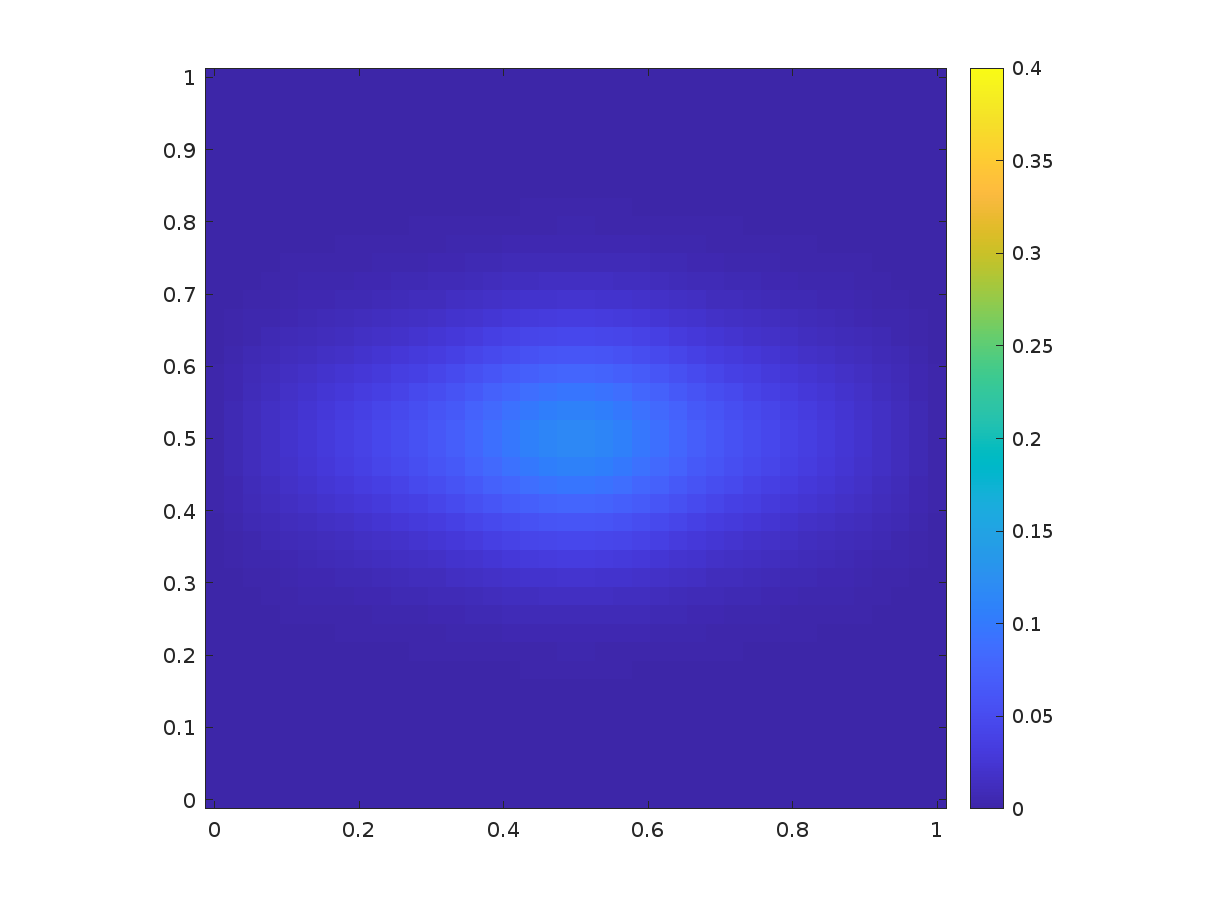}
  \quad
  \includegraphics[width=5.5cm]{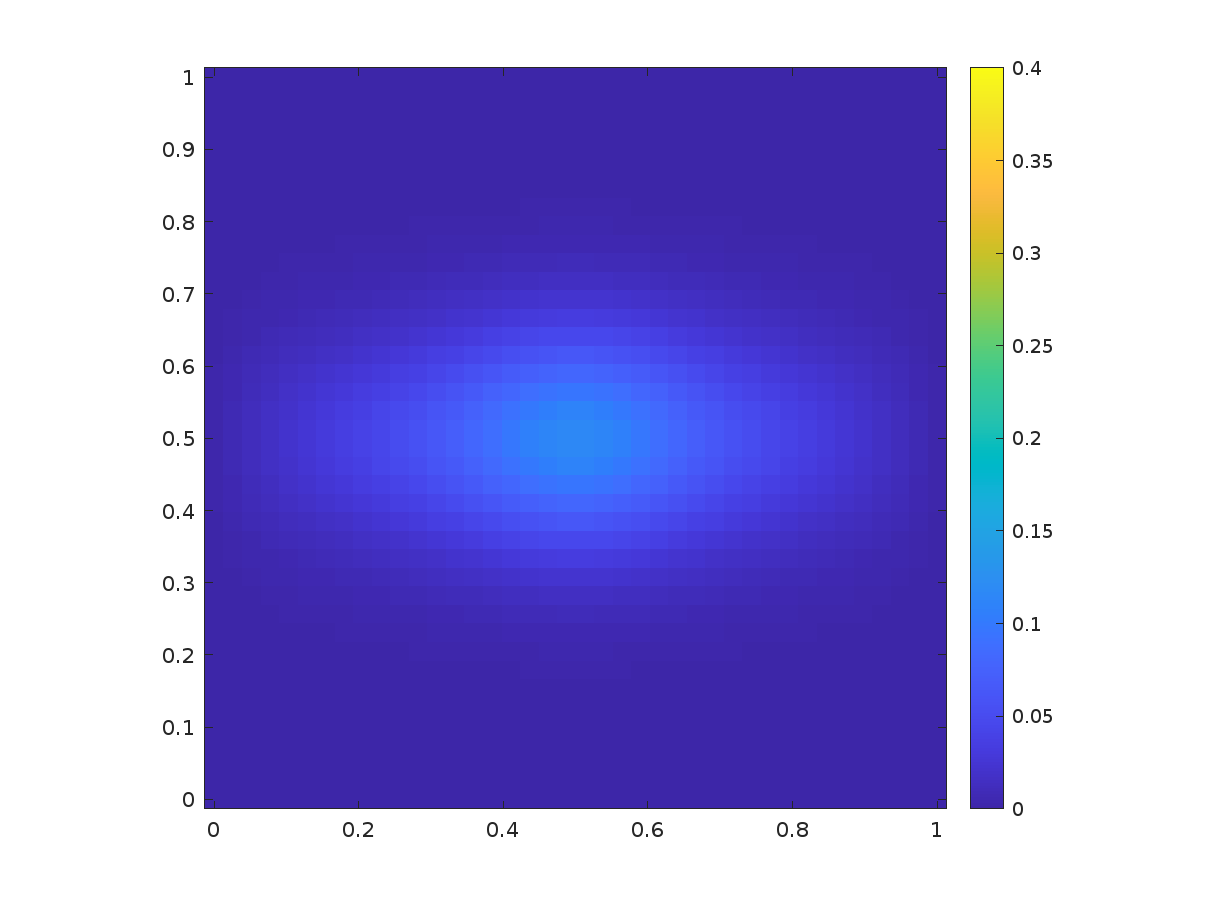}
  \caption{Solutions of concentration when $H=1/40$ for Case 1 in Example 1. First column: multiscale solution $C_1$ at $t=0.02$, $0.1$, $0.5$, $1$, $2$. Second column: reference averaged solution in $\Omega_1$ at the corresponding time instants.}
  \label{fig:Example1_Case1_U1}
\end{figure}

\begin{figure}
  \centering
  \includegraphics[width=5.5cm]{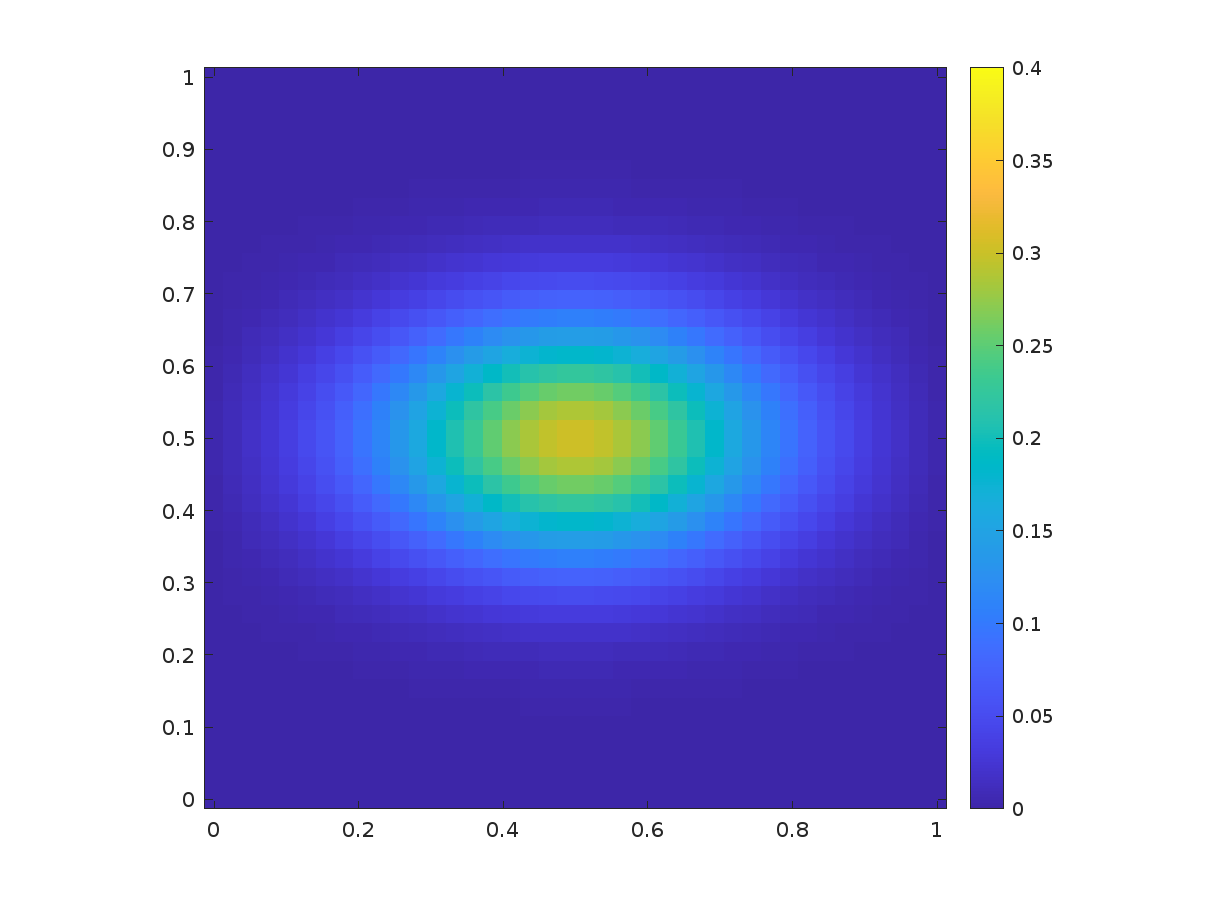}
  \quad
  \includegraphics[width=5.5cm]{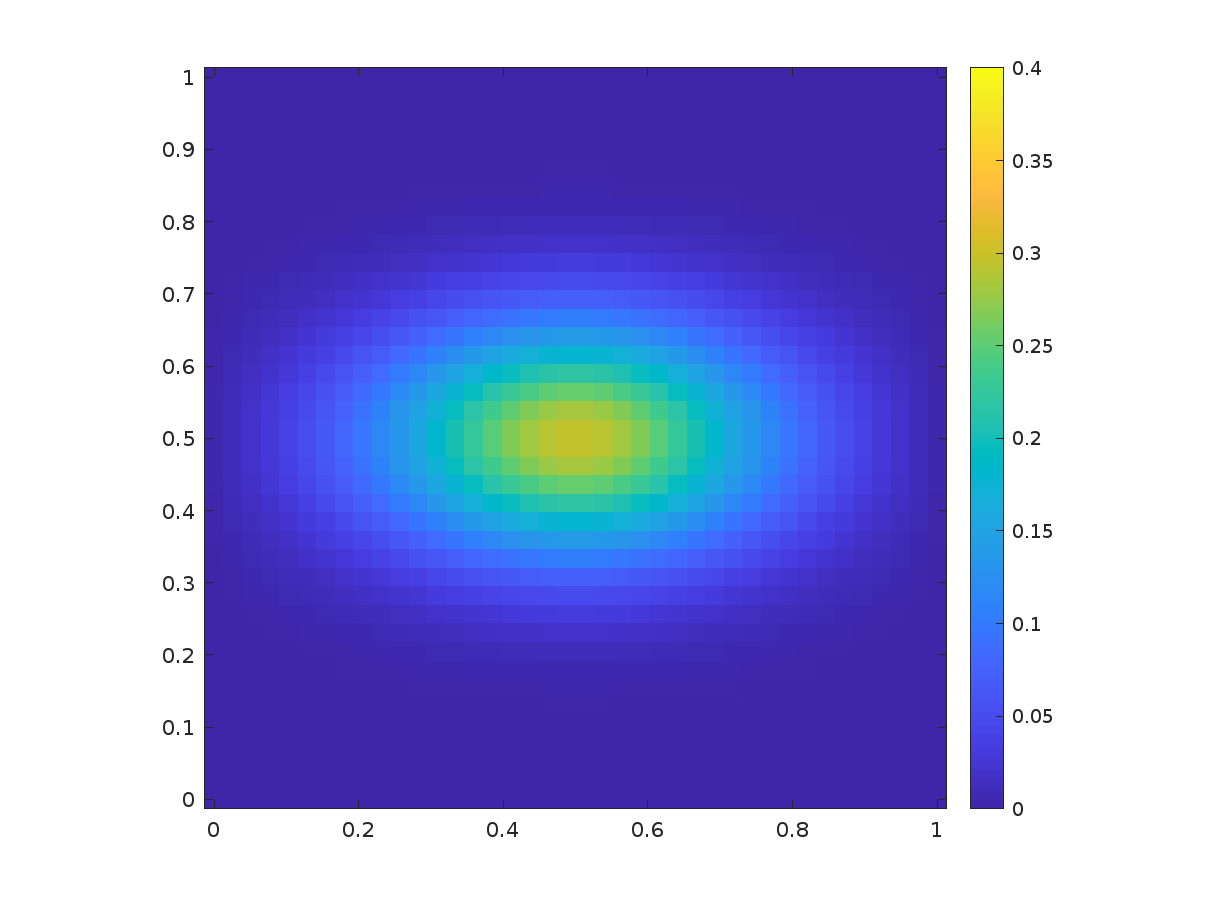}
  \quad
  \includegraphics[width=5.5cm]{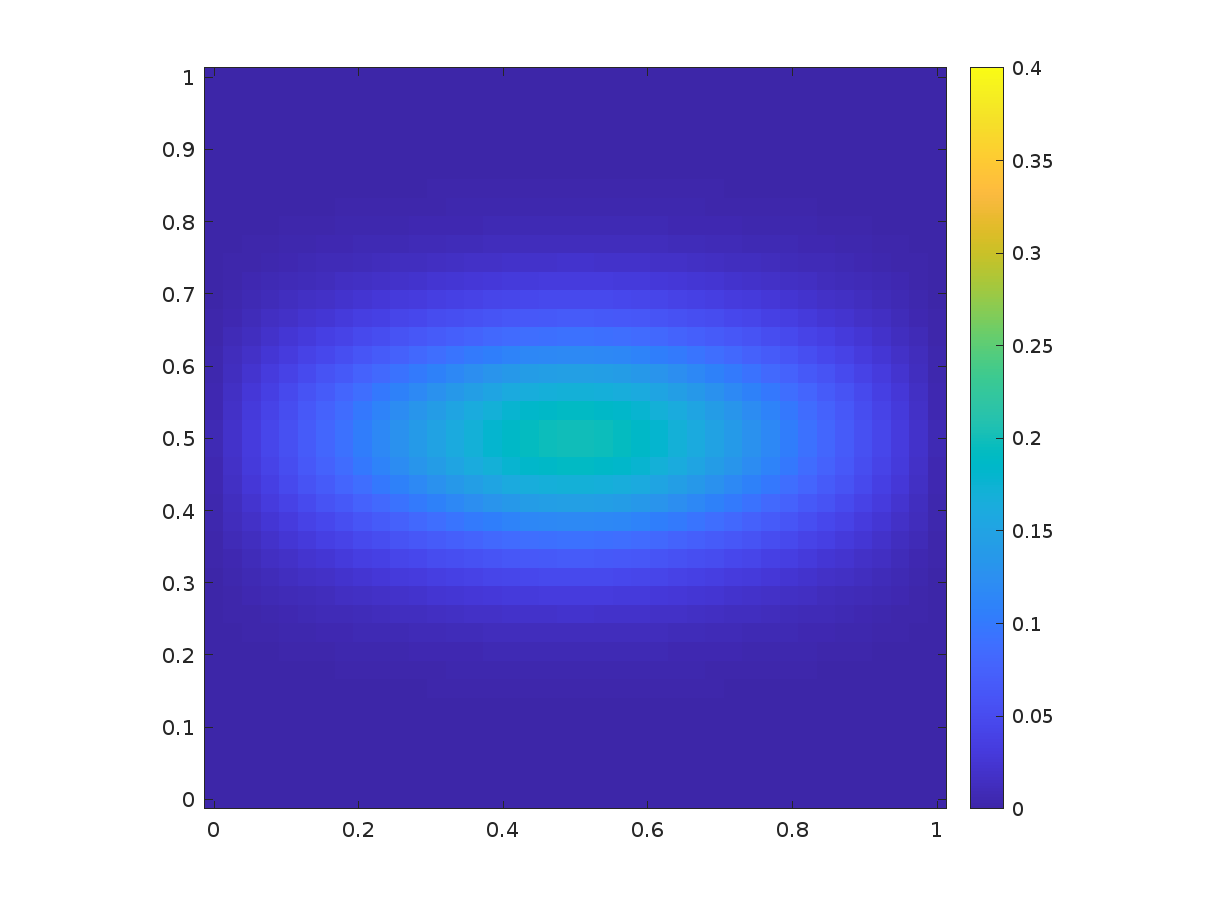}
  \quad
  \includegraphics[width=5.5cm]{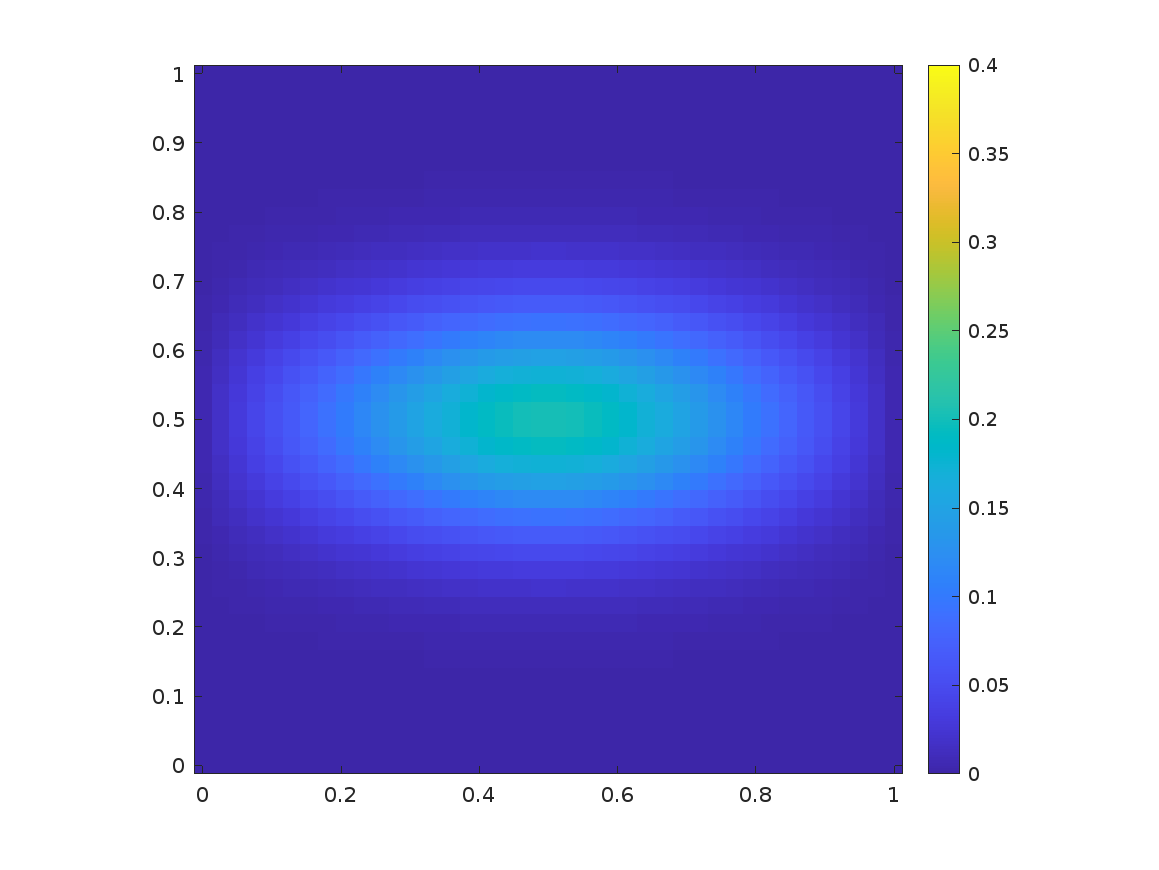}
  \quad
  \includegraphics[width=5.5cm]{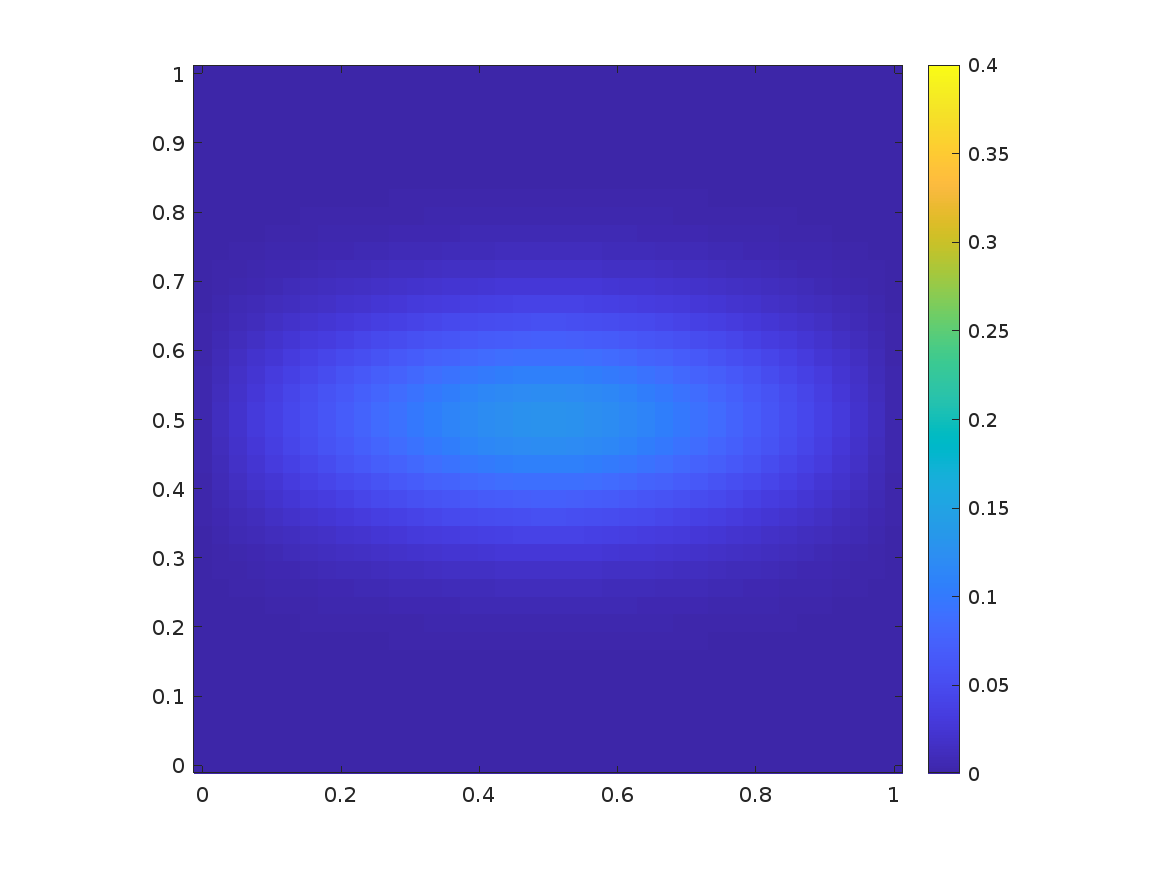}
  \quad
  \includegraphics[width=5.5cm]{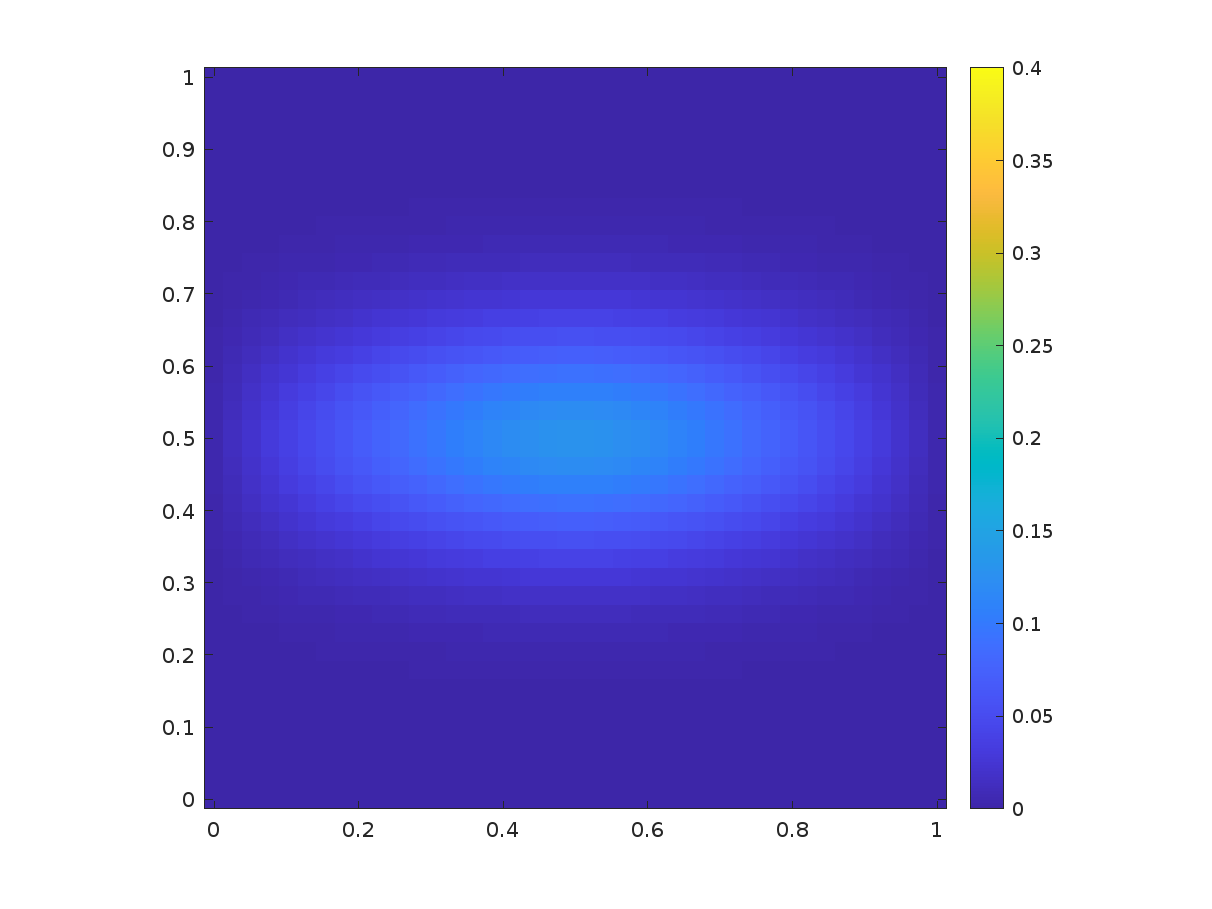}
  \quad
  \includegraphics[width=5.5cm]{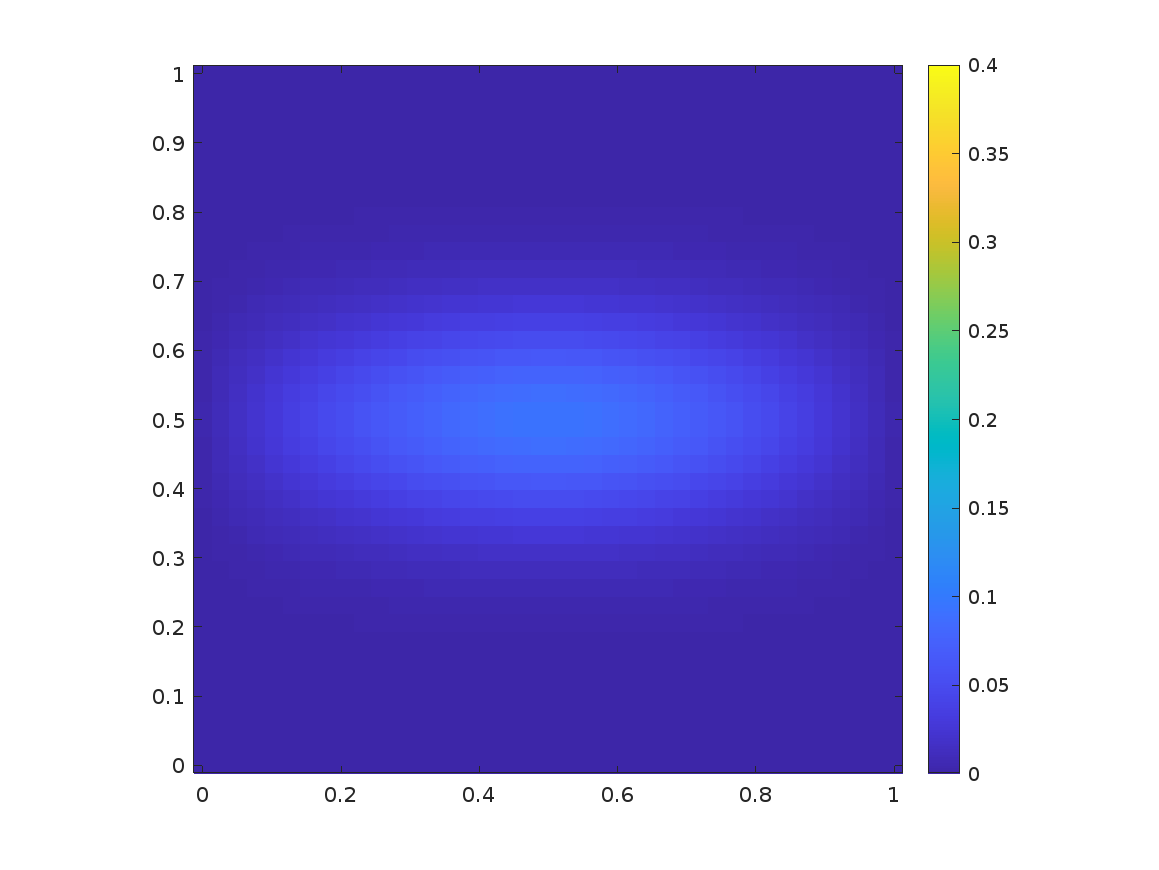}
  \quad
  \includegraphics[width=5.5cm]{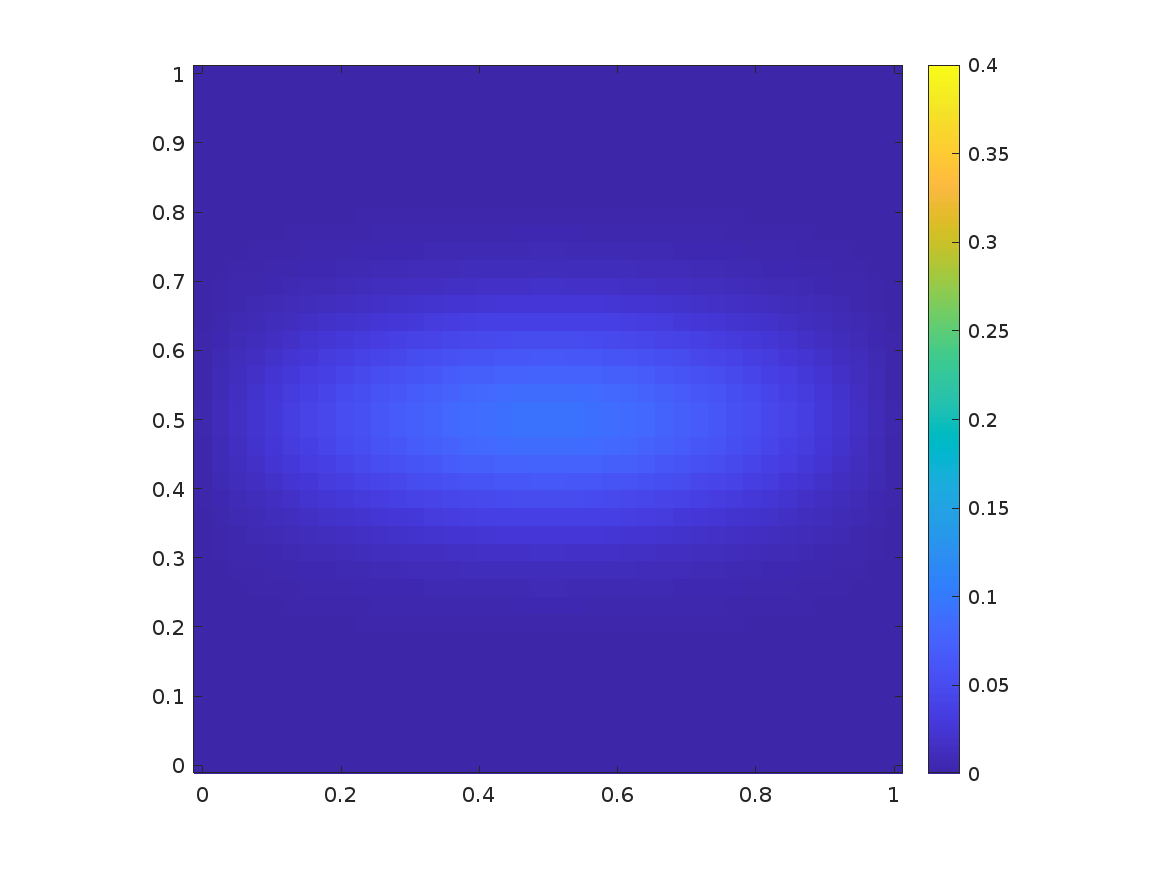}
  \quad
  \includegraphics[width=5.5cm]{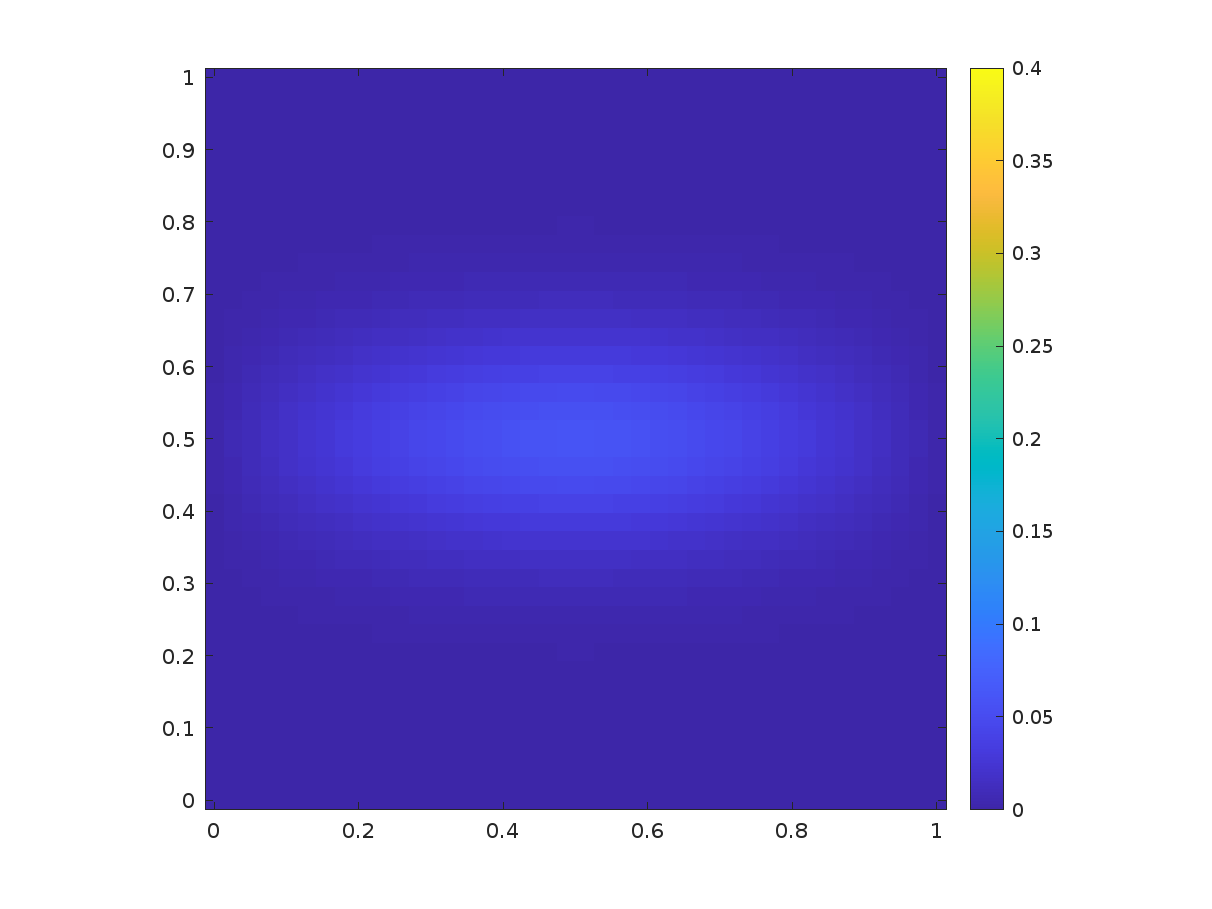}
  \quad
  \includegraphics[width=5.5cm]{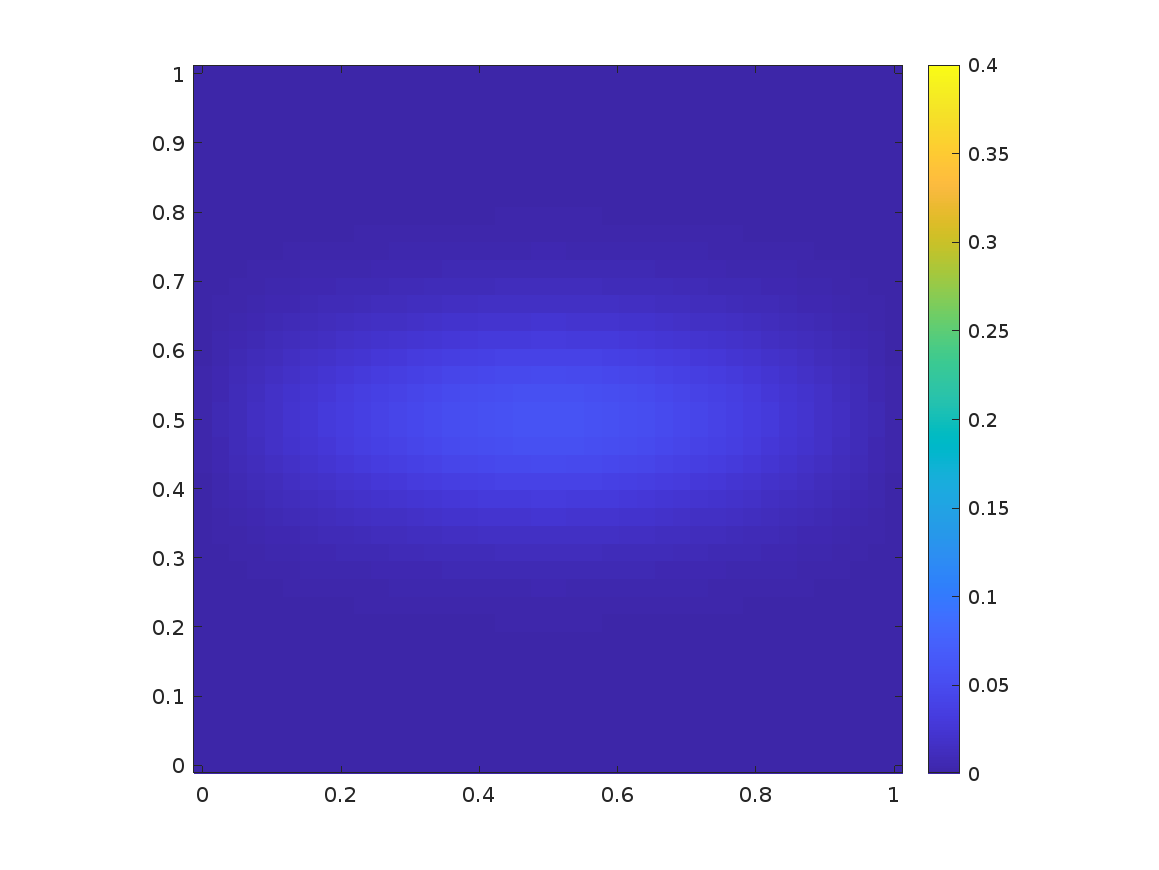}
  \caption{Solutions of concentration when $H=1/40$ for Case 1 in Example 1. First column: multiscale solution $C_2$ at $t=0.02$, $0.1$, $0.5$, $1$, $2$. Second column: reference averaged solution in $\Omega_2$ at the corresponding time instants.}
  \label{fig:Example1_Case1_U2}
\end{figure}

We depict the concentration at the initial time $t=0$ in Figure \ref{fig:Example1_init}. 
The relative $L^2$ errors of our proposed algorithm are shown in Table \ref{tab:Example1_Case1_error}. We also depict the results when $H=1/40$ in Figure \ref{fig:Example1_Case1_U1} and Figure \ref{fig:Example1_Case1_U2}. Figure \ref{fig:Example1_Case1_U1} demonstrates the multiscale solution $C_1$ and the reference averaged solution in $\Omega_1$ at $t=0.02$, $0.1$, $0.5$, $1$, $2$; Figure \ref{fig:Example1_Case1_U2} demonstrates the multiscale solution $C_2$ and its corresponding reference averaged solution. 

In this paragraph, we examine the obtained results. First, we can clearly see that the errors at different time $t$ between the multiscale and reference solutions are very small and our method can approximate the results accurately. Furthermore, decreasing the coarse mesh size $H$ makes the results more precise. 
Also, one can observe that the averaged concentration in $\Omega_1$ transports more slowly than that in $\Omega_2$. This is reasonable because the diffusion coefficient is much higher in $\Omega_2$ than in $\Omega_1$, as defined in (\ref{eq:Example1_coefficient}).

\begin{table}
\caption{Relative $L^2$ errors at $t=0.02$, $0.1$, $0.5$, $1$, $2$ for Case 1 in Example 1. Left: $H=1/20$ and $l=6$. Right: $H=1/40$ and $l=8$.}
\centering
\begin{tabular}{ccc}
   \toprule
   $t$      & $e^{(1)}(t)$    &   $e^{(2)}(t)$ \\
   \midrule
   $0.02$ & $4.63\%$ & $2.13\%$ \\
   $0.1$ & $3.03\%$ & $2.89\%$ \\
   $0.5$ & $2.54\%$ & $1.90\%$ \\
   $1.0$ & $2.39\%$ & $1.53\%$ \\
   $2.0$ & $1.77\%$ & $1.20\%$ \\
   \bottomrule 
\end{tabular}
\qquad
\begin{tabular}{ccc}
   \toprule
   $t$      & $e^{(1)}(t)$    &   $e^{(2)}(t)$ \\
   \midrule
   $0.02$ & $2.43\%$ & $2.20\%$ \\
   $0.1$ & $0.76\%$ & $1.32\%$ \\
   $0.5$ & $0.39\%$ & $0.22\%$ \\
   $1.0$ & $0.55\%$ & $0.51\%$ \\
   $2.0$ & $0.60\%$ & $0.82\%$ \\
   \bottomrule 
\end{tabular}
\label{tab:Example1_Case1_error}
\end{table}

\subsubsection{Case 2}

For Case 2, we change the boundary condition for the pressure equation to be an inhomogeneous Dirichlet boundary condition while keeping the homogeneous Dirichlet condition for the concentration equation, that is,
\begin{equation}
\label{eq:Example1_Case2_bc}
    \begin{split}
        p(x,t) = x, \quad & x\in \partial \Omega, \quad \forall t\in [0,\infty),\\
        c(x,t) = 0, \quad & x\in \partial \Omega, \quad \forall t\in [0,\infty).
    \end{split}
\end{equation}
All the other assumptions are the same as in Case 1.

\begin{figure}
    \centering
    \includegraphics[width=5.5cm]{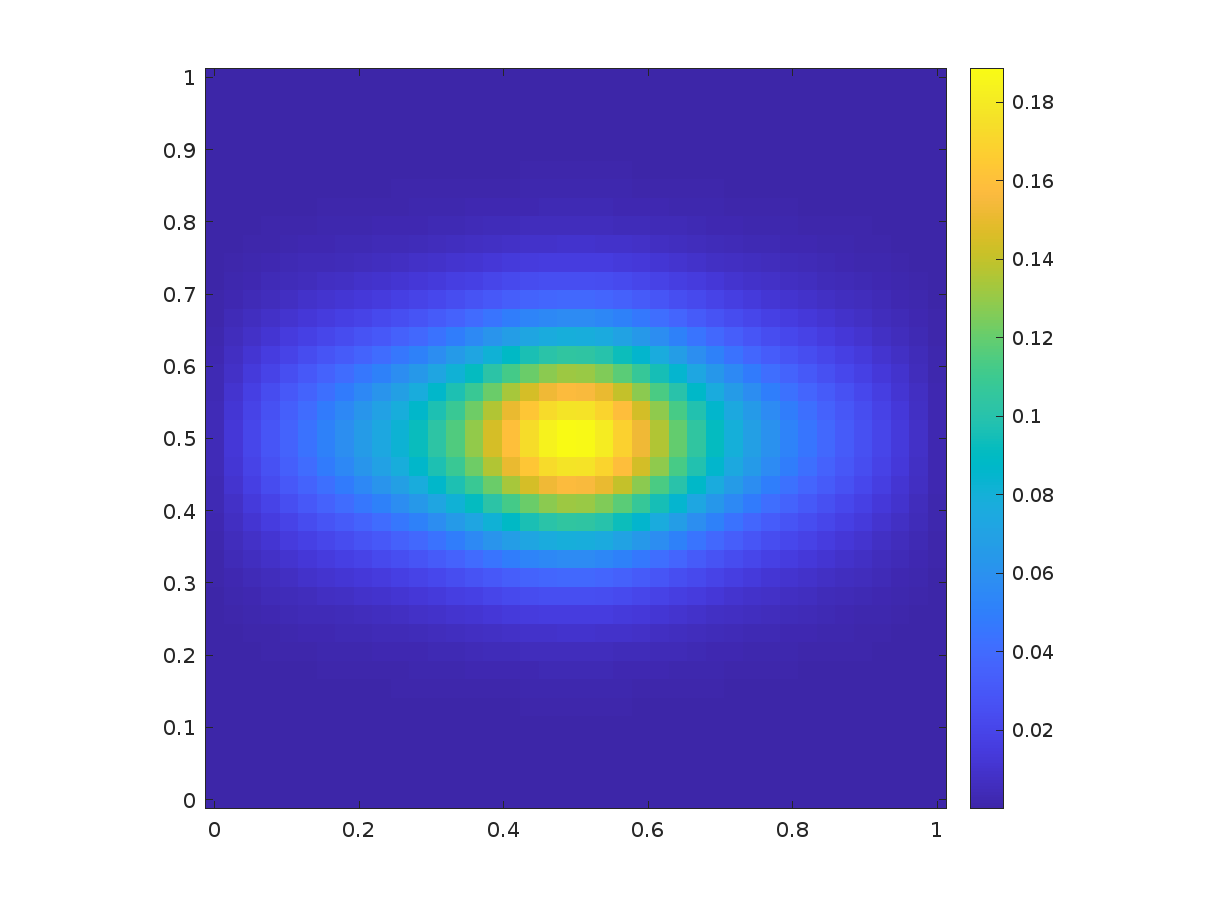}
    \quad
    \includegraphics[width=5.5cm]{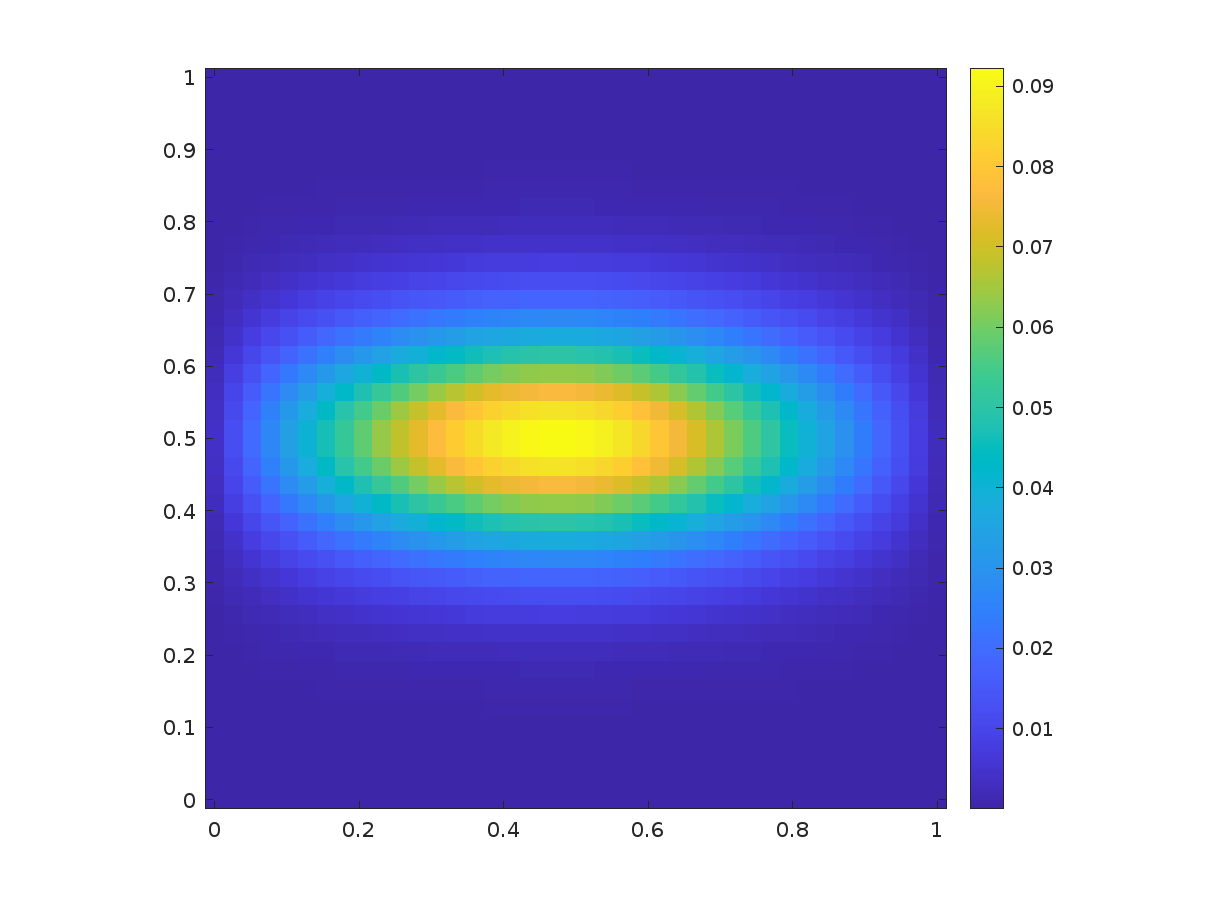}
    \caption{Snapshots of $C_1$ and $C_2$ when $t=1$ for Case 2 in Example}
    \label{fig:Example1_Case2_snapshot}
\end{figure}

\begin{figure}
  \centering
  \includegraphics[width=5.5cm]{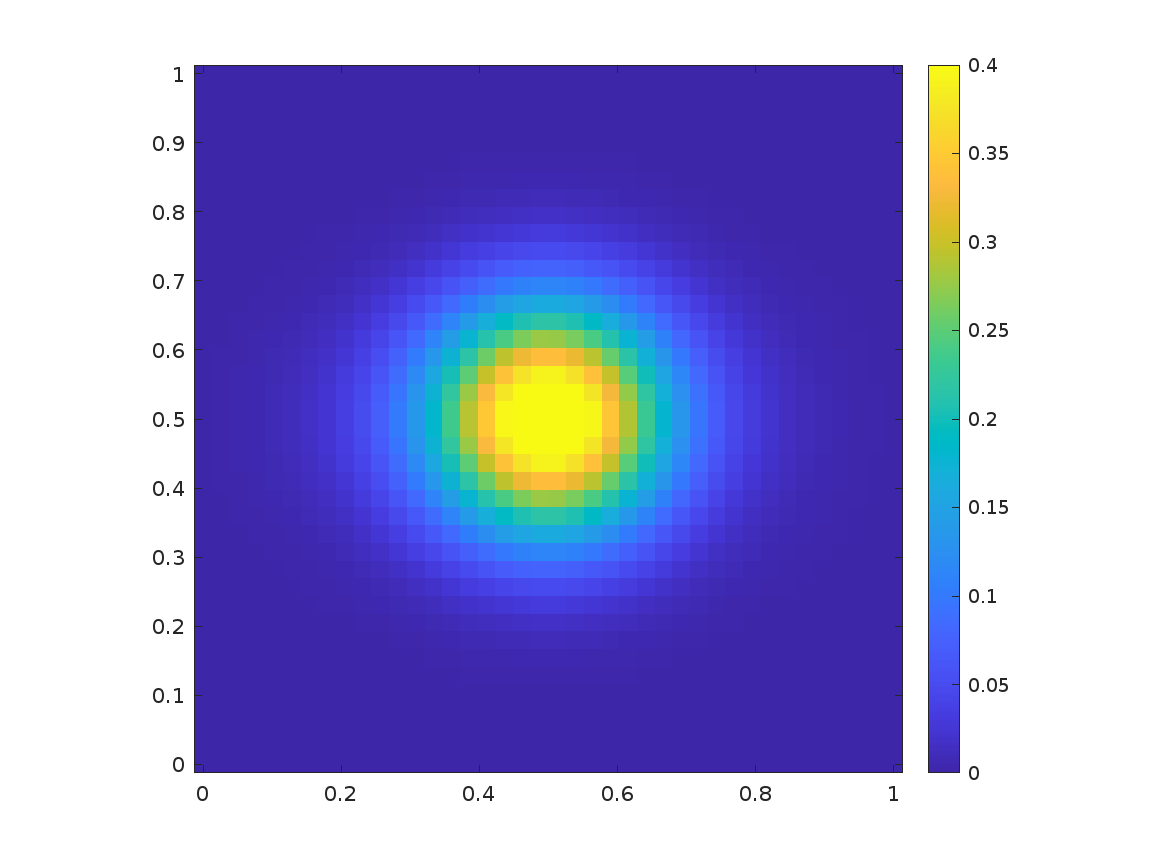}
  \quad
  \includegraphics[width=5.5cm]{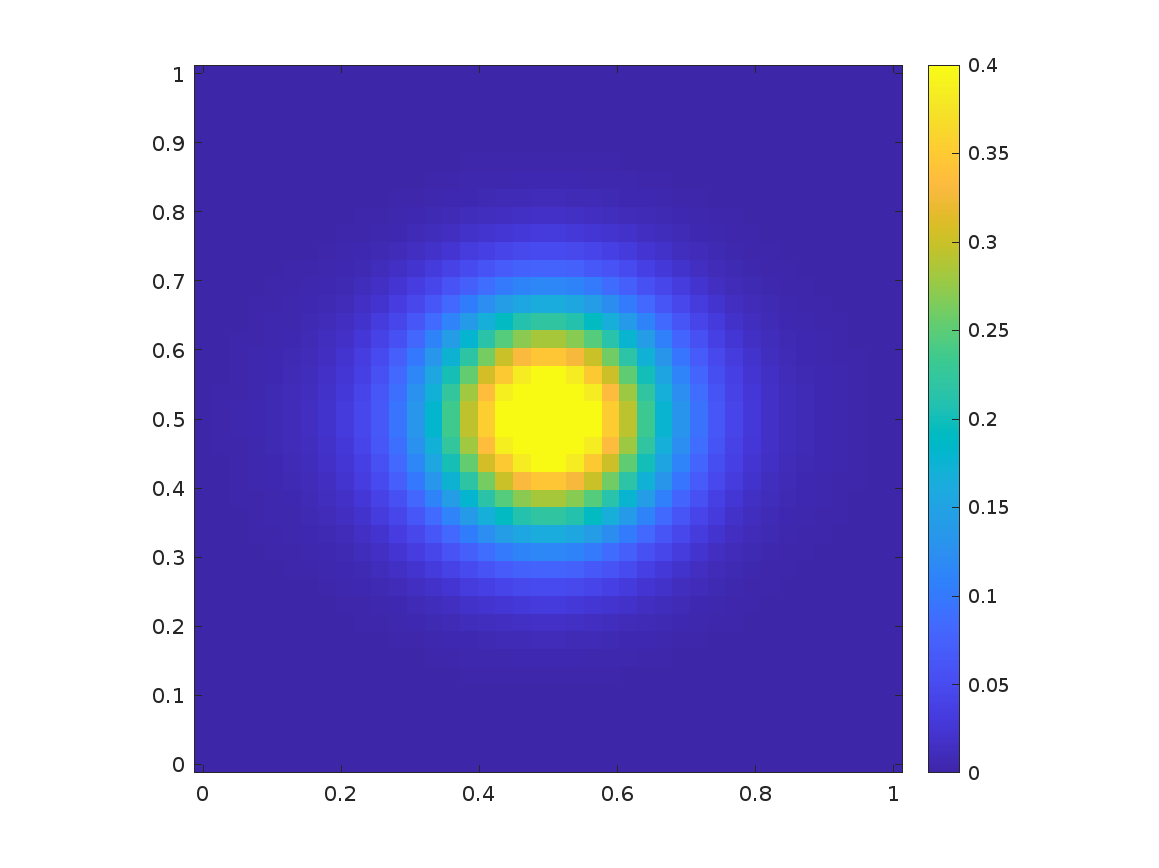}
  \quad
  \includegraphics[width=5.5cm]{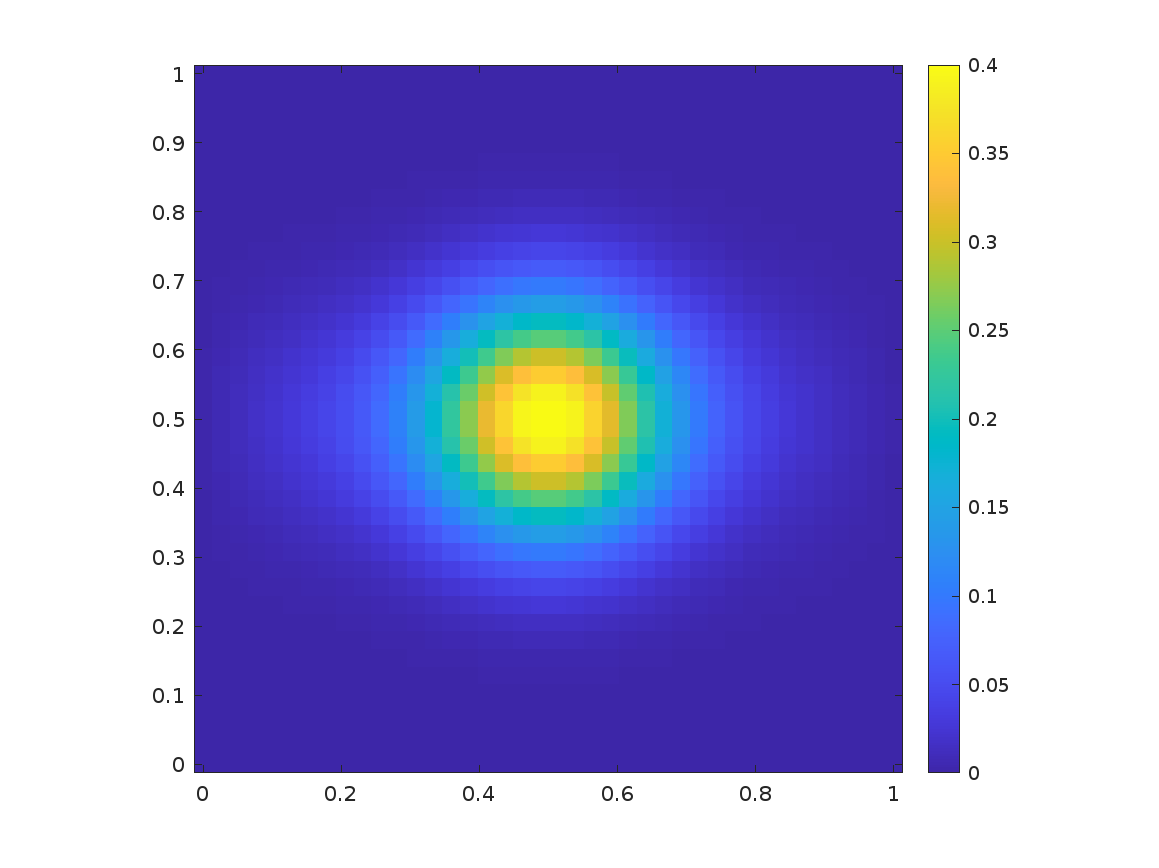}
  \quad
  \includegraphics[width=5.5cm]{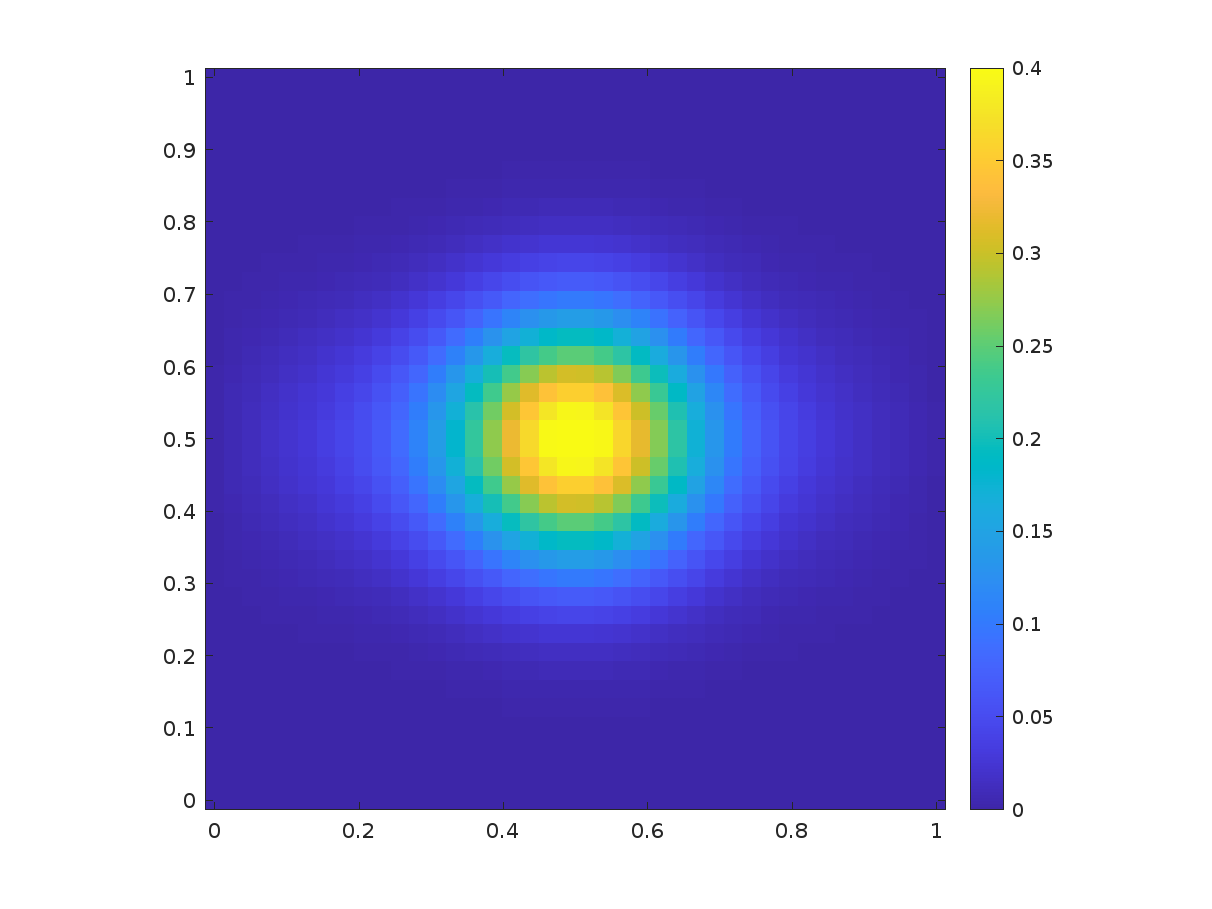}
  \quad
  \includegraphics[width=5.5cm]{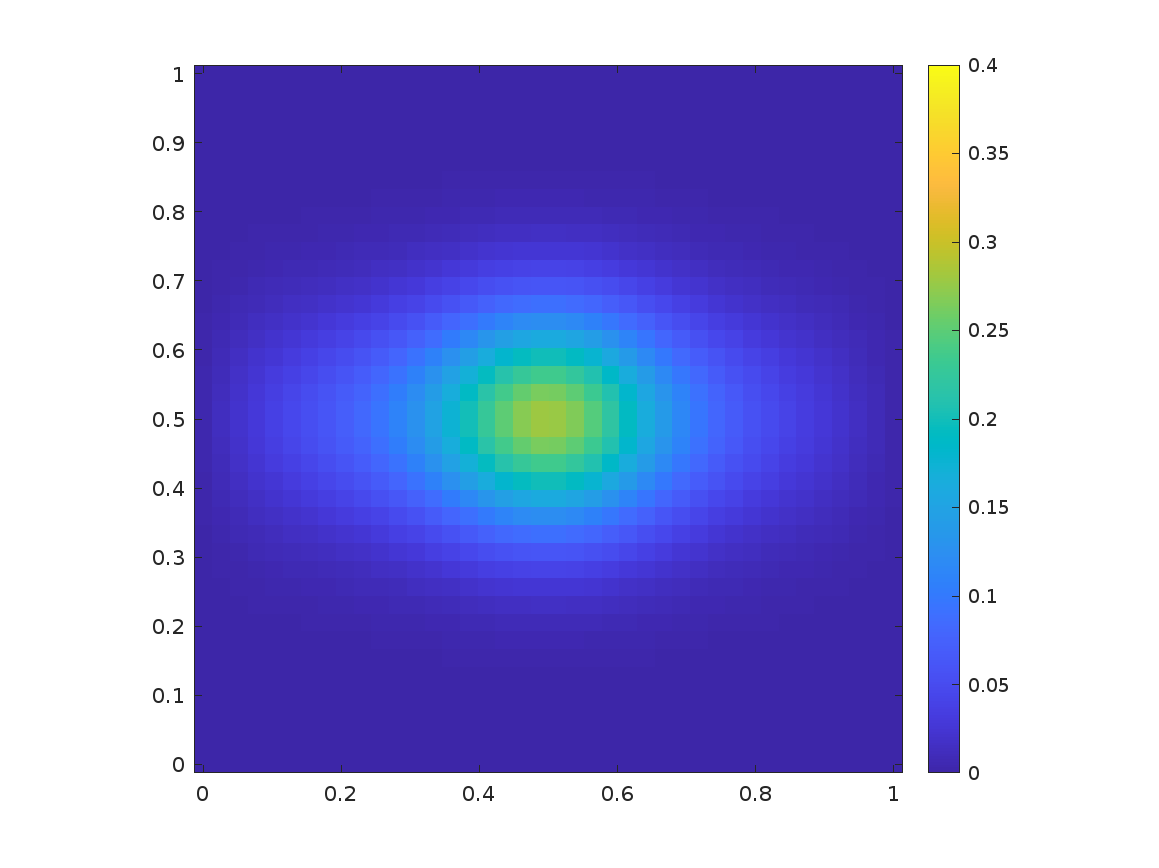}
  \quad
  \includegraphics[width=5.5cm]{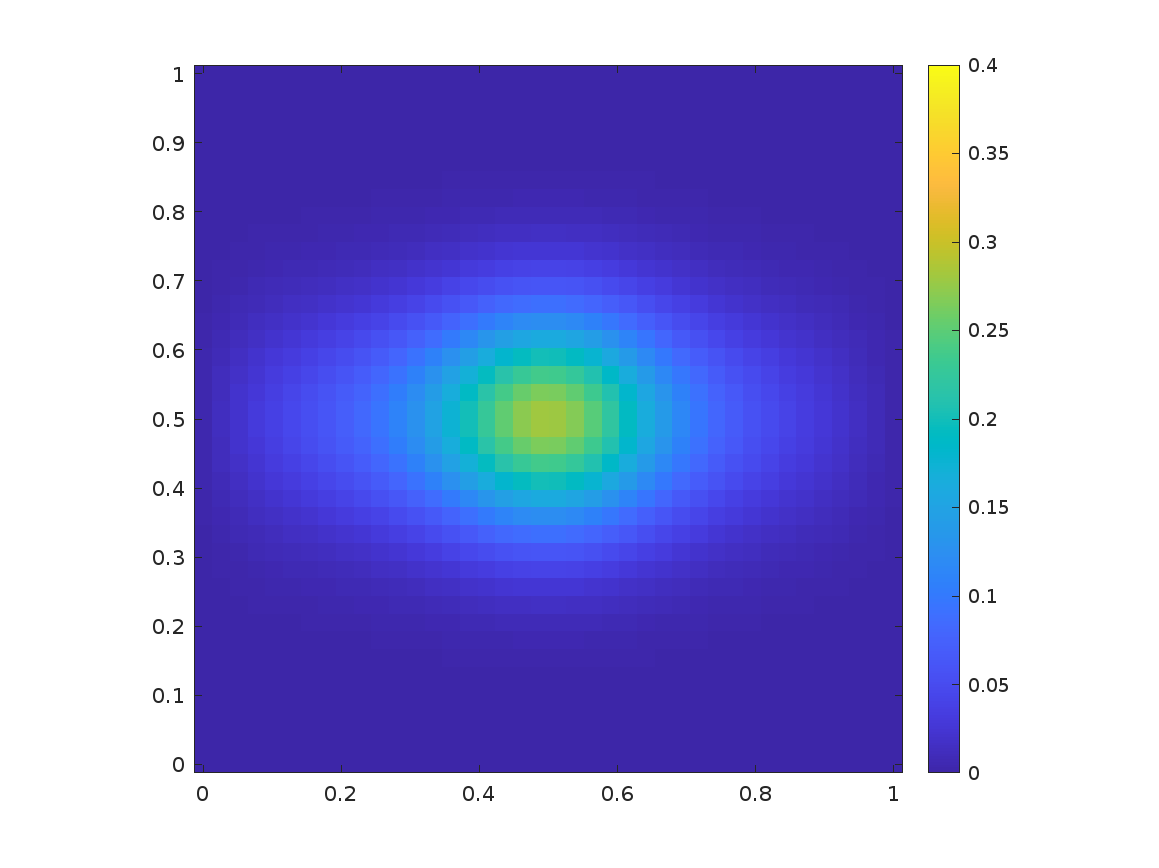}
  \quad
  \includegraphics[width=5.5cm]{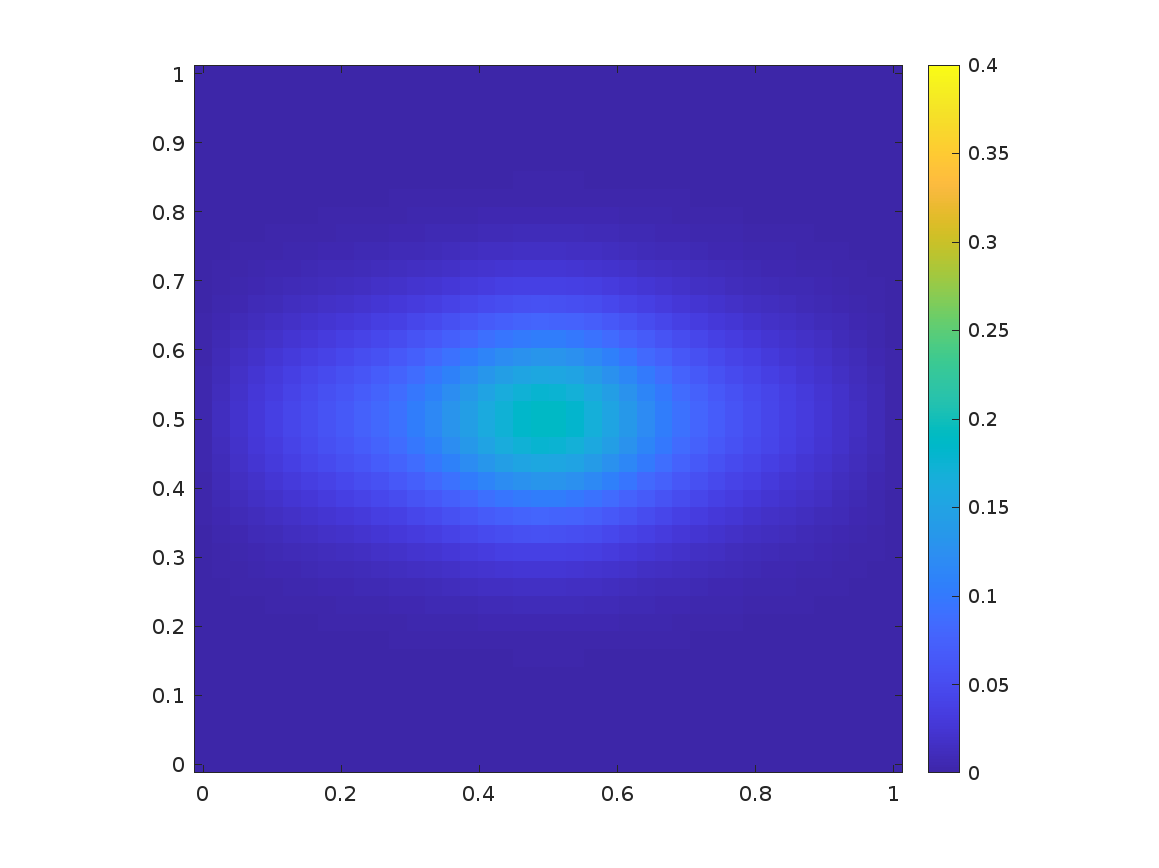}
  \quad
  \includegraphics[width=5.5cm]{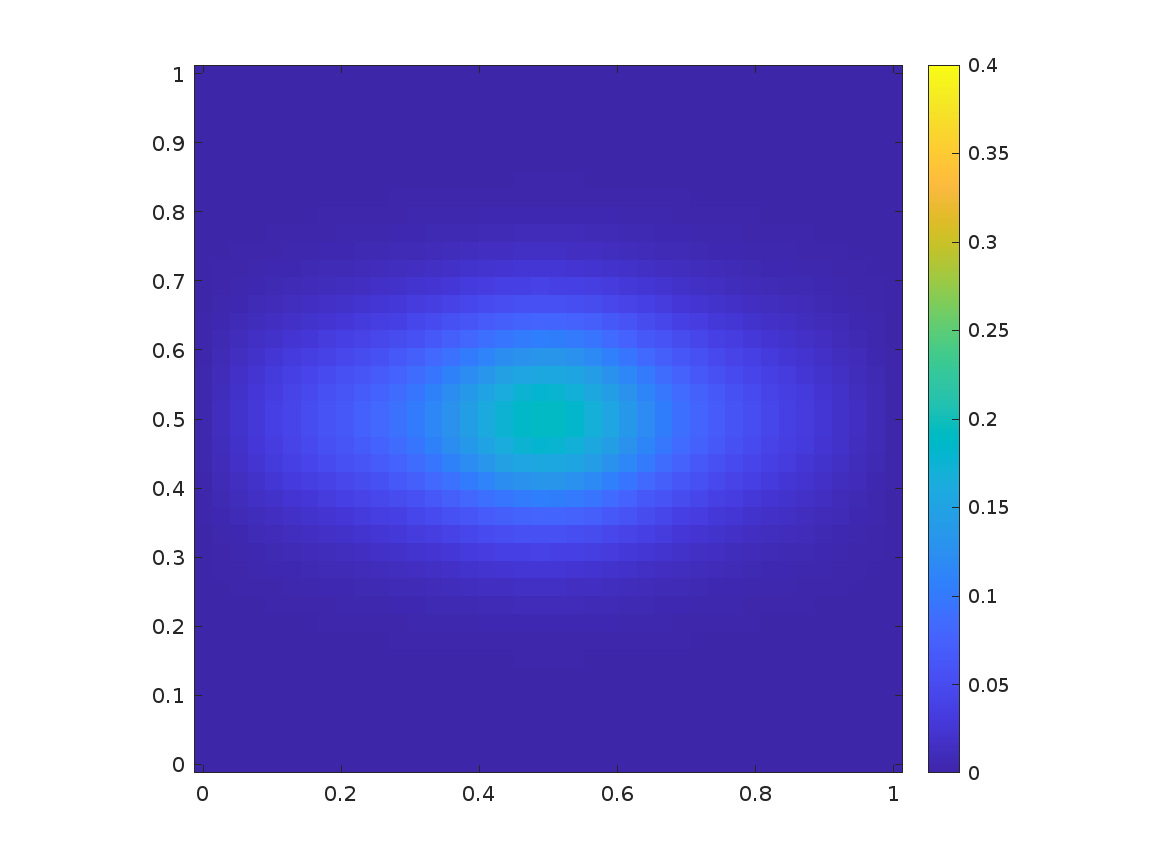}
  \quad
  \includegraphics[width=5.5cm]{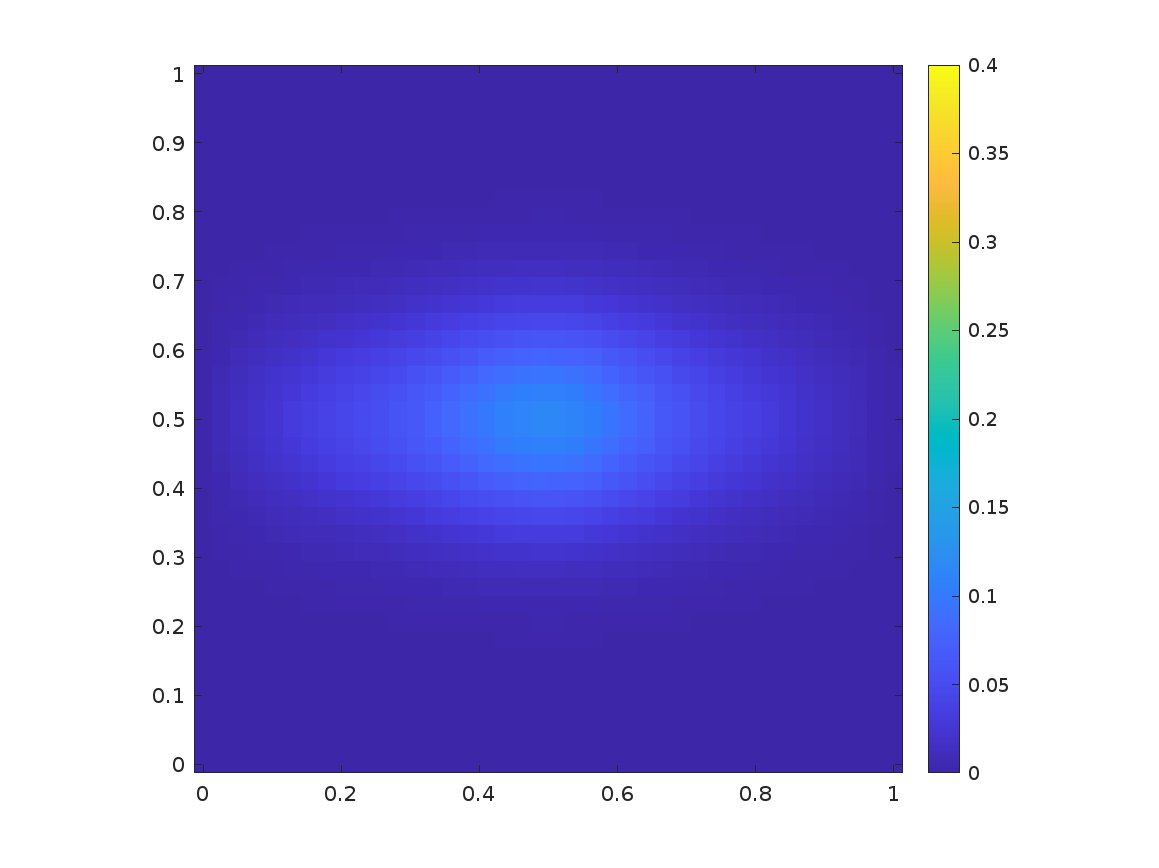}
  \quad
  \includegraphics[width=5.5cm]{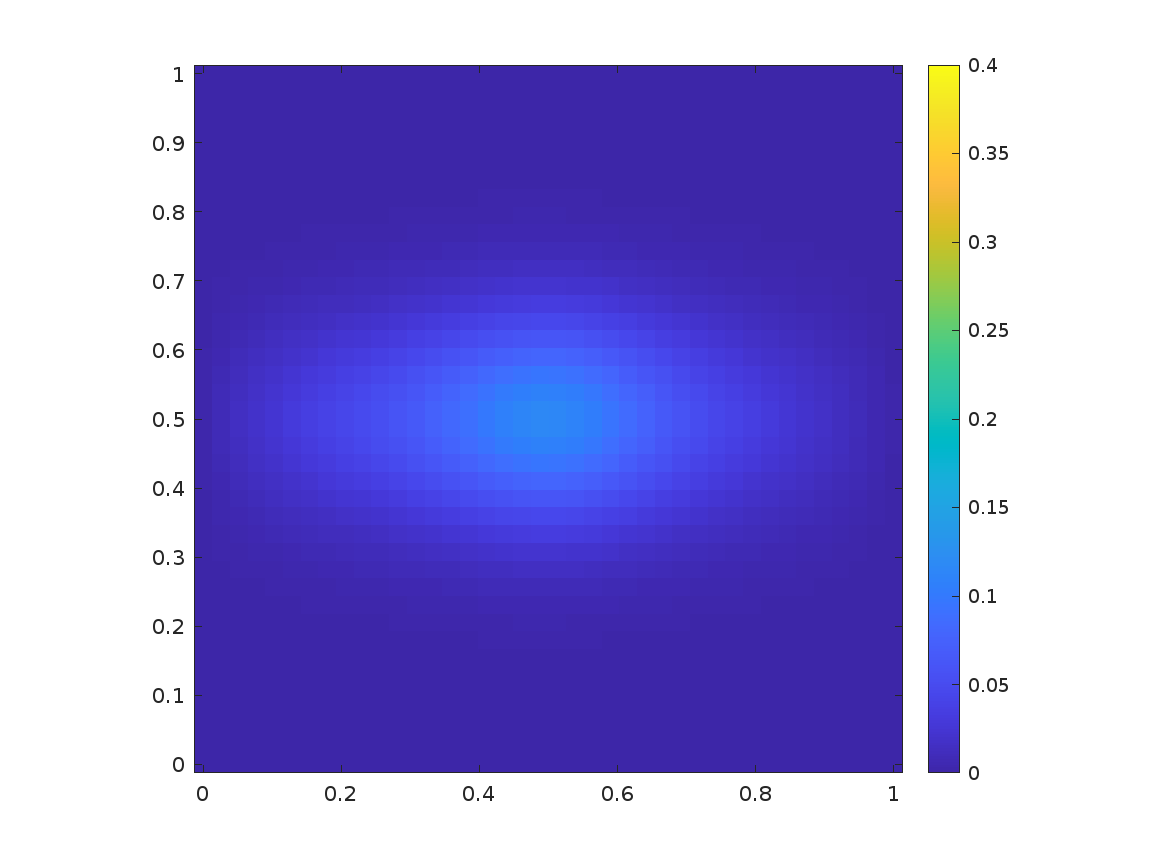}
  \caption{Solutions of concentration when $H=1/40$ for Case 2 in Example 1. First column: multiscale solution $C_1$ at $t=0.02$, $0.1$, $0.5$, $1$, $2$. Second column: reference averaged solution in $\Omega_1$ at the corresponding time instants.}
  \label{fig:Example1_Case2_U1}
\end{figure}

\begin{figure}
  \centering
  \includegraphics[width=5.5cm]{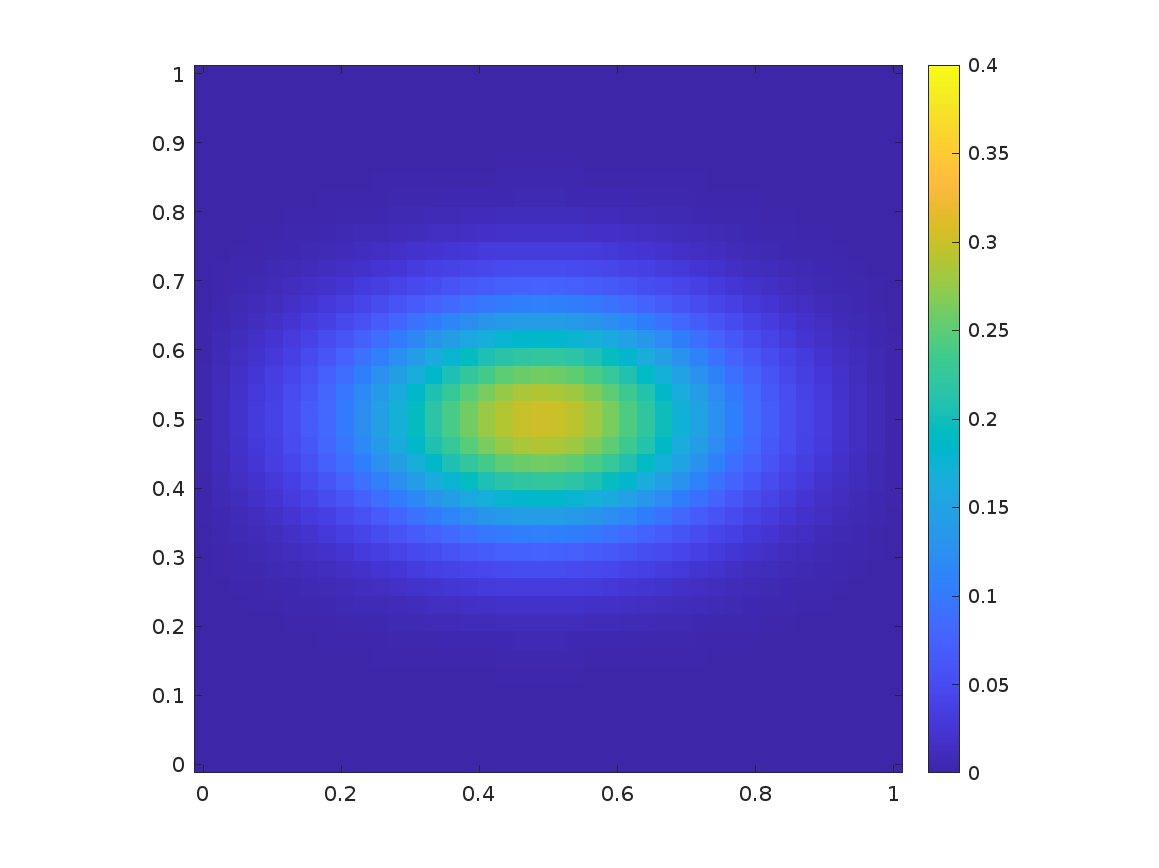}
  \quad
  \includegraphics[width=5.5cm]{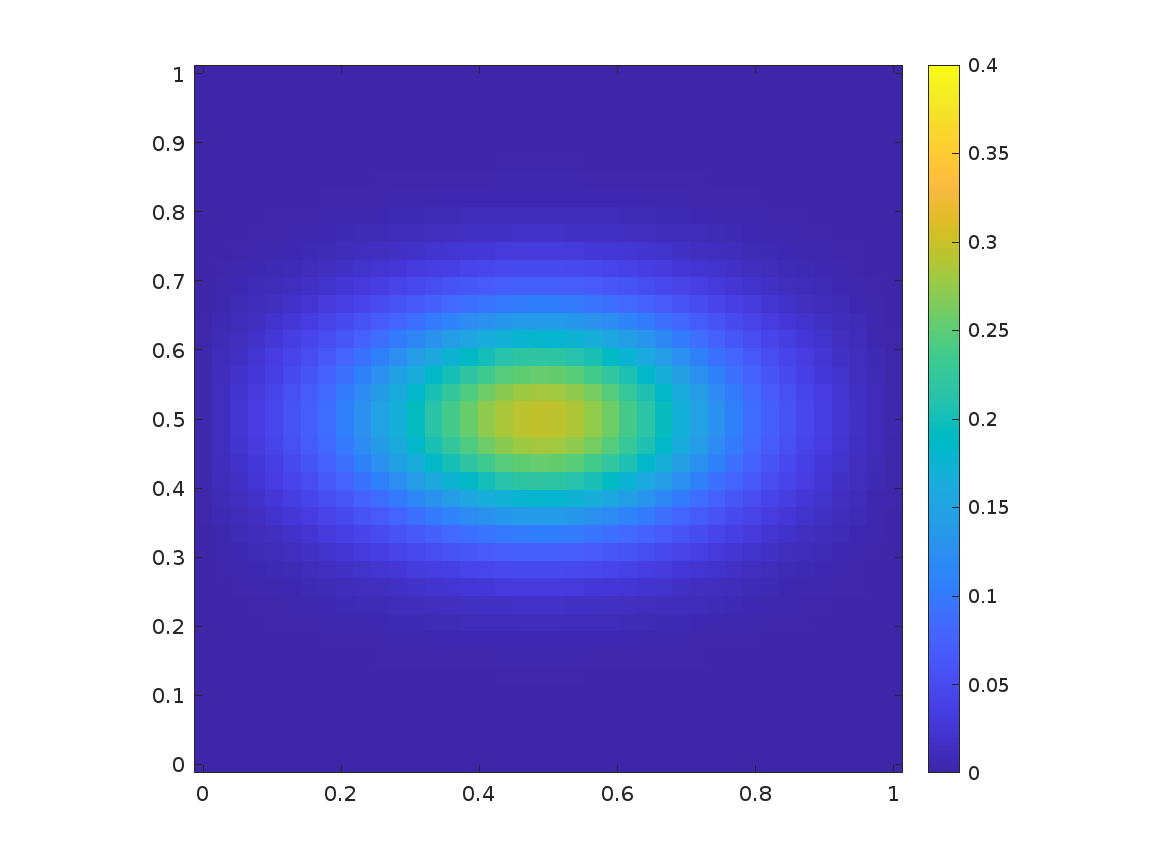}
  \quad
  \includegraphics[width=5.5cm]{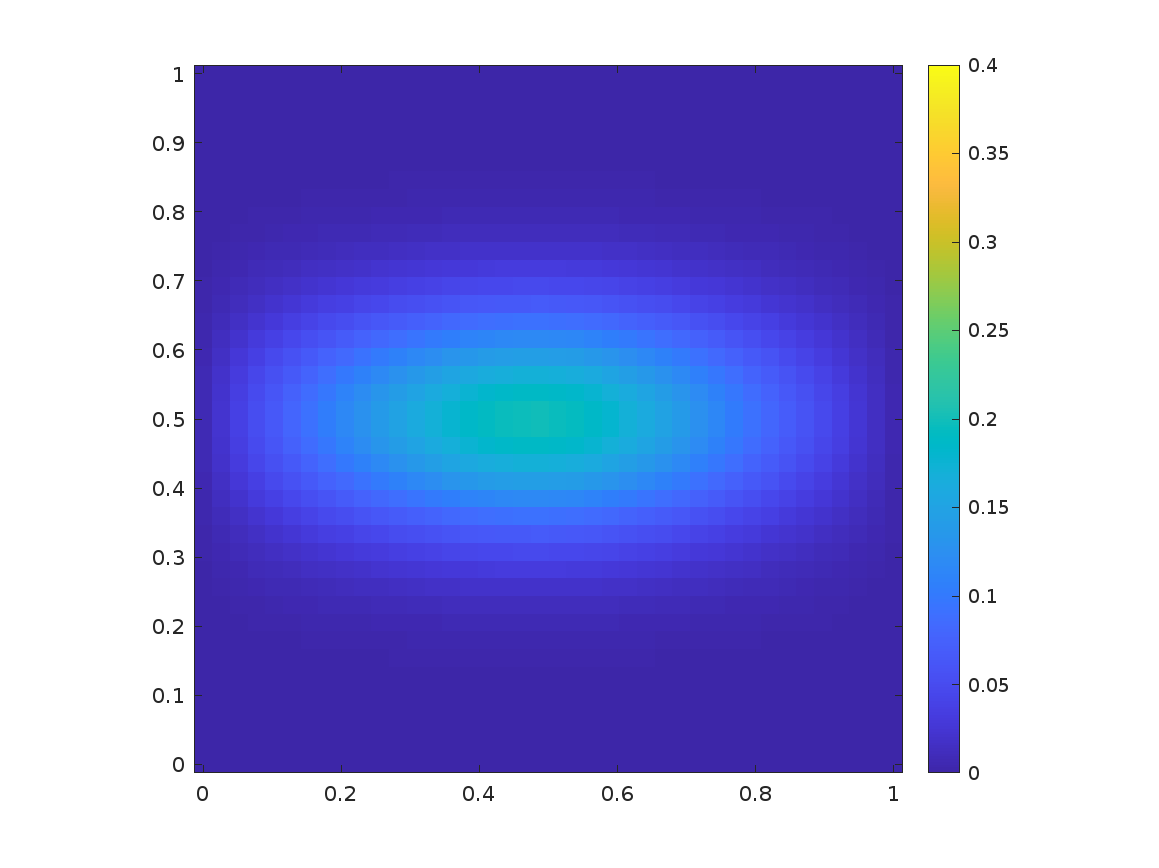}
  \quad
  \includegraphics[width=5.5cm]{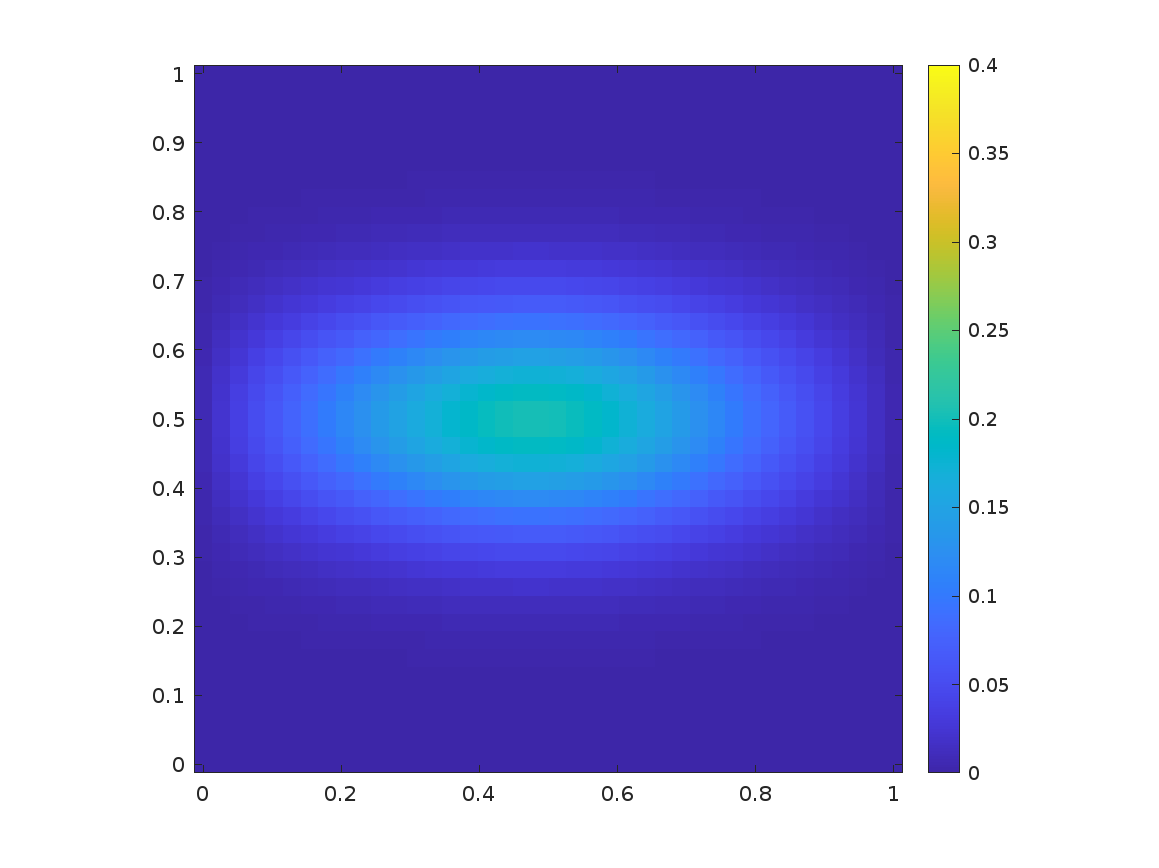}
  \quad
  \includegraphics[width=5.5cm]{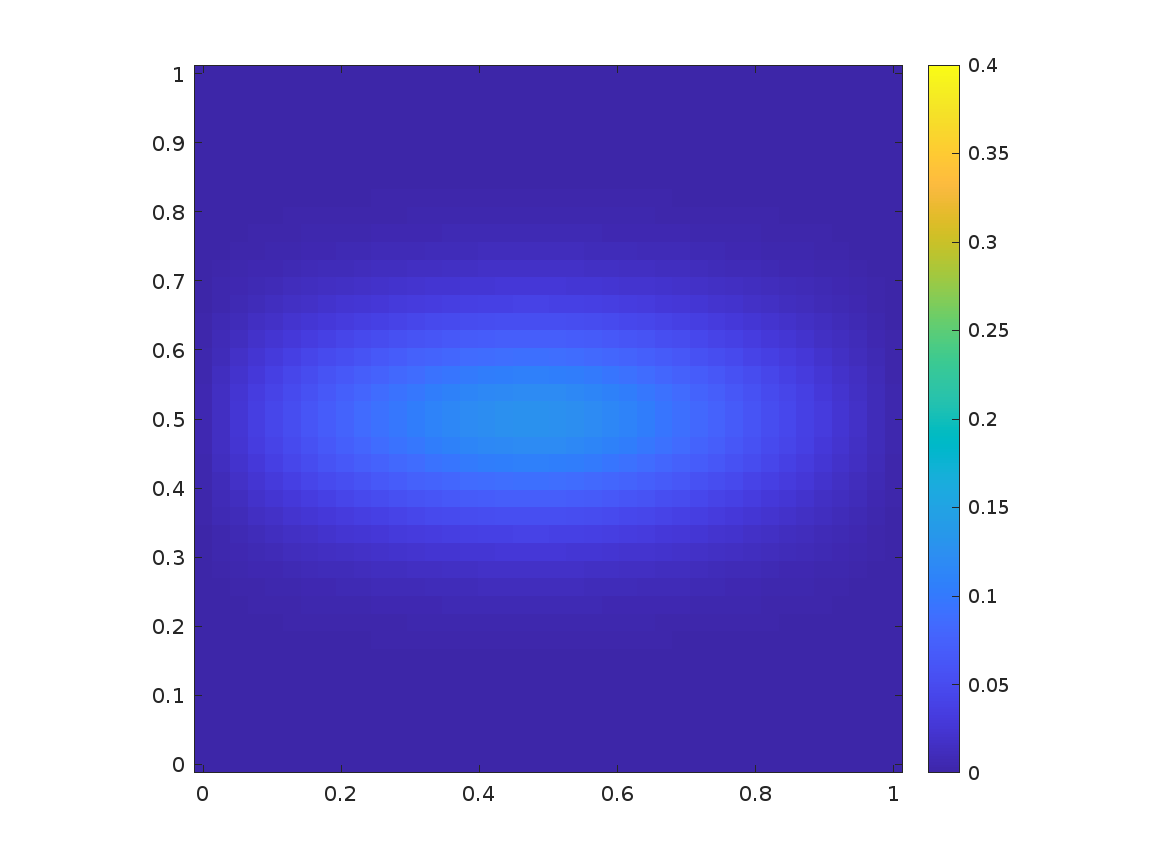}
  \quad
  \includegraphics[width=5.5cm]{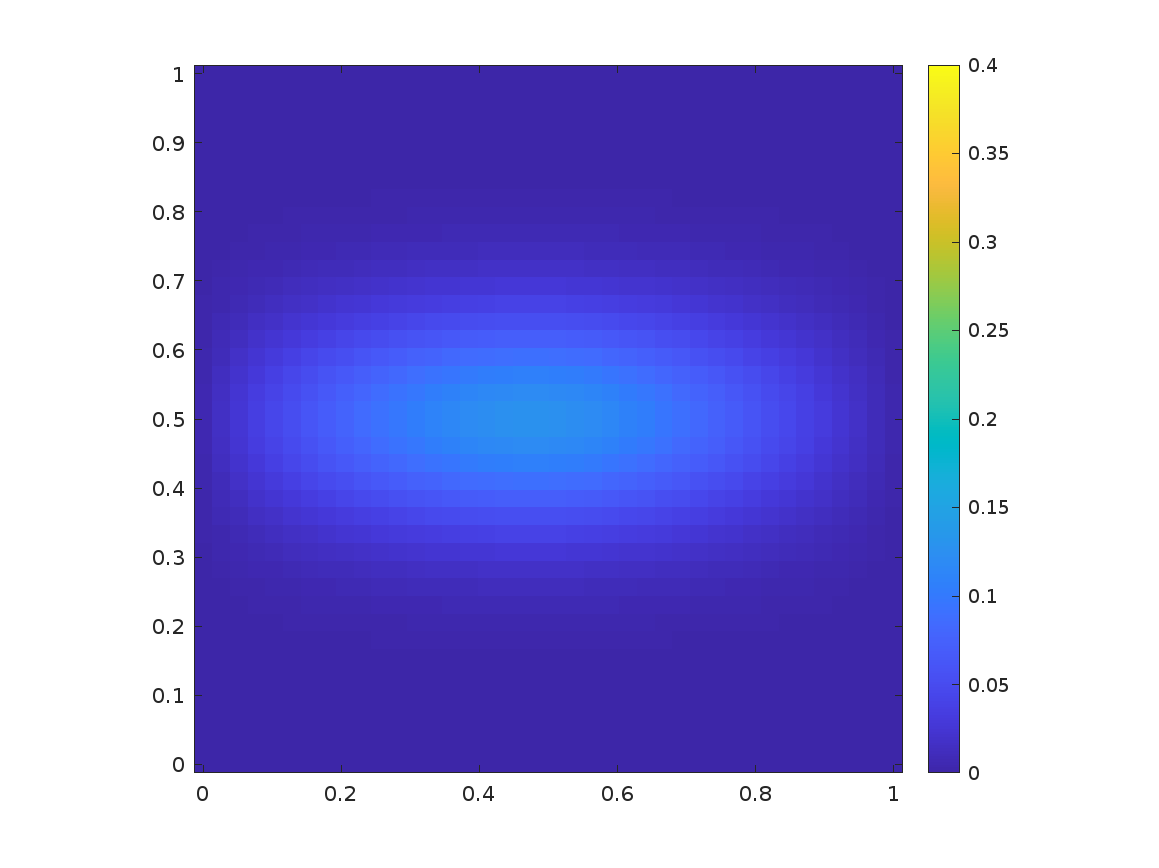}
  \quad
  \includegraphics[width=5.5cm]{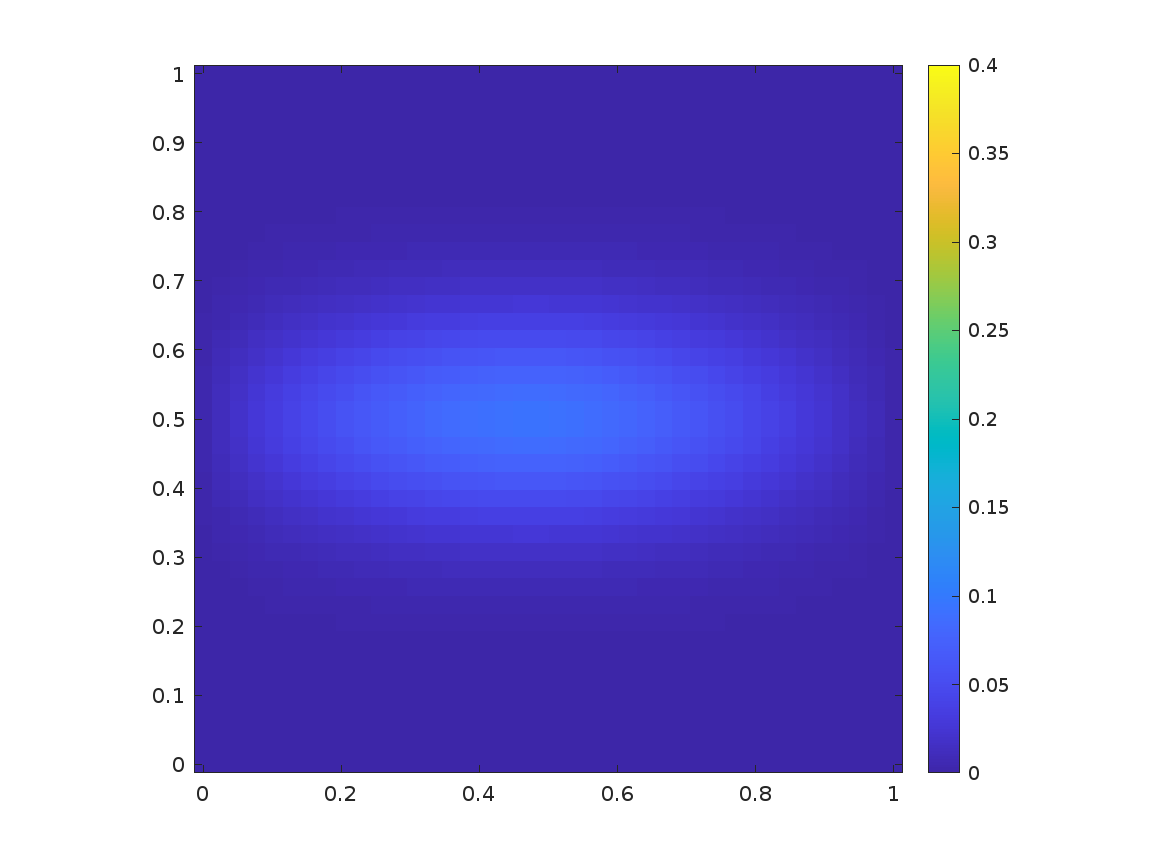}
  \quad
  \includegraphics[width=5.5cm]{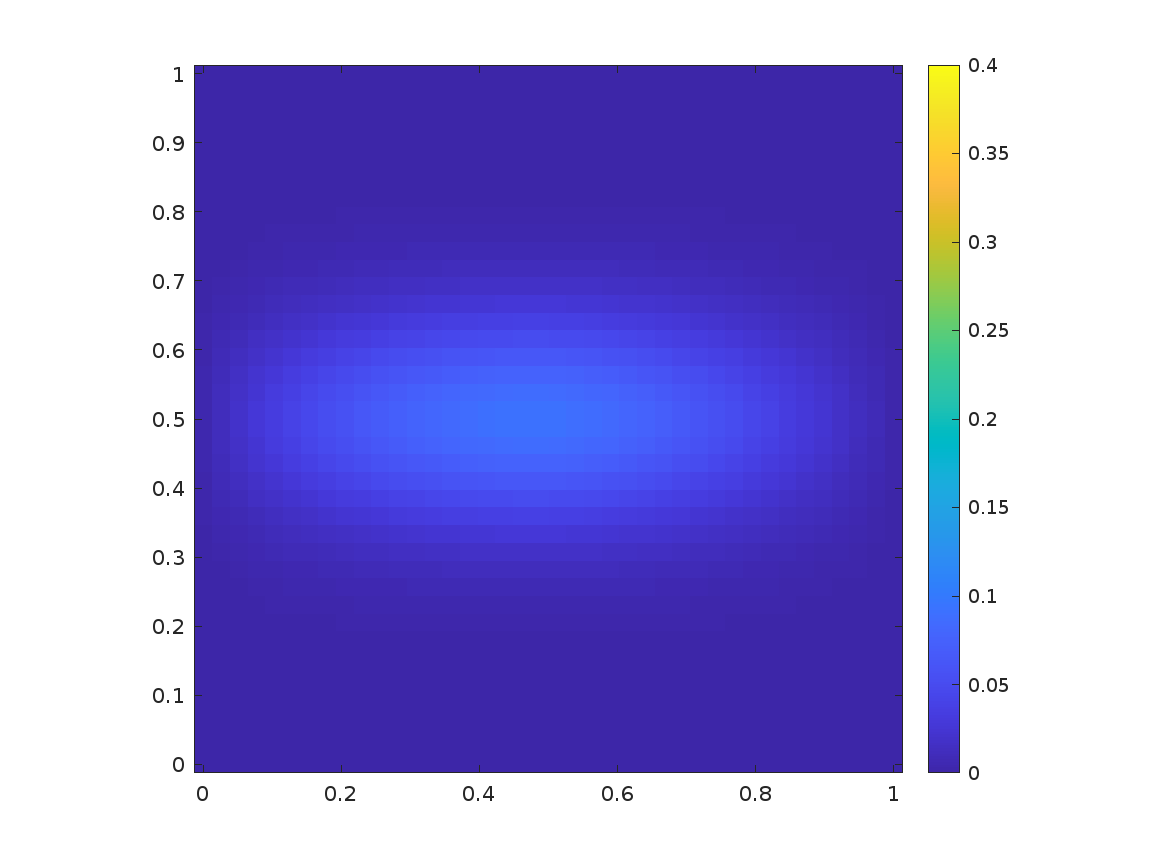}
  \quad
  \includegraphics[width=5.5cm]{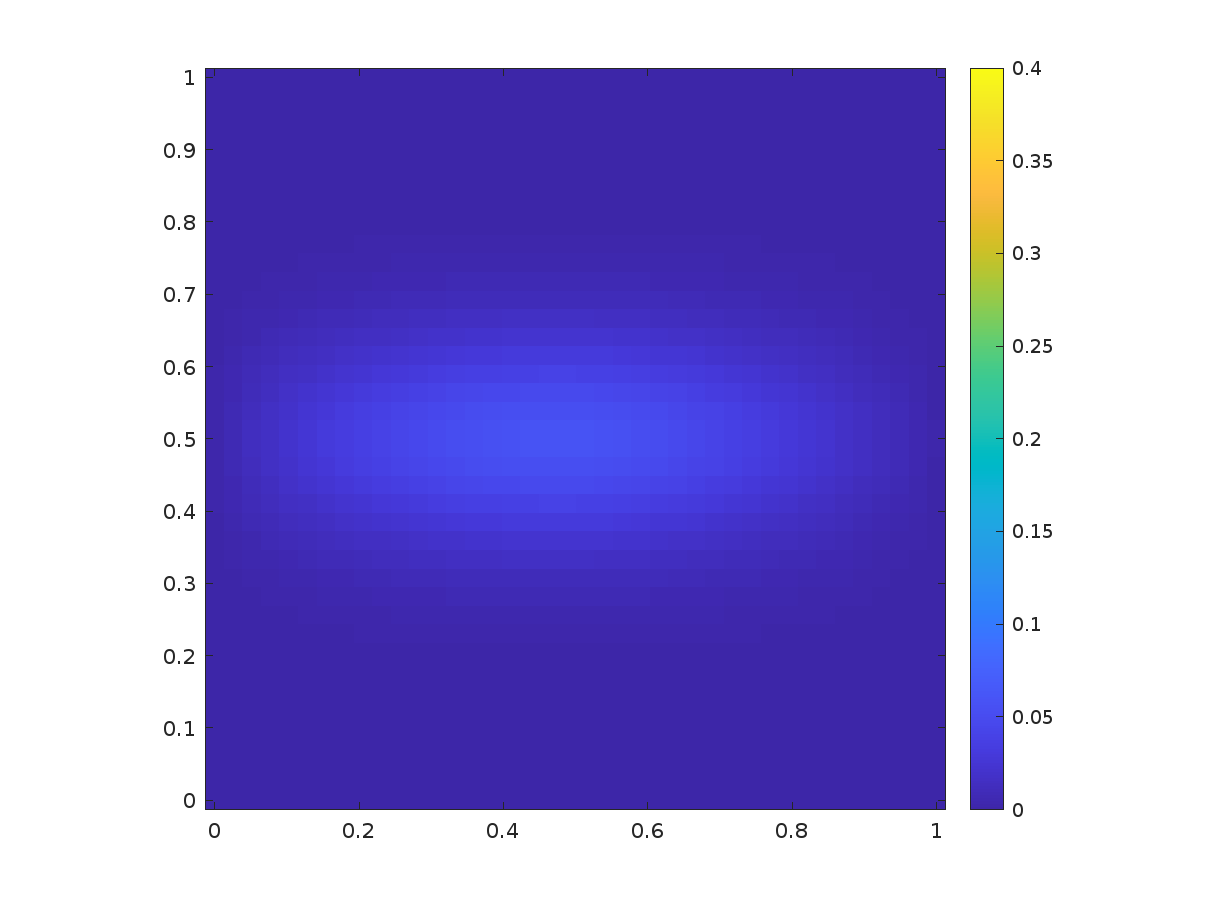}
  \quad
  \includegraphics[width=5.5cm]{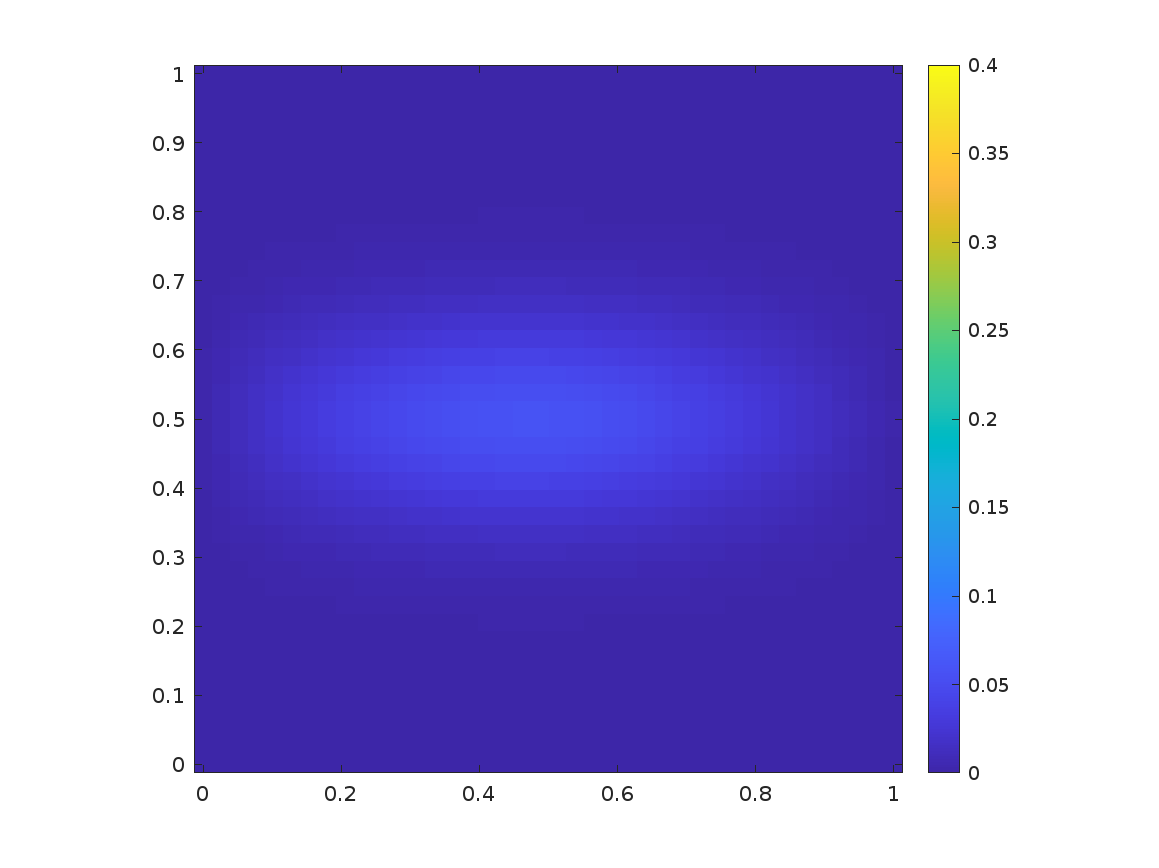}
  \caption{Solutions of concentration when $H=1/40$ for Case 2 in Example 1. First column: multiscale solution $C_2$ at $t=0.02$, $0.1$, $0.5$, $1$, $2$. Second column: reference averaged solution in $\Omega_2$ at the corresponding time instants.}
  \label{fig:Example1_Case2_U2}
\end{figure}

We list the relative $L^2$ errors when $H=1/20$ and $H=1/40$ in Table \ref{tab:Example1_Case2_error}. The numerical solutions when $H=1/40$ are presented in Figures \ref{fig:Example1_Case2_snapshot}, \ref{fig:Example1_Case2_U1}, and \ref{fig:Example1_Case2_U2}. 
One can notice that our method yields results that closely align with the reference solutions in this case. 
We can also observe from Figure \ref{fig:Example1_Case2_snapshot} that the concentration exhibits faster transport in the $-x$-direction compared to the $x$-direction. 
Indeed, due to the imposed boundary condition (\ref{eq:Example1_Case2_bc}), the pressure generally increases along the $x$-direction on average. According to Darcy's law, the fluid predominantly flows toward regions of decreasing pressure, thereby carrying substances in the $-x$-direction and offsetting the effects of diffusion.

\begin{table}
\caption{Relative $L^2$ errors at $t=0.02$, $0.1$, $0.5$, $1$, $2$ for Case 2 in Example 1. Left: $H=1/20$ and $l=6$. Right: $H=1/40$ and $l=8$.}
\centering
\begin{tabular}{ccc}
   \toprule
   $t$      & $e^{(1)}(t)$    &   $e^{(2)}(t)$ \\
   \midrule
   $0.02$ & $4.64\%$ & $2.17\%$ \\
   $0.1$ & $3.03\%$ & $2.89\%$ \\
   $0.5$ & $2.54\%$ & $1.87\%$ \\
   $1.0$ & $2.39\%$ & $1.50\%$ \\
   $2.0$ & $1.75\%$ & $1.16\%$ \\
   \bottomrule 
\end{tabular}
\qquad
\begin{tabular}{ccc}
   \toprule
   $t$      & $e^{(1)}(t)$    &   $e^{(2)}(t)$ \\
   \midrule
   $0.02$ & $2.44\%$ & $2.25\%$ \\
   $0.1$ & $0.76\%$ & $1.34\%$ \\
   $0.5$ & $0.39\%$ & $0.23\%$ \\
   $1.0$ & $0.55\%$ & $0.52\%$ \\
   $2.0$ & $0.60\%$ & $0.82\%$ \\
   \bottomrule 
\end{tabular}
\label{tab:Example1_Case2_error}
\end{table}

\subsubsection{Case 3}

For Case 3, we consider a homogeneous Neumann condition for the concentration equation and the boundary conditions will be
\begin{equation}
    \begin{split}
        p(x,t) = x, \quad & x\in \partial \Omega, \quad \forall t\in [0,\infty),\\
        \nabla c(x,t) \cdot \nu = 0, \quad & x\in \partial \Omega, \quad \forall t\in [0,\infty),
    \end{split}
\end{equation}
where $\nu$ is the outward normal vector defined on $\partial \Omega$. We keep all the other assumptions the same as in Case 1.

\begin{figure}
  \centering
  \includegraphics[width=5.5cm]{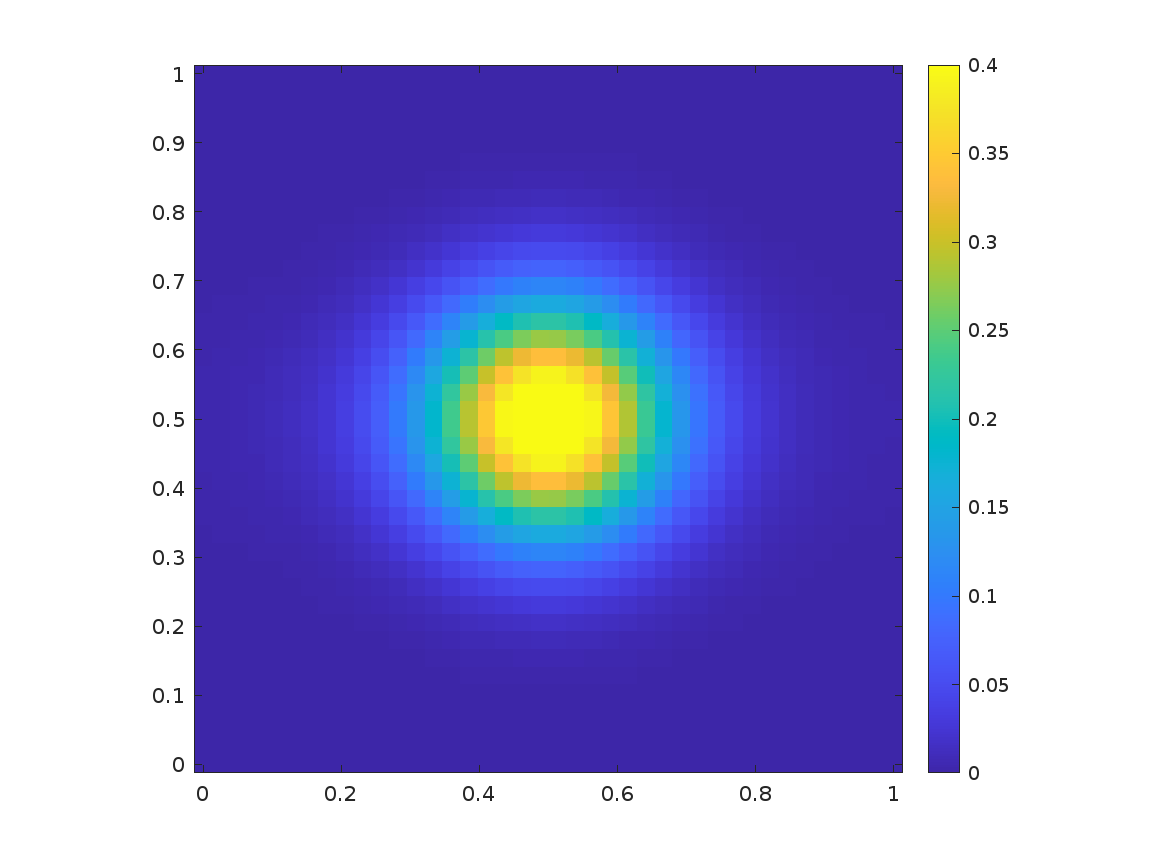}
  \quad
  \includegraphics[width=5.5cm]{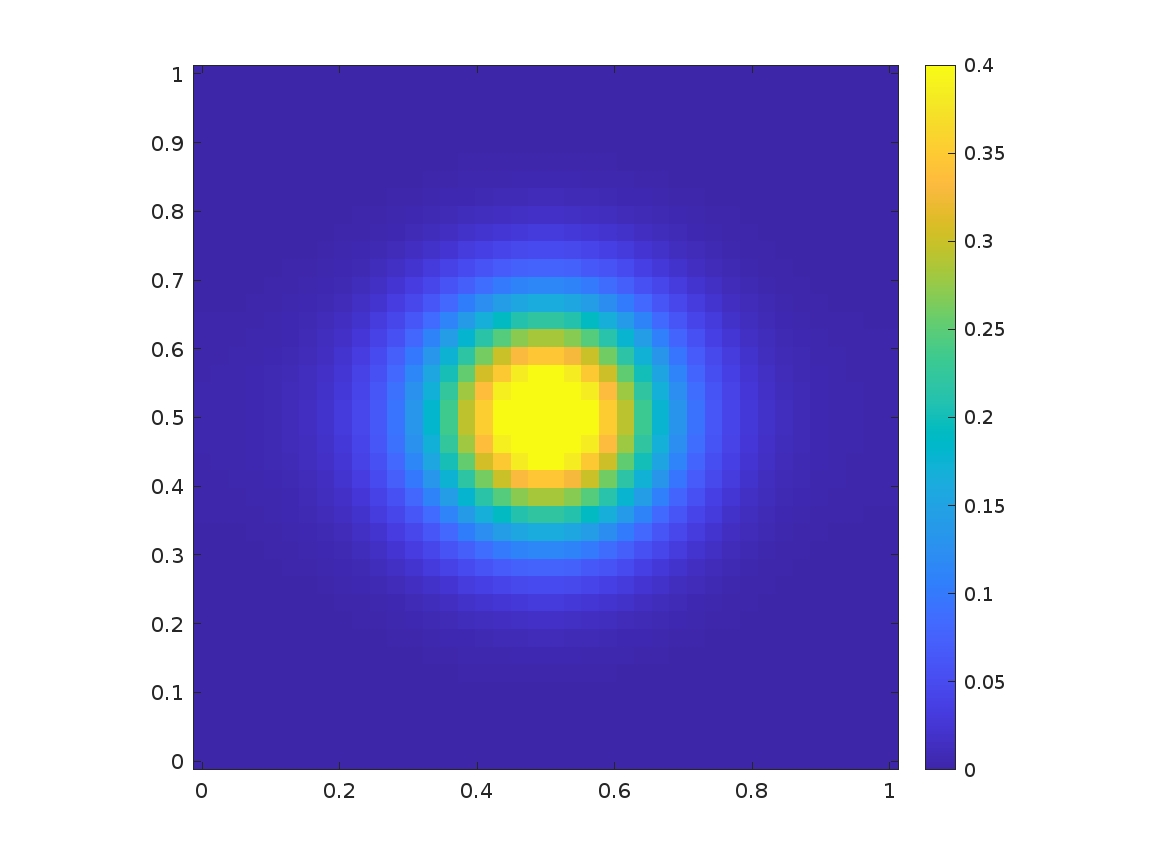}
  \quad
  \includegraphics[width=5.5cm]{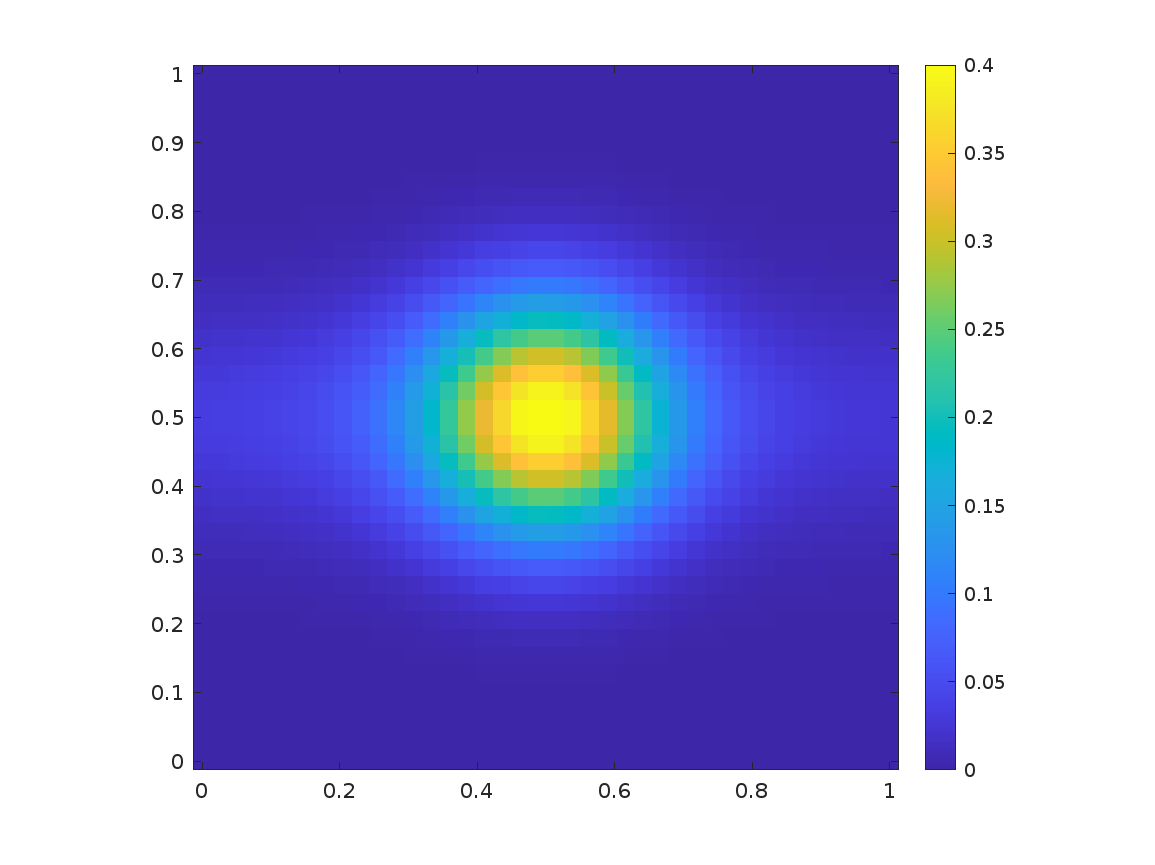}
  \quad
  \includegraphics[width=5.5cm]{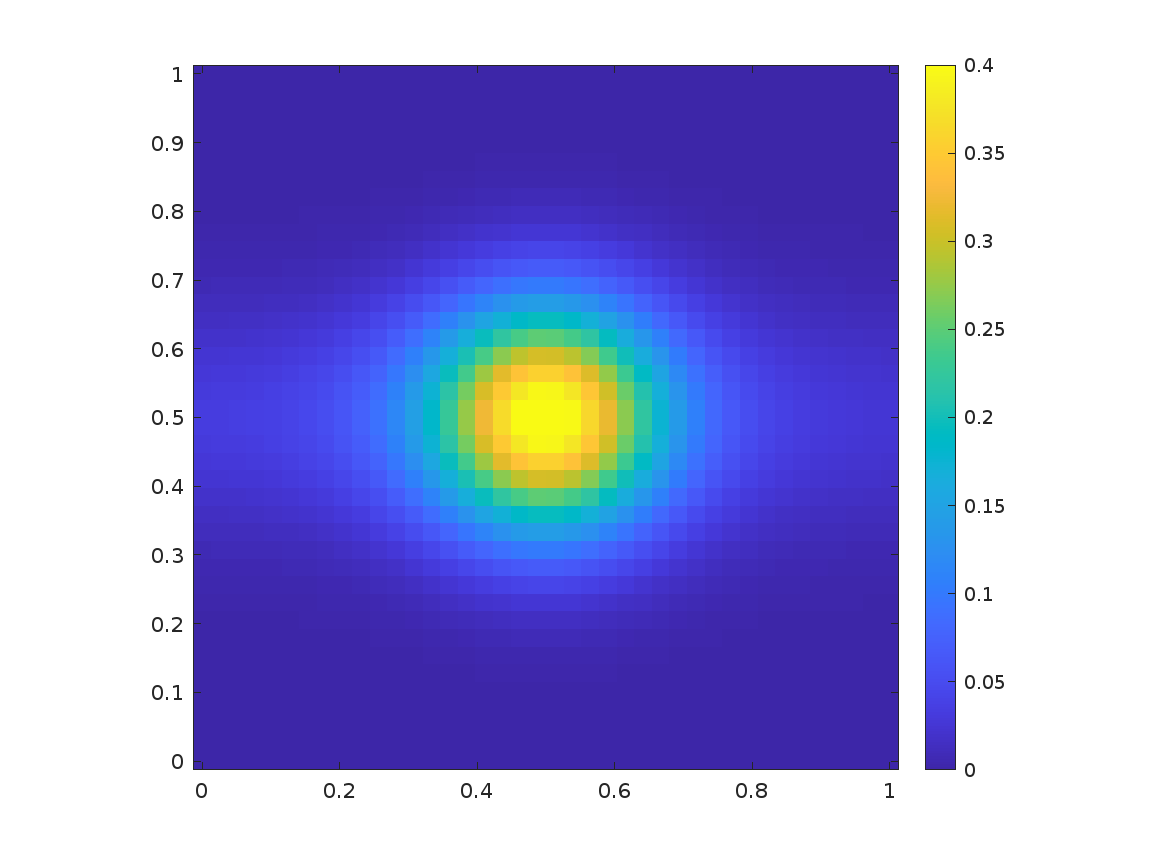}
  \quad
  \includegraphics[width=5.5cm]{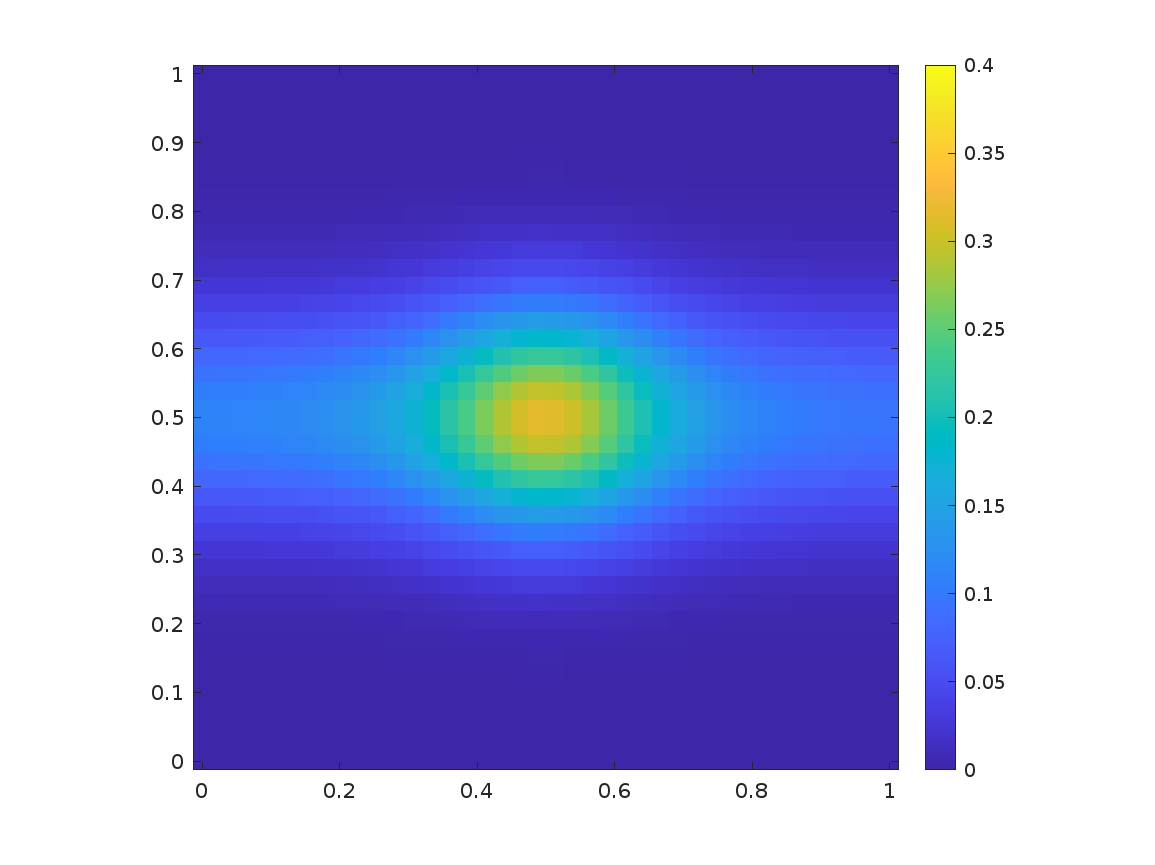}
  \quad
  \includegraphics[width=5.5cm]{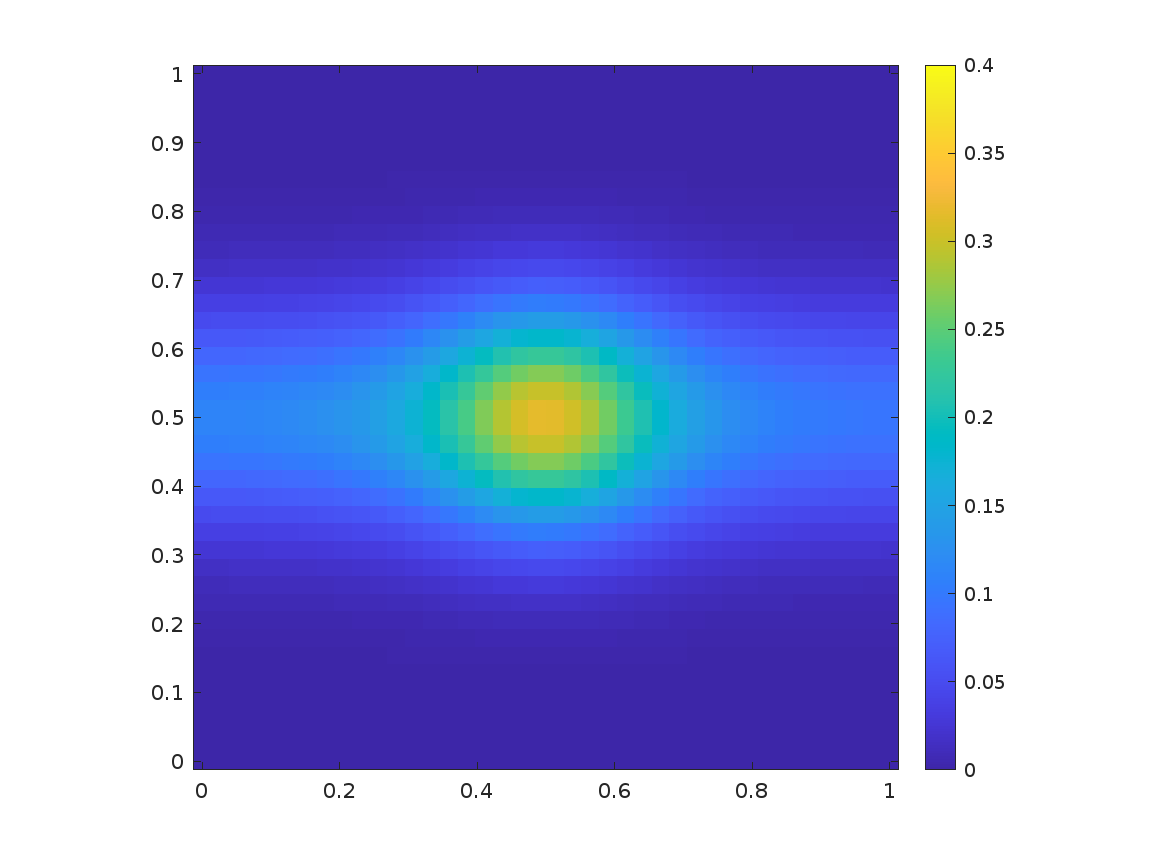}
  \quad
  \includegraphics[width=5.5cm]{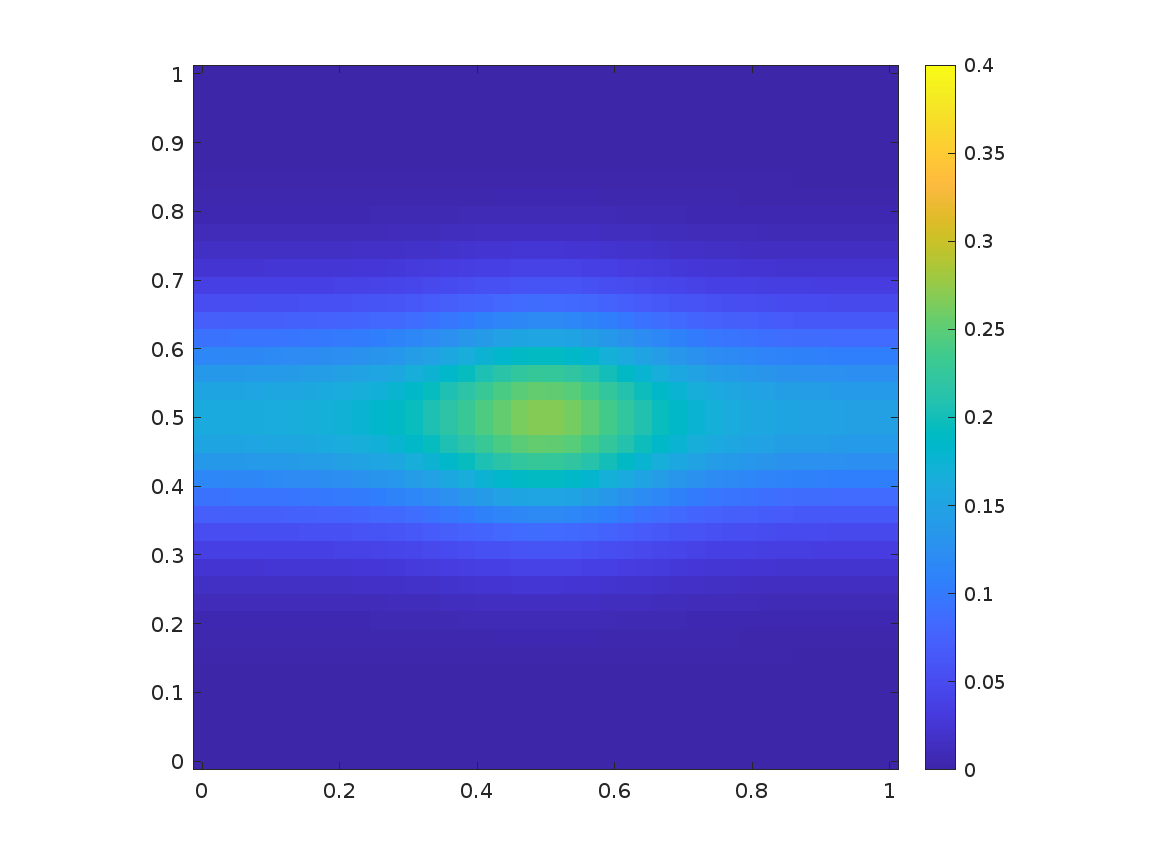}
  \quad
  \includegraphics[width=5.5cm]{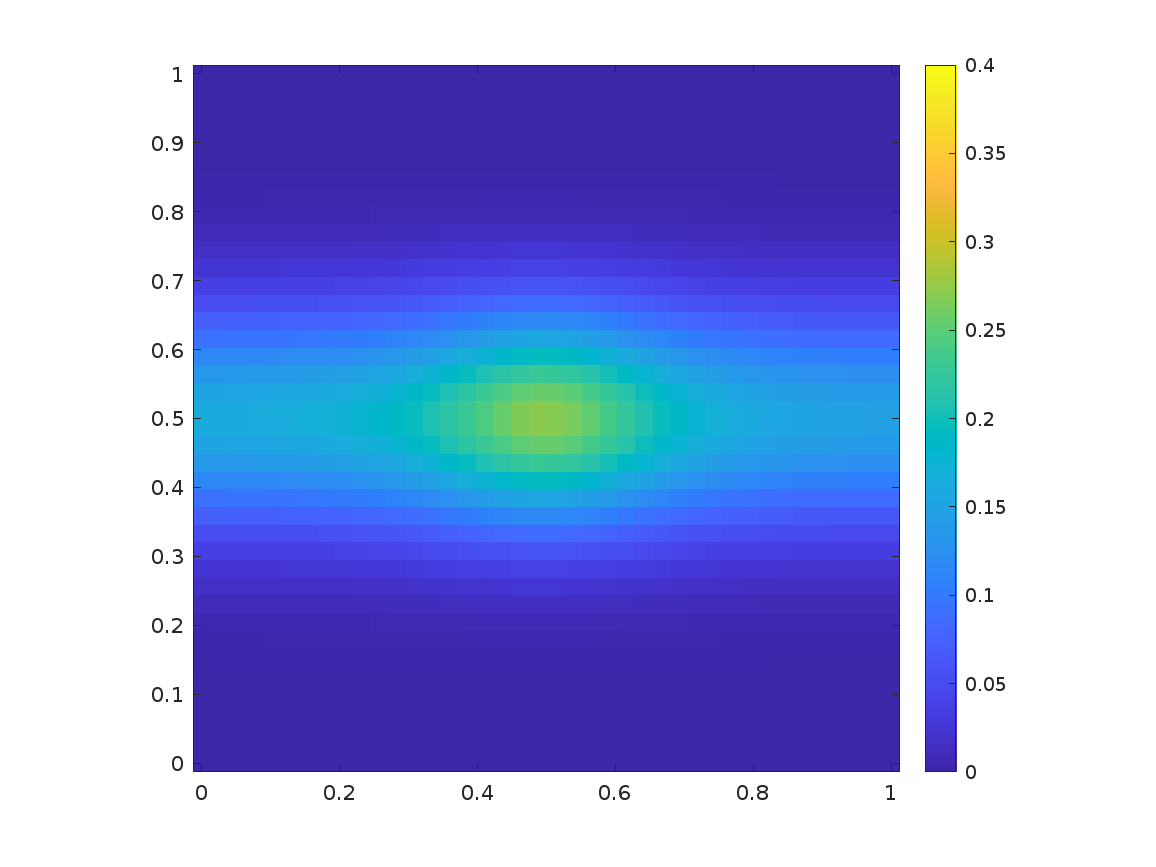}
  \quad
  \includegraphics[width=5.5cm]{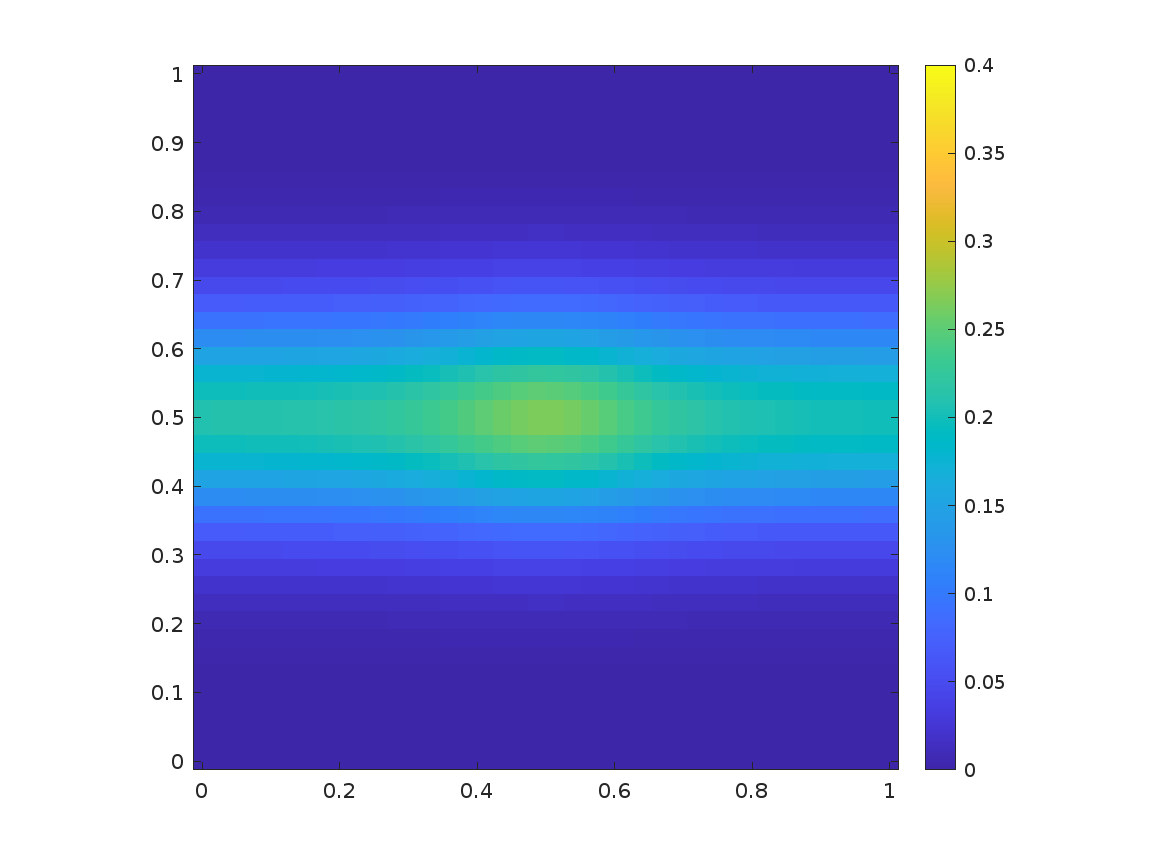}
  \quad
  \includegraphics[width=5.5cm]{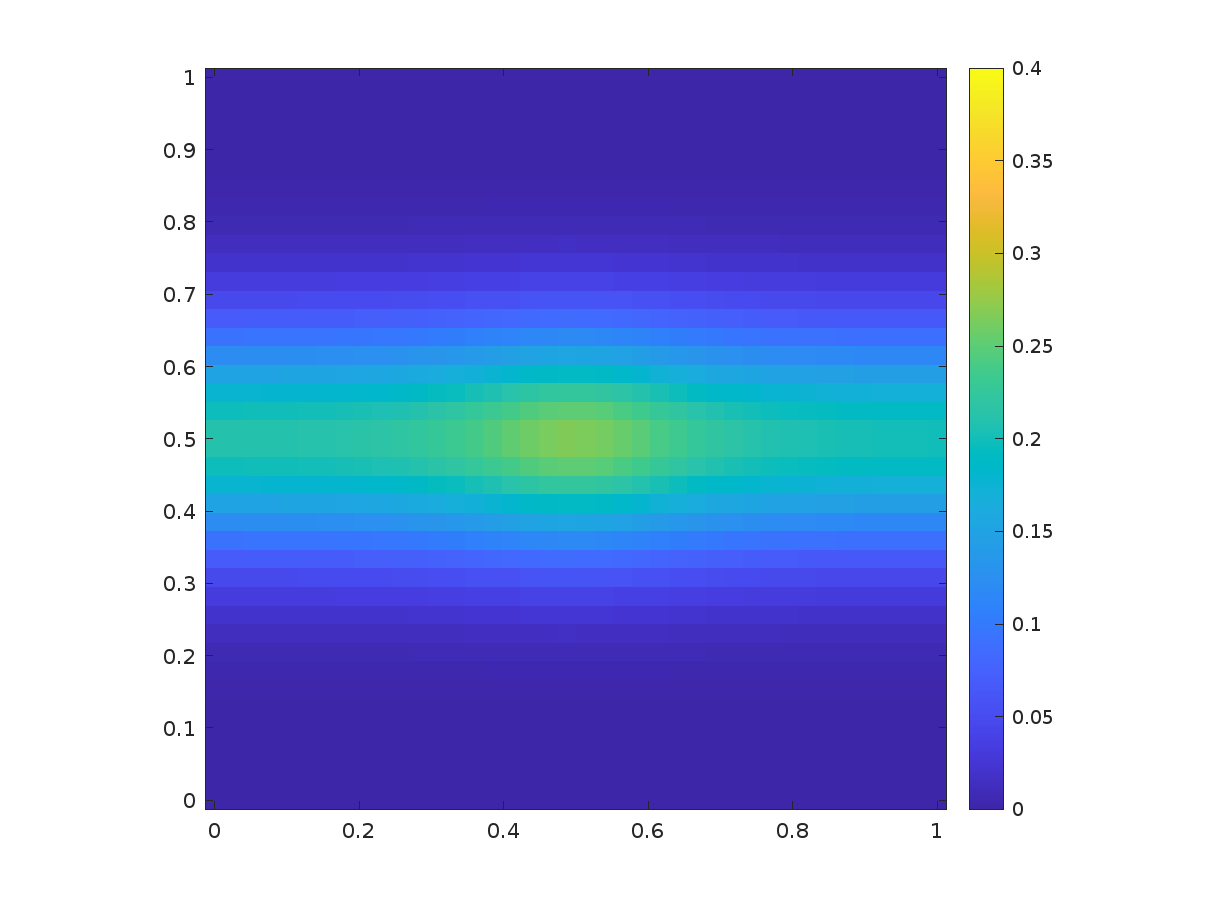}
  \caption{Solutions of concentration when $H=1/40$ for Case 3 in Example 1. First column: multiscale solution $C_1$ at $t=0.02$, $0.1$, $0.5$, $1$, $2$. Second column: reference averaged solution in $\Omega_1$ at the corresponding time instants.}
  \label{fig:Example1_Case3_U1}
\end{figure}

\begin{figure}
  \centering
  \includegraphics[width=5.5cm]{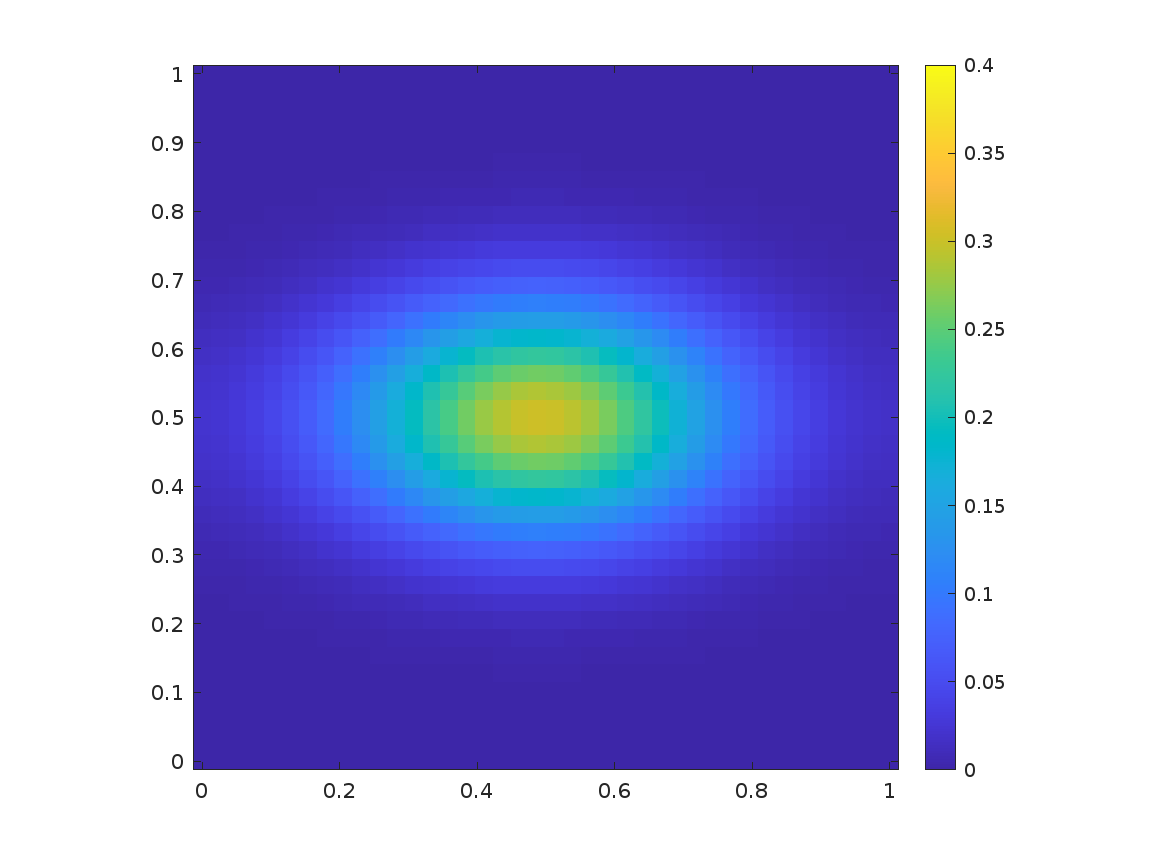}
  \quad
  \includegraphics[width=5.5cm]{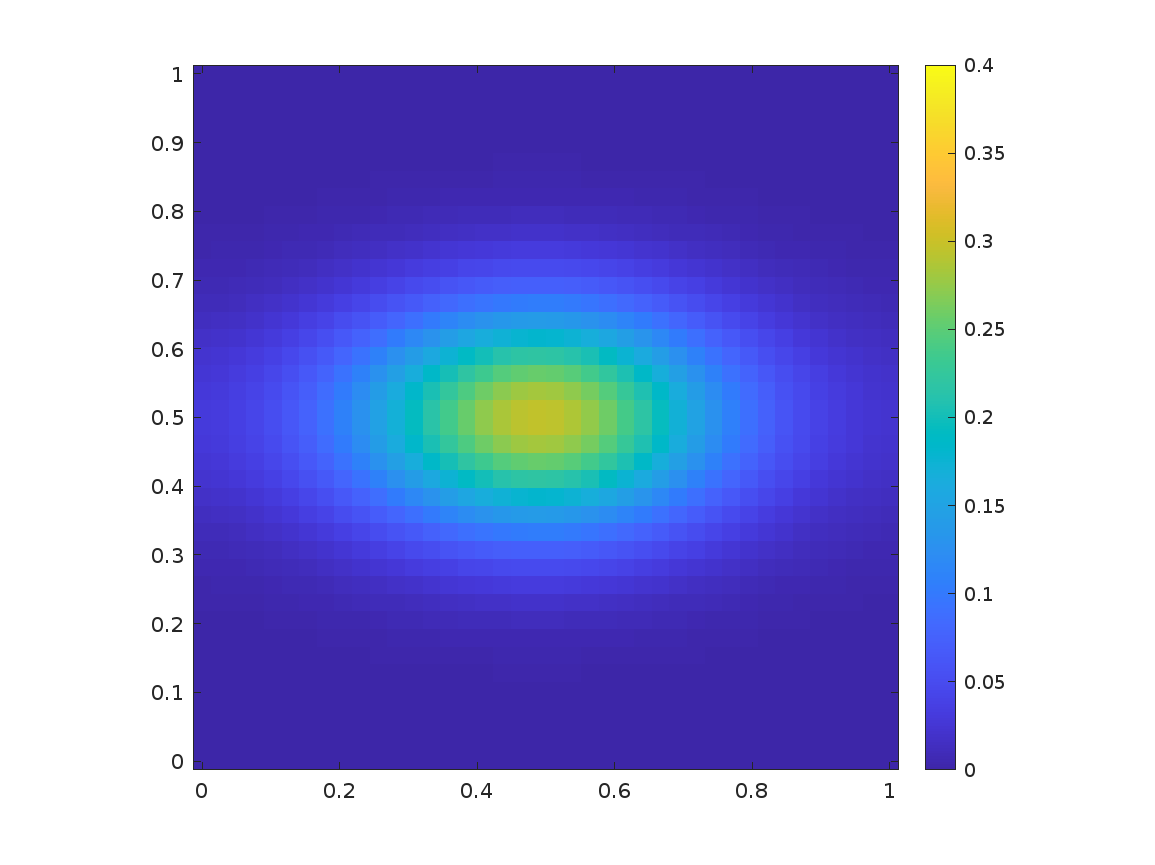}
  \quad
  \includegraphics[width=5.5cm]{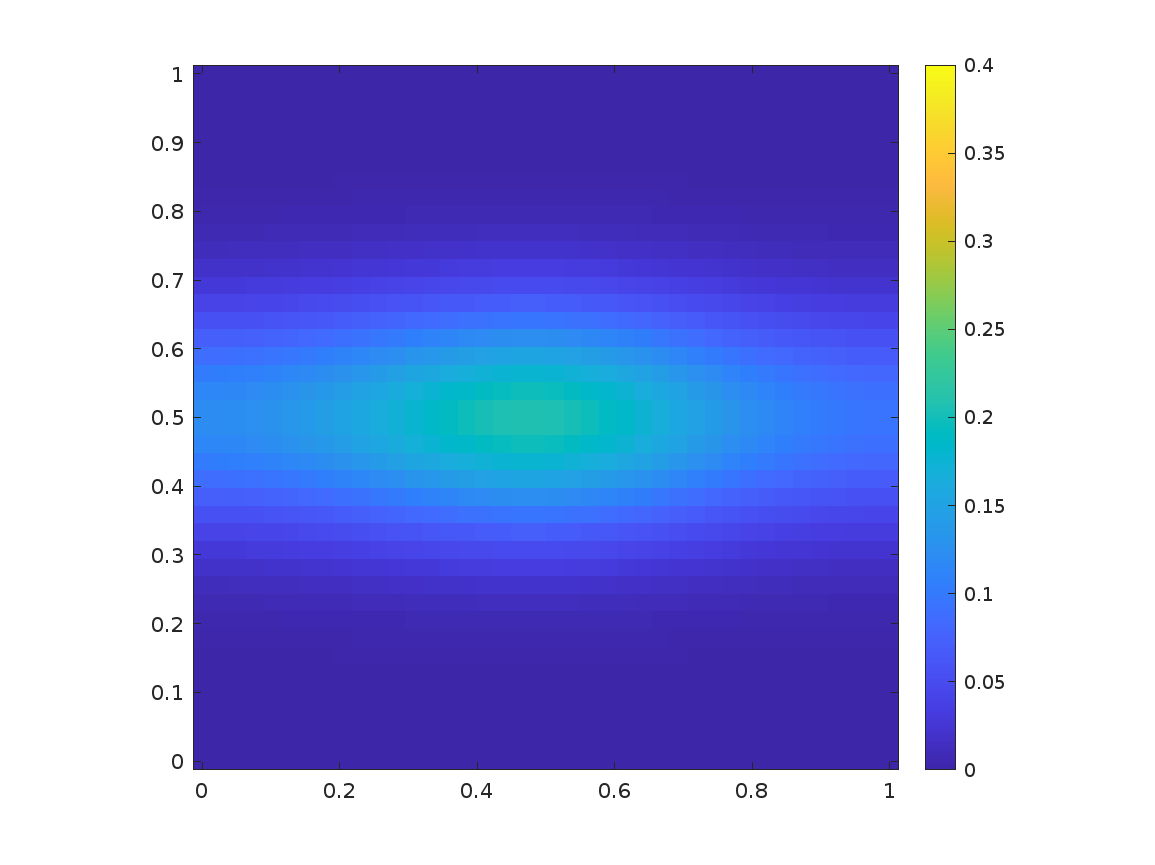}
  \quad
  \includegraphics[width=5.5cm]{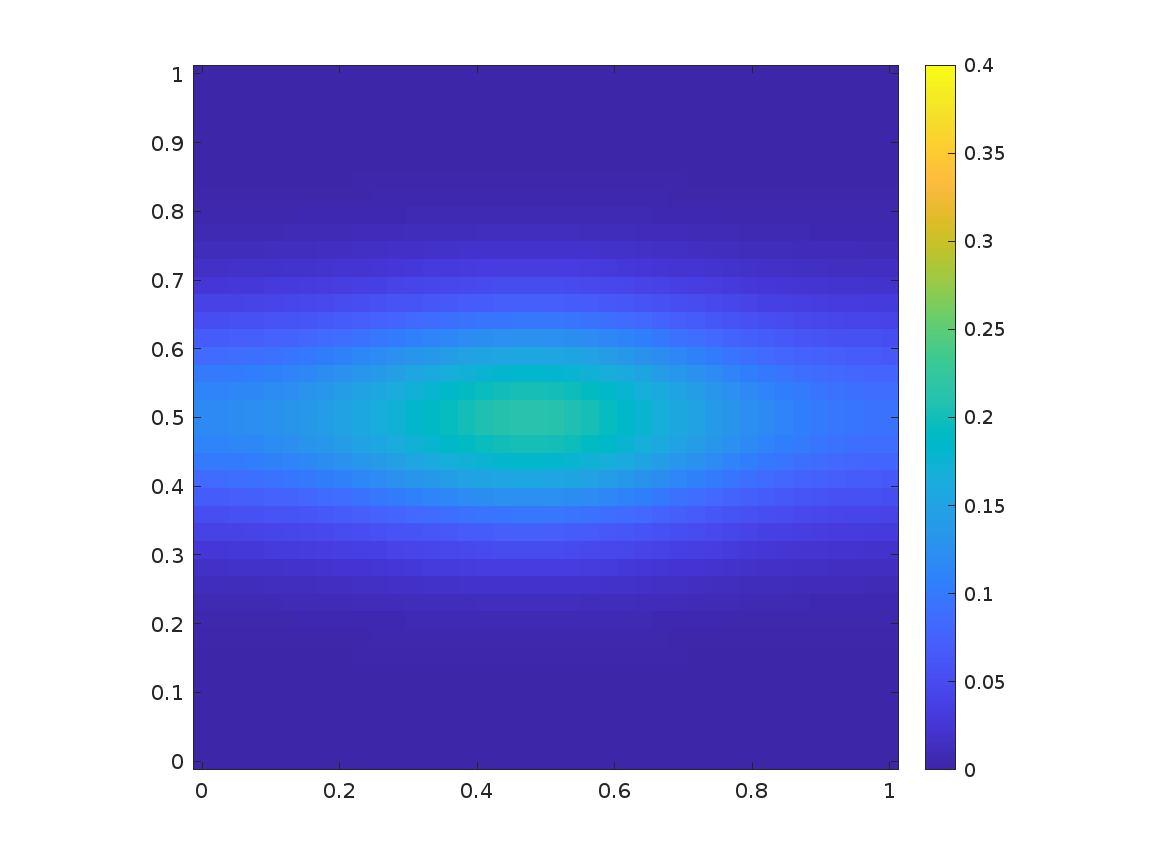}
  \quad
  \includegraphics[width=5.5cm]{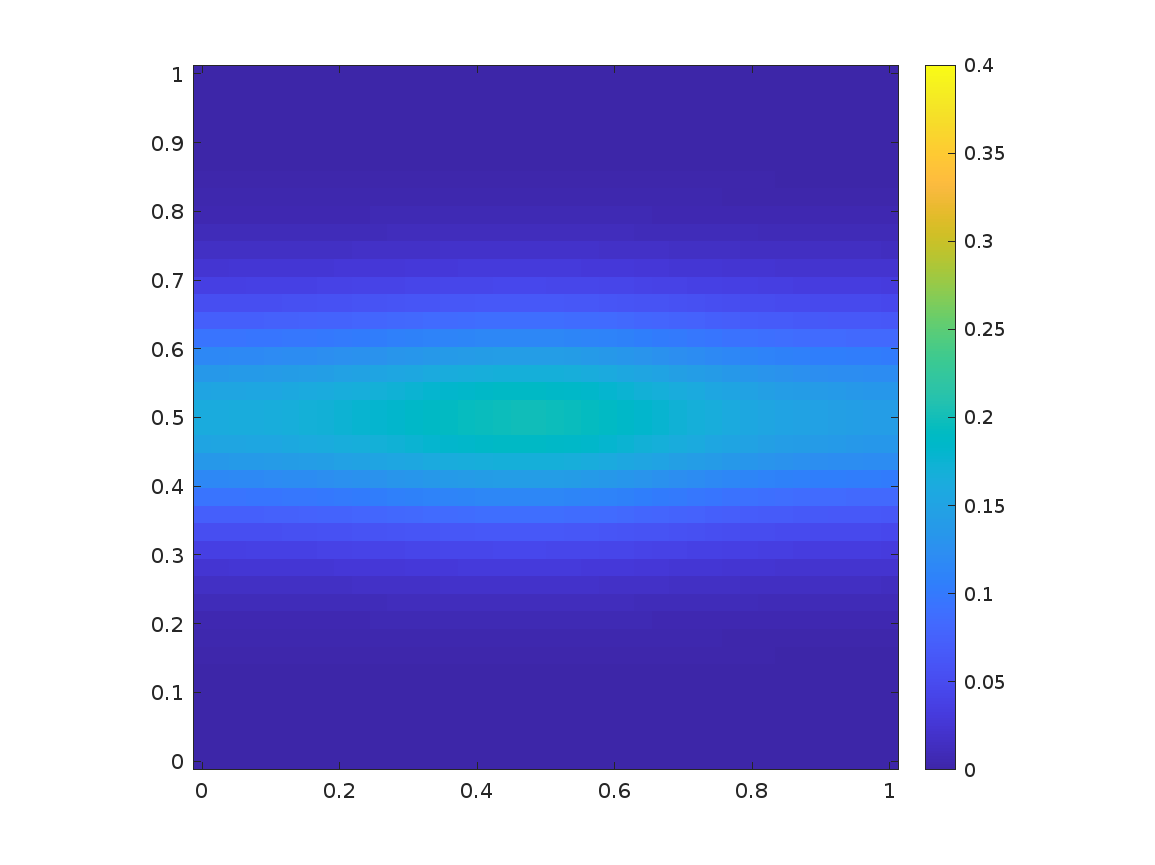}
  \quad
  \includegraphics[width=5.5cm]{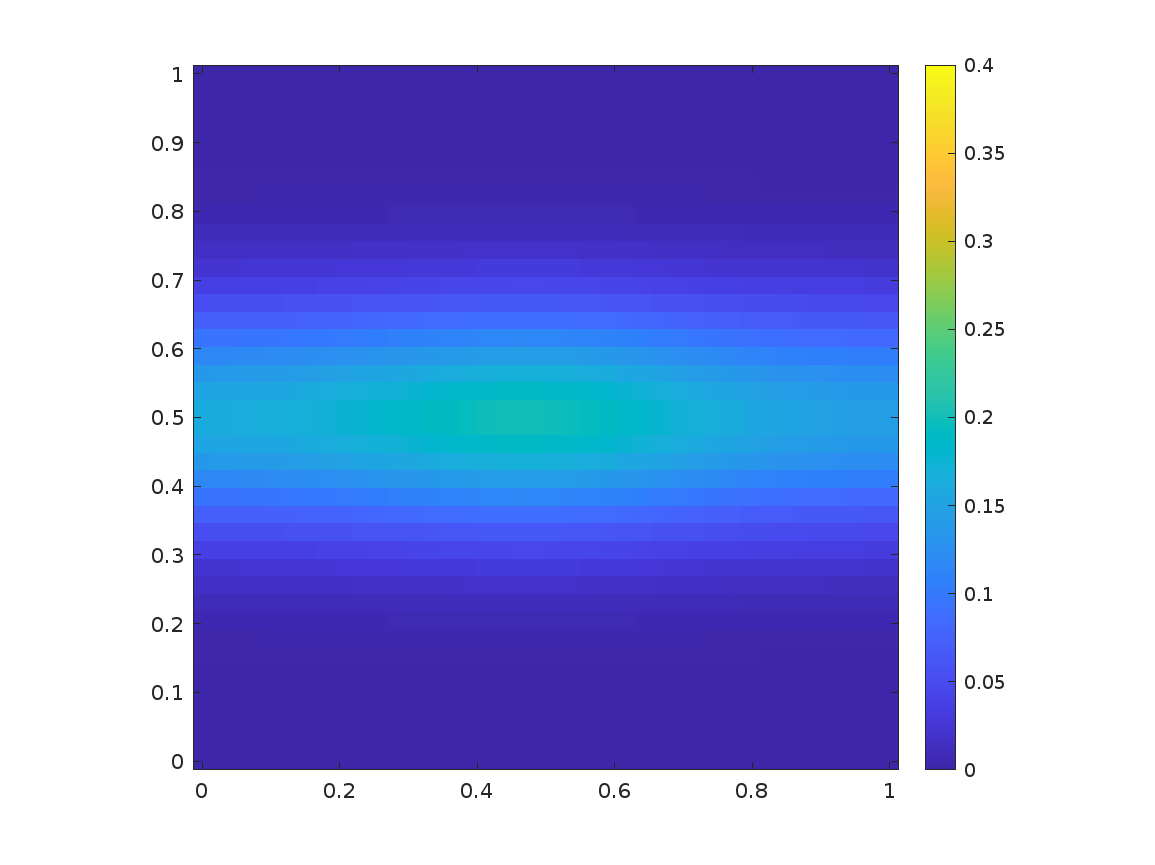}
  \quad
  \includegraphics[width=5.5cm]{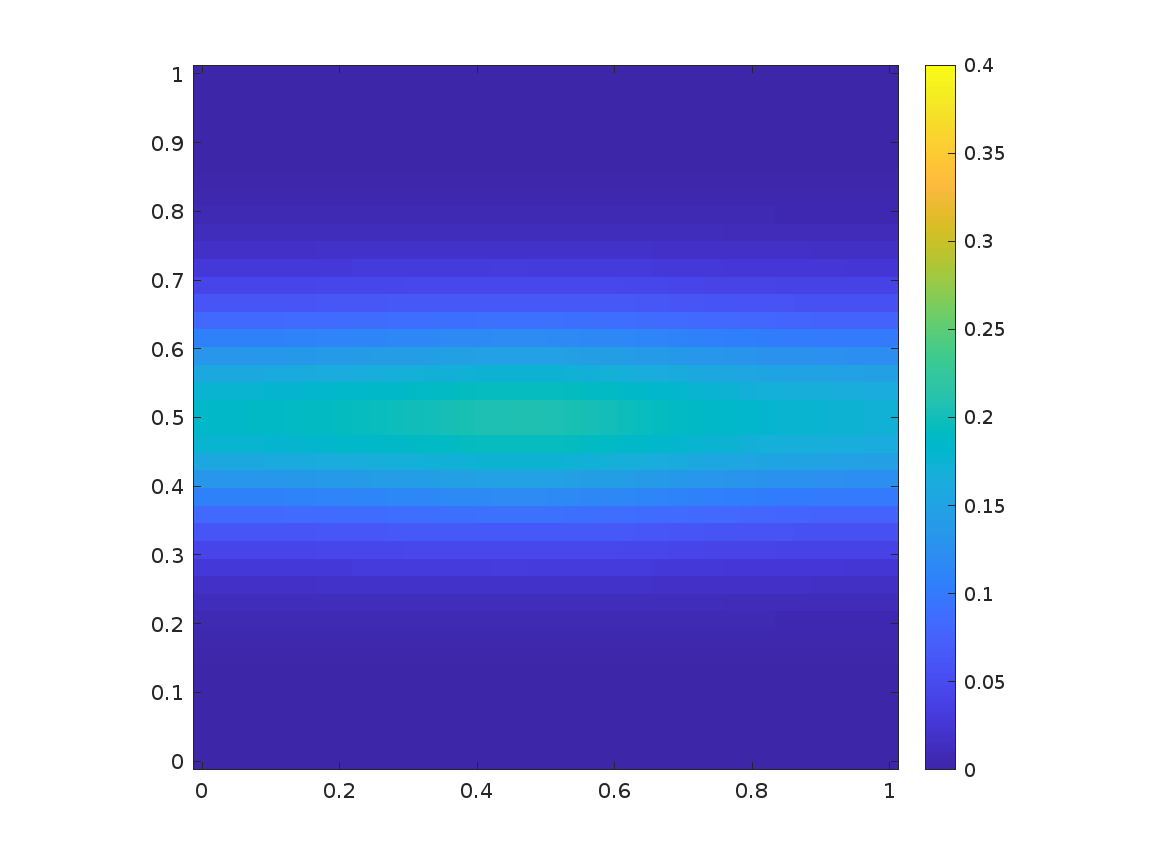}
  \quad
  \includegraphics[width=5.5cm]{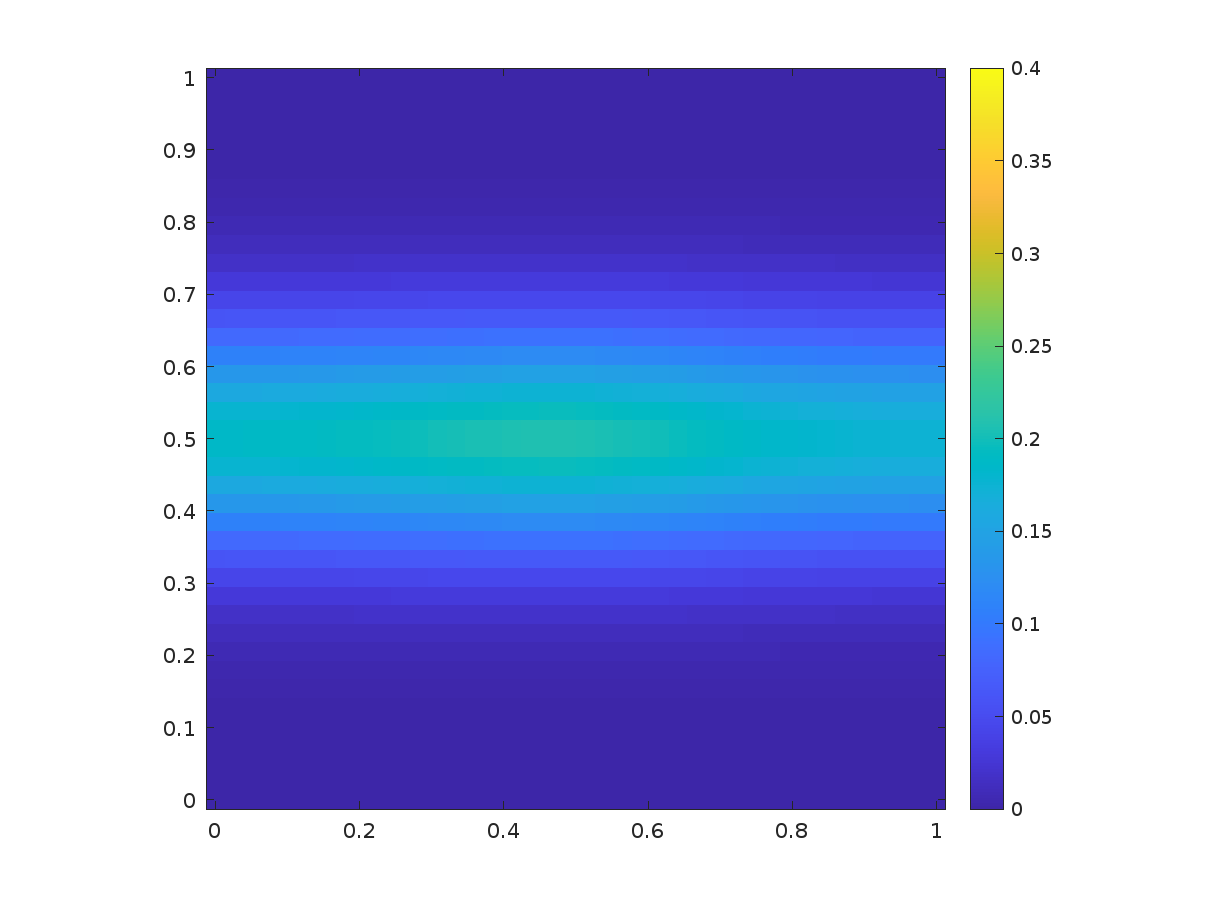}
  \quad
  \includegraphics[width=5.5cm]{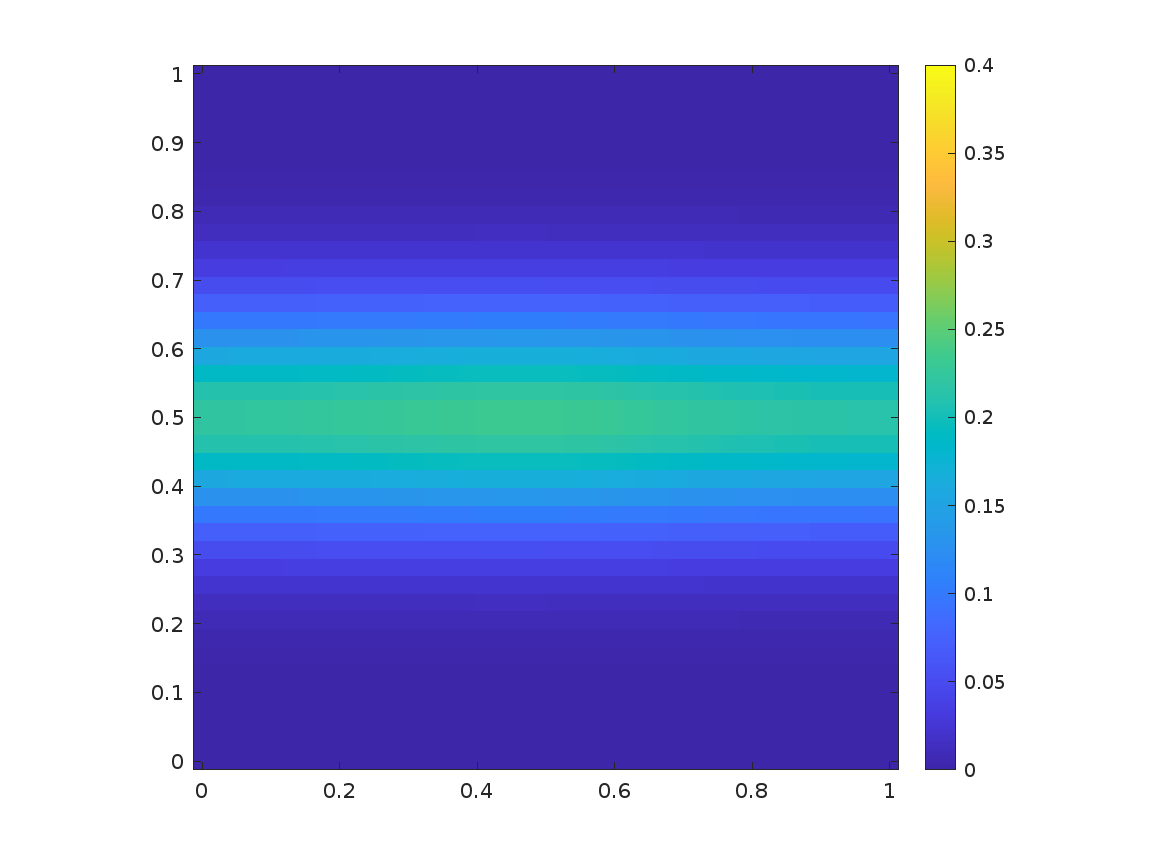}
  \quad
  \includegraphics[width=5.5cm]{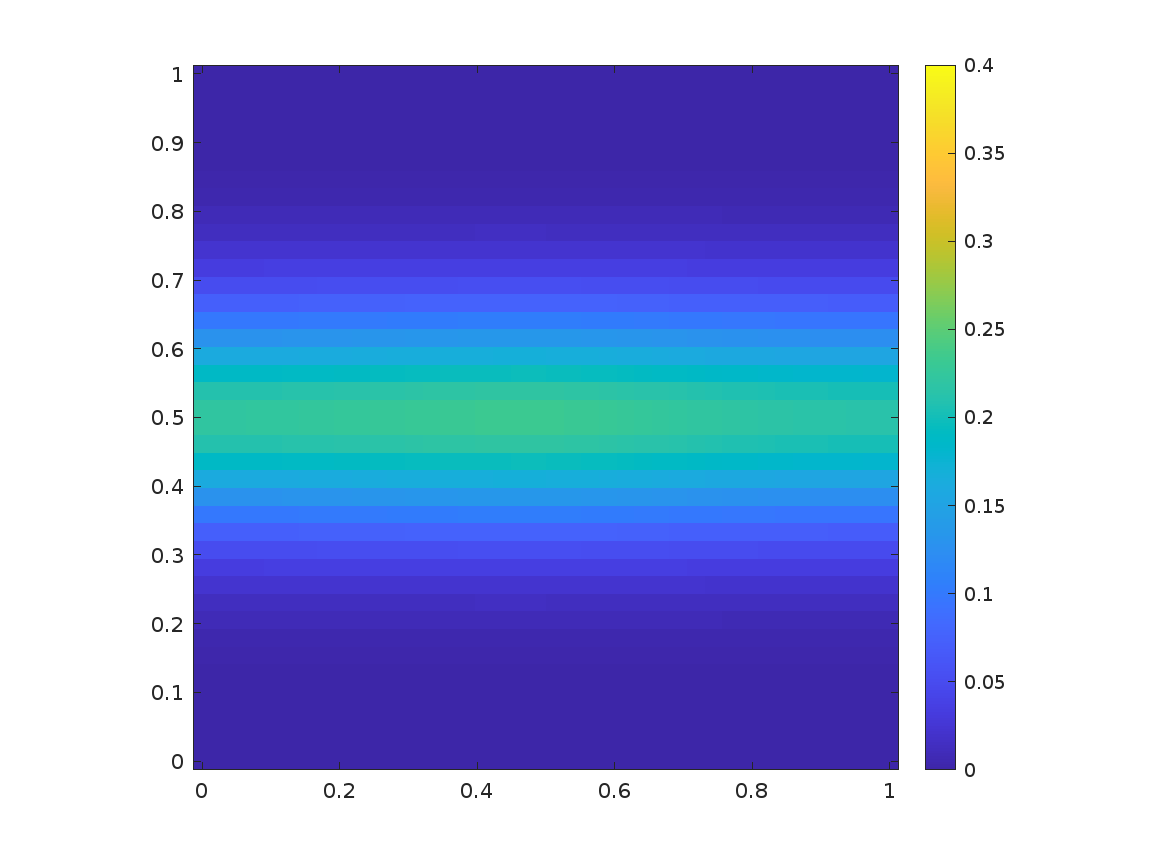}
  \caption{Solutions of concentration when $H=1/40$ for Case 3 in Example 1. First column: multiscale solution $C_2$ at $t=0.02$, $0.1$, $0.5$, $1$, $2$. Second column: reference averaged solution in $\Omega_2$ at the corresponding time instants.}
  \label{fig:Example1_Case3_U2}
\end{figure}

The relative $L^2$ errors when $H=1/20$ and $H=1/40$ are presented in Table \ref{tab:Example1_Case3_error}. The numerical results when $H=1/40$ are shown in Figures \ref{fig:Example1_Case3_U1} and \ref{fig:Example1_Case3_U2}. Our method ensures negligible discrepancy from the reference solution, and the convection can still be observed.

\begin{table}
\caption{Relative $L^2$ errors at $t=0.02$, $0.1$, $0.5$, $1$, $2$ for Case 3 in Example 1. Left: $H=1/20$ and $l=6$. Right: $H=1/40$ and $l=8$.}
\centering
\begin{tabular}{ccc}
   \toprule
   $t$      & $e^{(1)}(t)$    &   $e^{(2)}(t)$ \\
   \midrule
   $0.02$ & $4.64\%$ & $2.90\%$ \\
   $0.1$ & $3.04\%$ & $3.23\%$ \\
   $0.5$ & $2.40\%$ & $1.99\%$ \\
   $1.0$ & $2.11\%$ & $1.92\%$ \\
   $2.0$ & $1.77\%$ & $1.72\%$ \\
   \bottomrule 
\end{tabular}
\qquad
\begin{tabular}{ccc}
   \toprule
   $t$      & $e^{(1)}(t)$    &   $e^{(2)}(t)$ \\
   \midrule
   $0.02$ & $2.44\%$ & $3.00\%$ \\
   $0.1$ & $0.83\%$ & $1.82\%$ \\
   $0.5$ & $0.46\%$ & $0.24\%$ \\
   $1.0$ & $0.41\%$ & $0.20\%$ \\
   $2.0$ & $0.24\%$ & $0.20\%$ \\
   \bottomrule 
\end{tabular}
\label{tab:Example1_Case3_error}
\end{table}

\subsection{Example 2: Circular field}

\begin{figure}[htbp]
    \centering
    \includegraphics[width = 8cm]{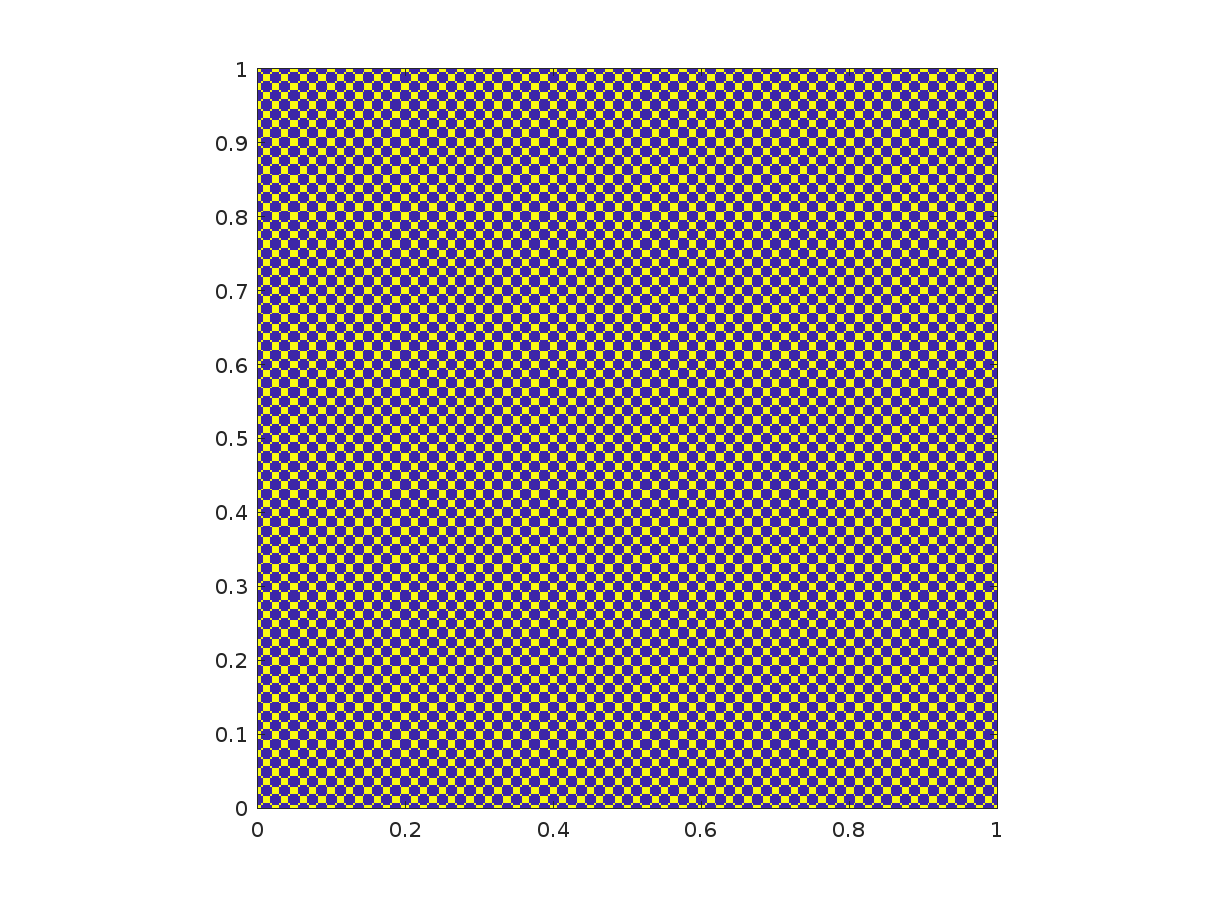}
    \caption{Circular field for Example 2 ($\Omega_1$: blue regions; $\Omega_2$: yellow regions)}
    \label{fig:Example2_field}
\end{figure}

In this example, we consider a different structure of permeability field $\kappa$ and the diffusion field $D$ as illustrated in Figure \ref{fig:Example2_field} and 
\begin{equation}
\label{eq:Example2_coefficient}
    \kappa(x) = D(x) = 
    \begin{cases}
        10^{-4}, \quad & x\in \Omega_1,\\
        1, \quad & x \in \Omega_2.
    \end{cases}
\end{equation}
All other assumptions remain unchanged. 

\begin{figure}
  \centering
  \includegraphics[width=5.5cm]{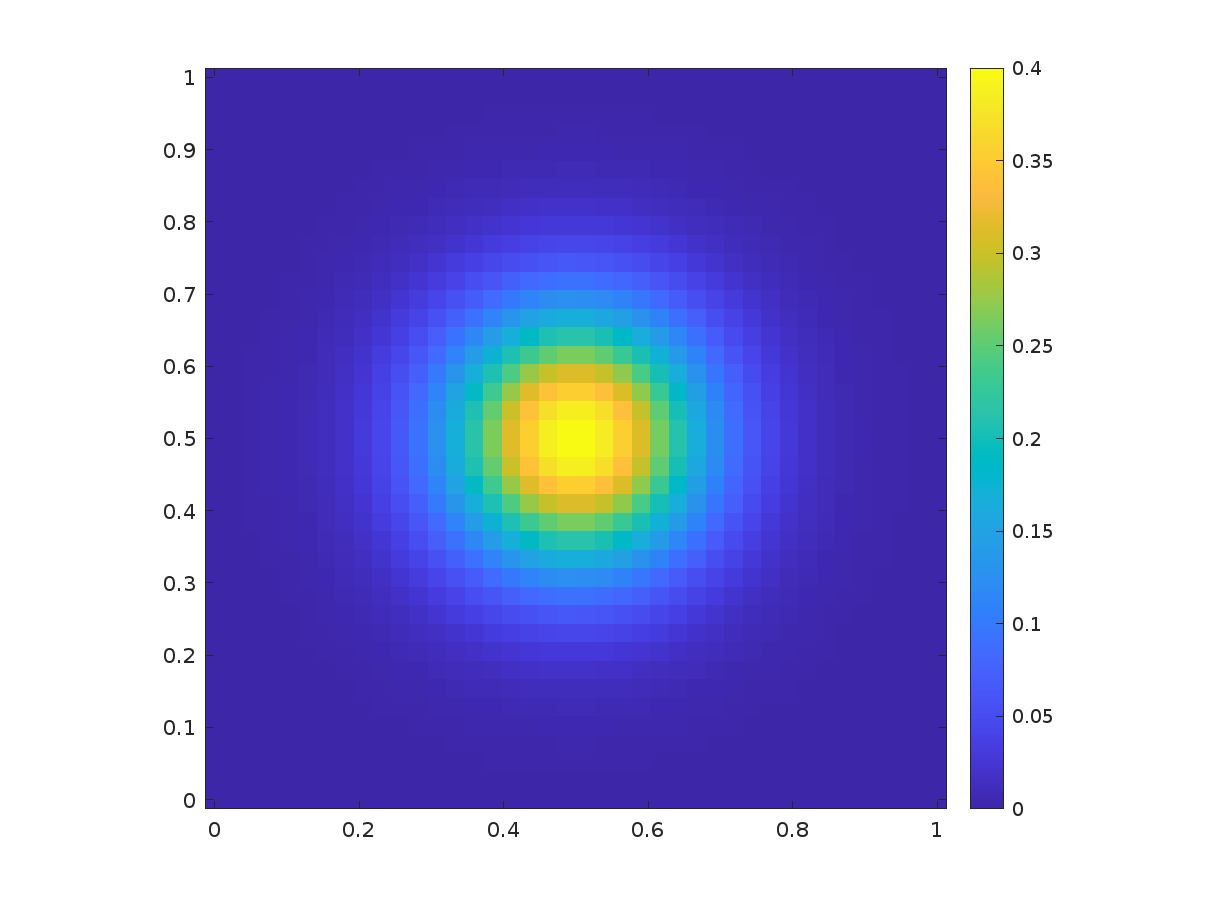}
  \quad
  \includegraphics[width=5.5cm]{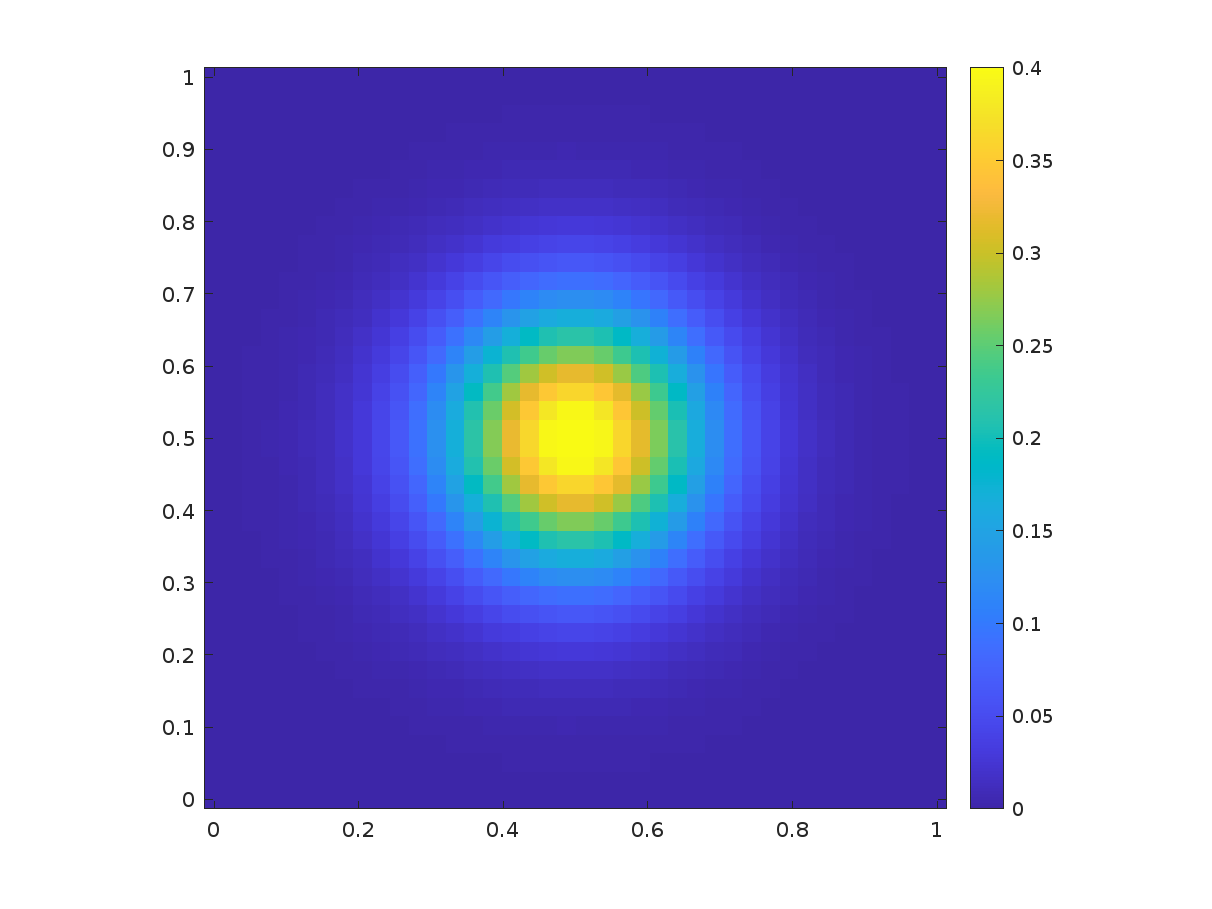}
  \quad
  \includegraphics[width=5.5cm]{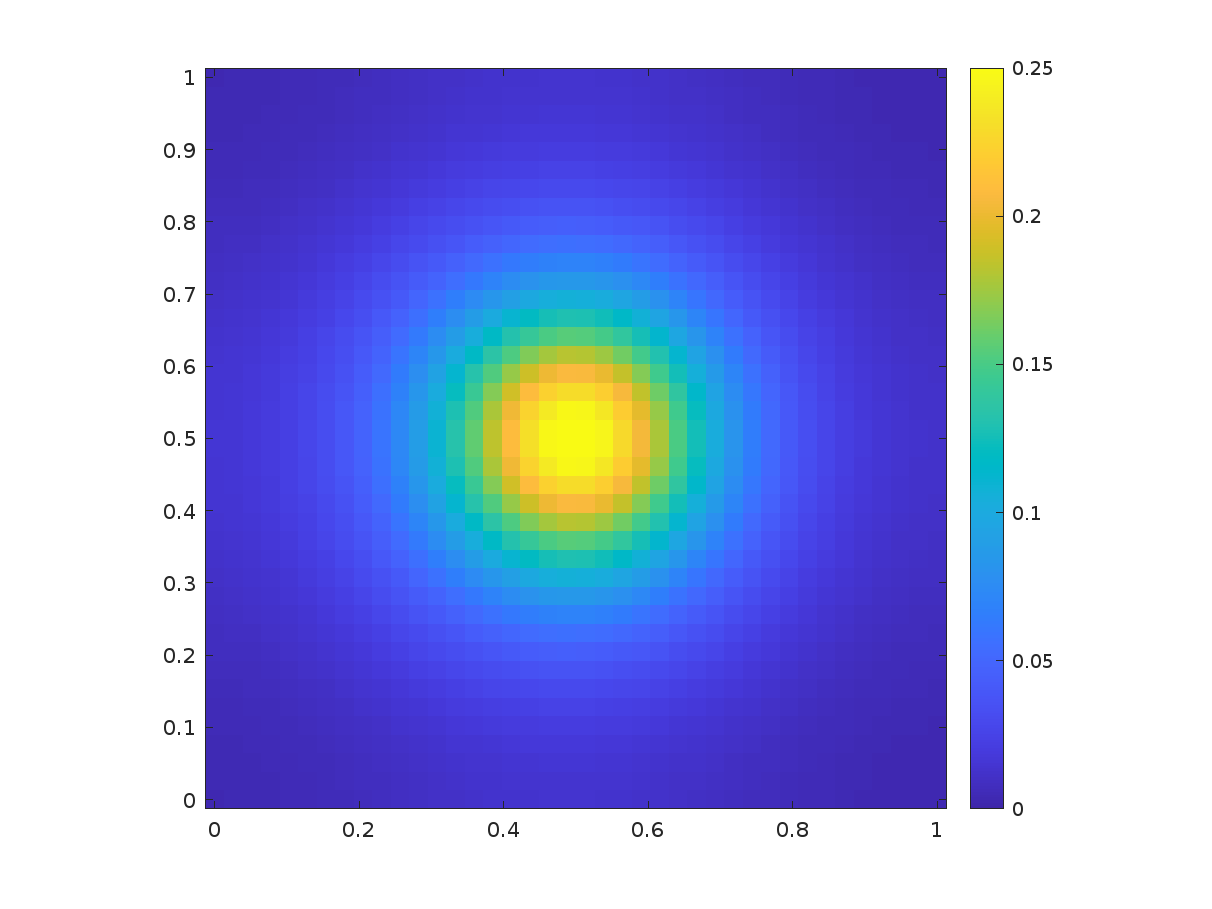}
  \quad
  \includegraphics[width=5.5cm]{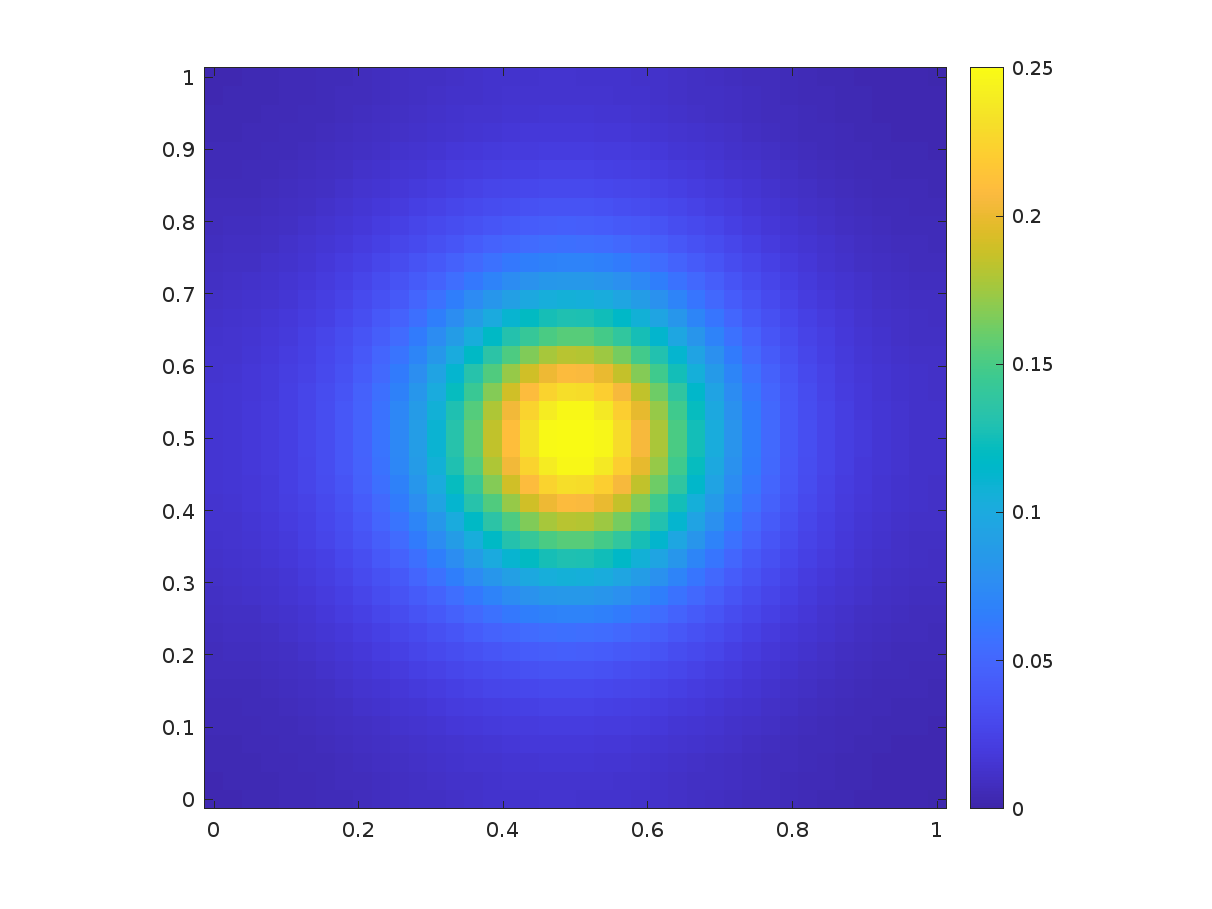}
  \quad
  \includegraphics[width=5.5cm]{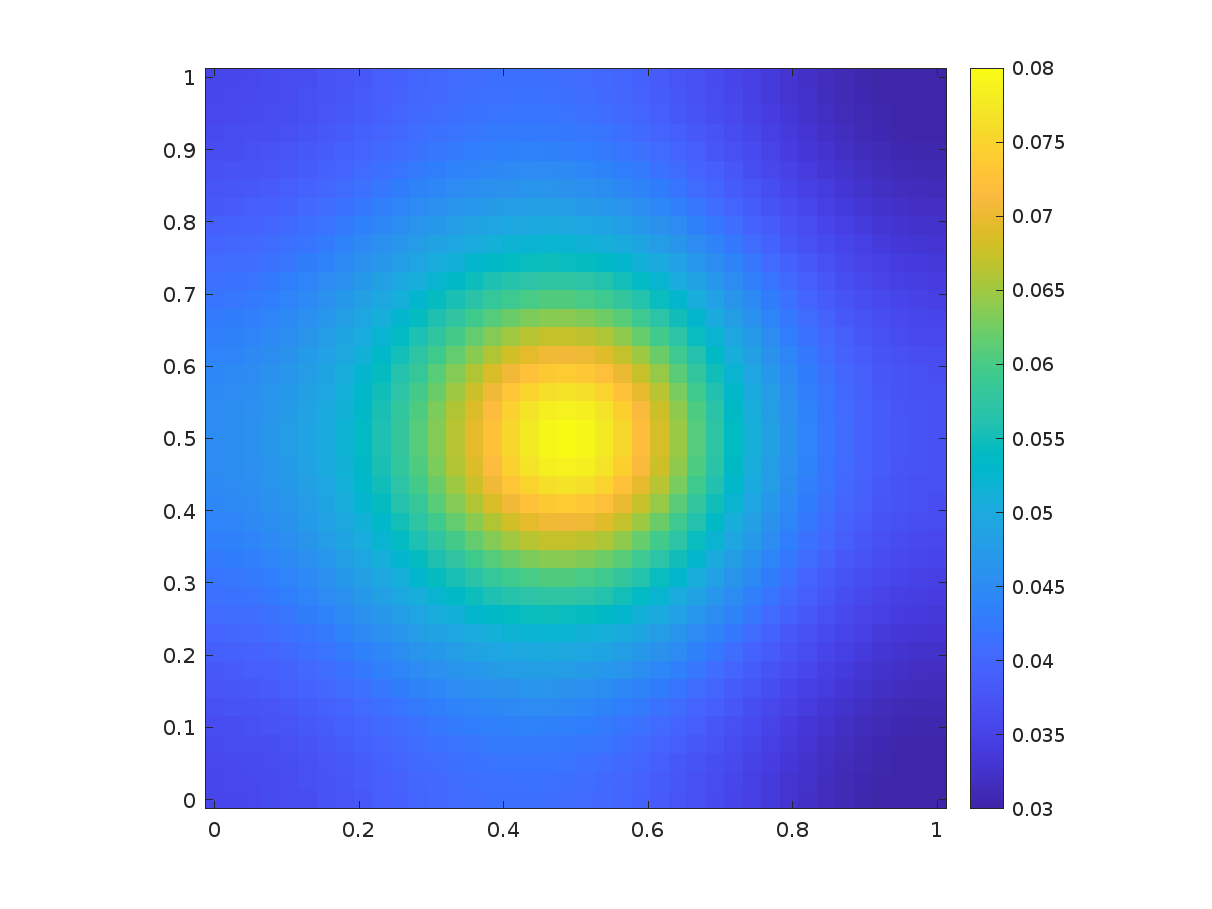}
  \quad
  \includegraphics[width=5.5cm]{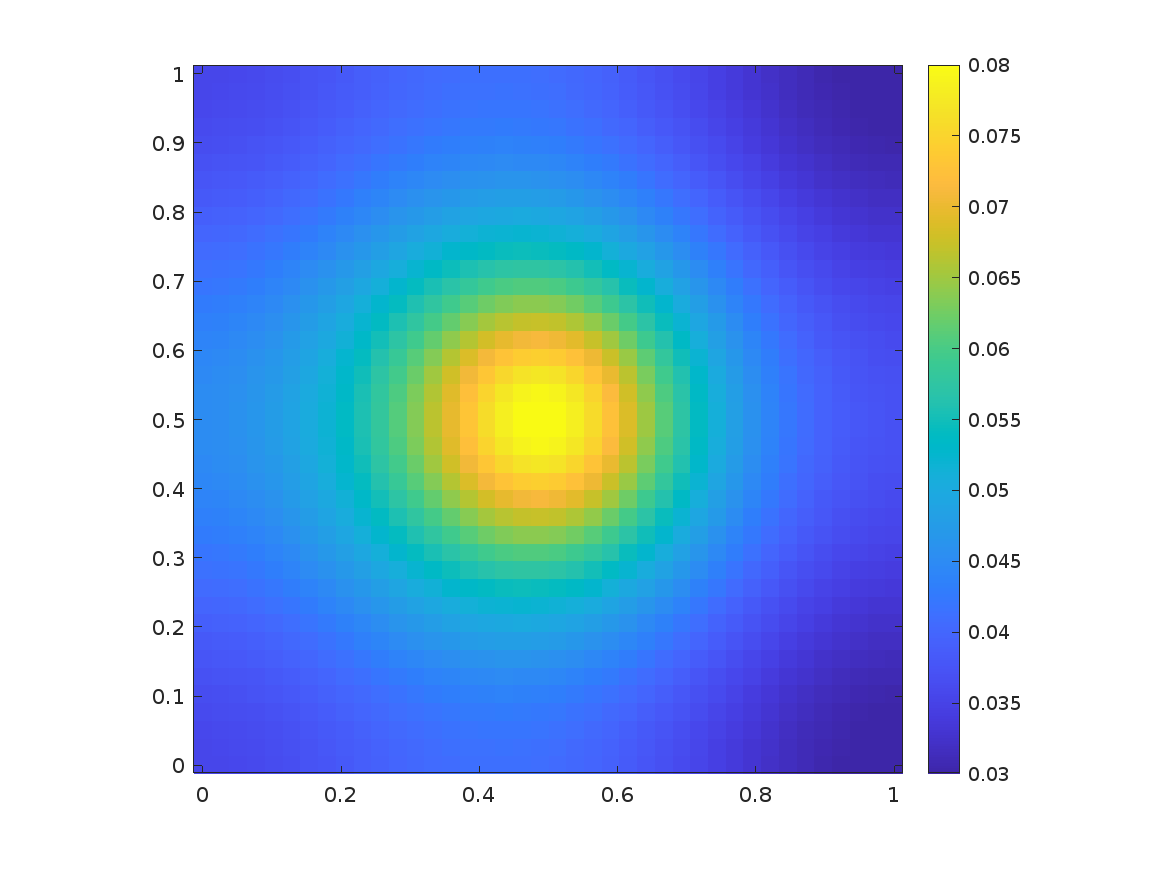}
  \quad
  \includegraphics[width=5.5cm]{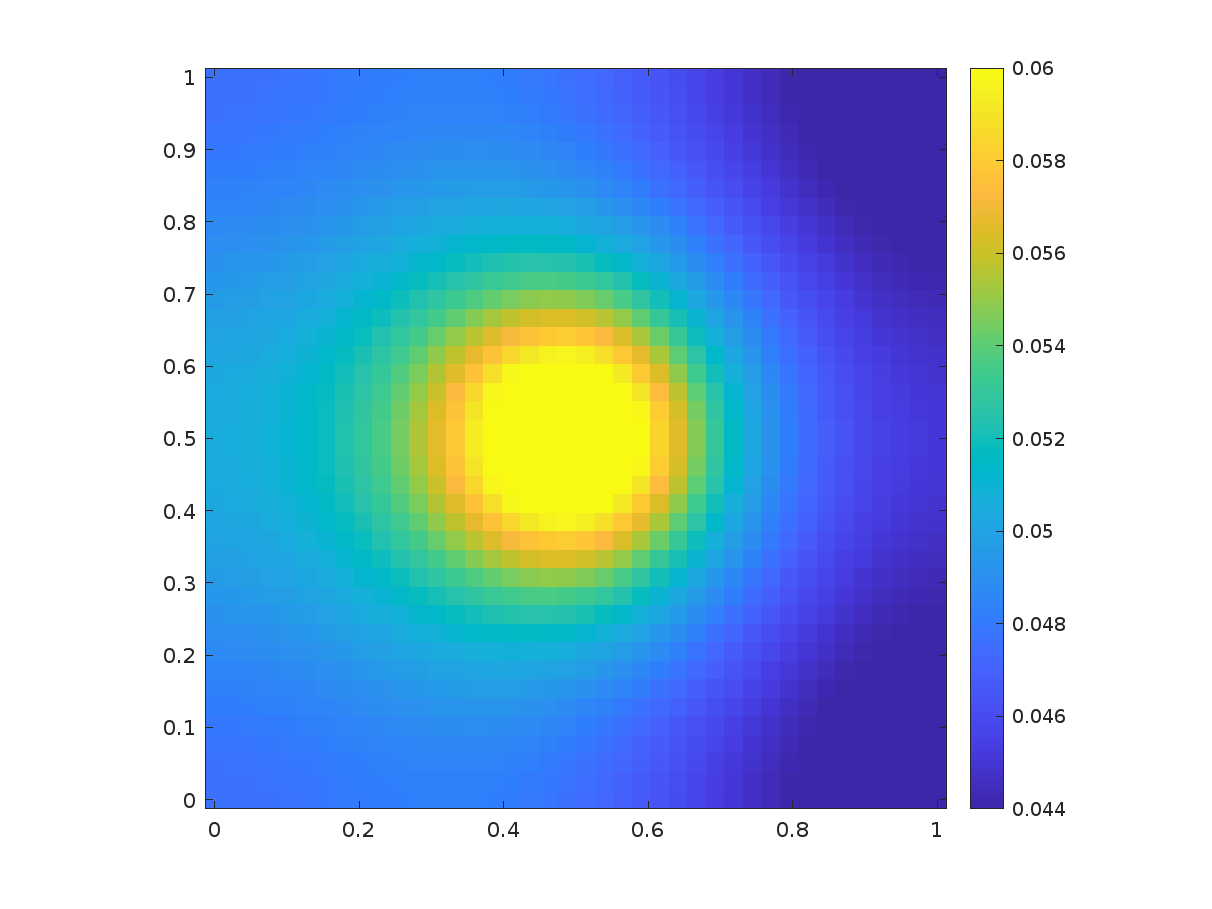}
  \quad
  \includegraphics[width=5.5cm]{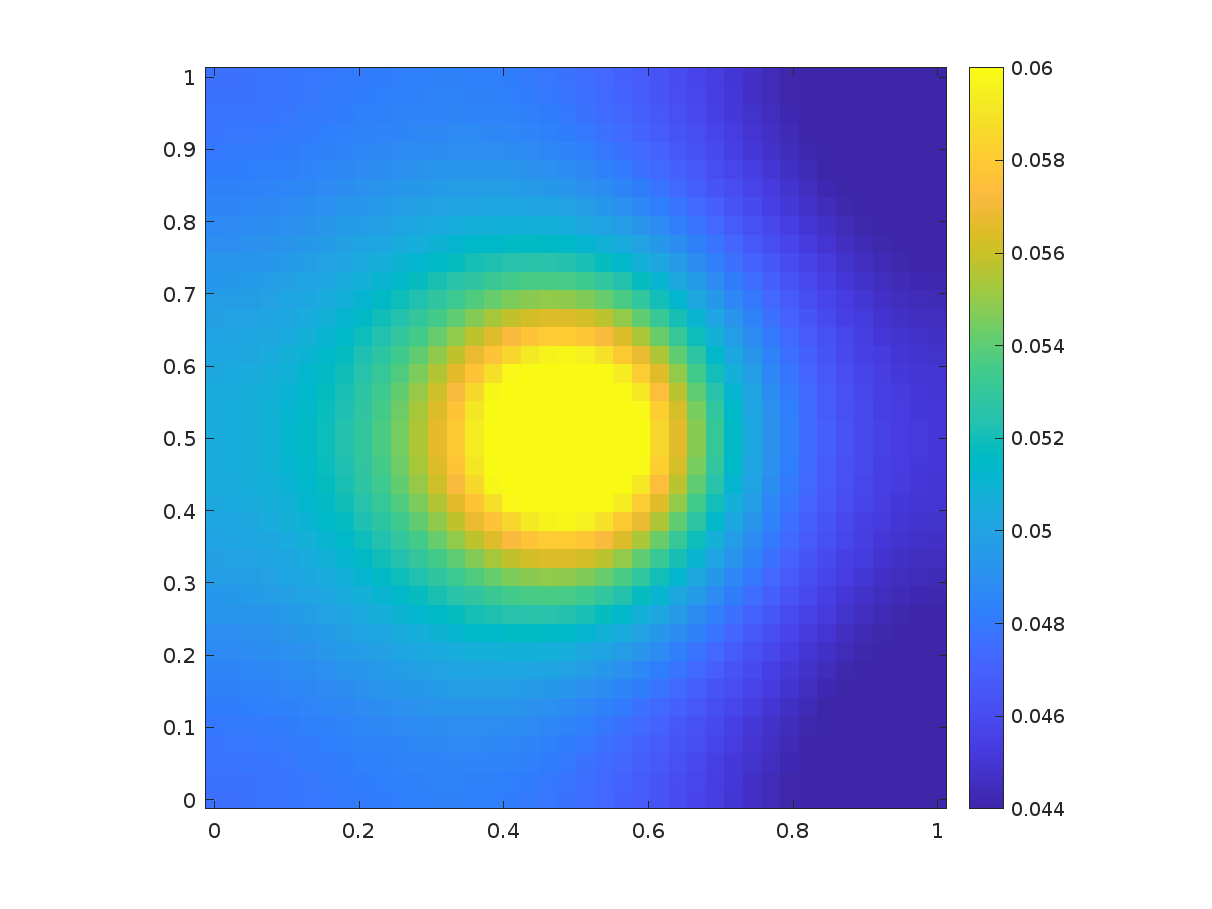}
  \quad
  \includegraphics[width=5.5cm]{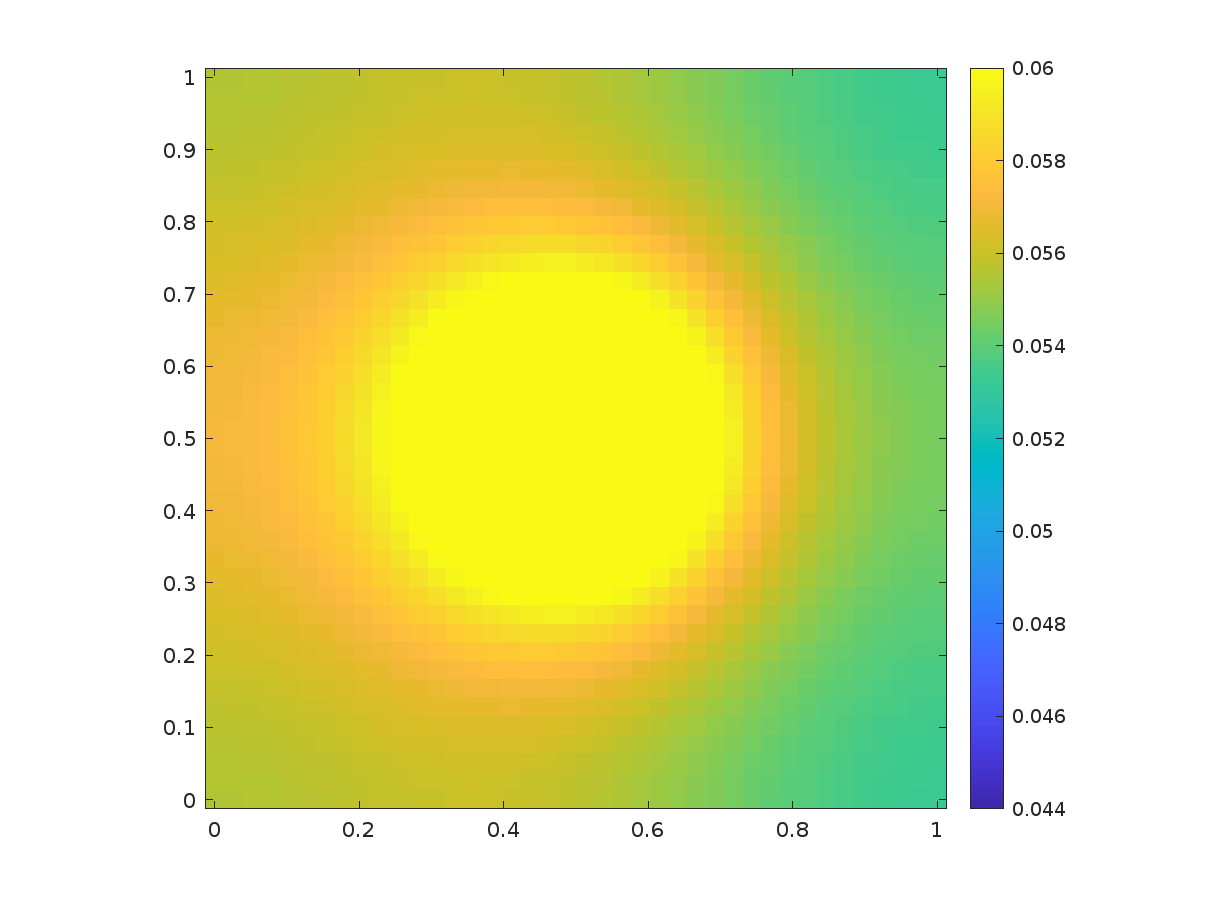}
  \quad
  \includegraphics[width=5.5cm]{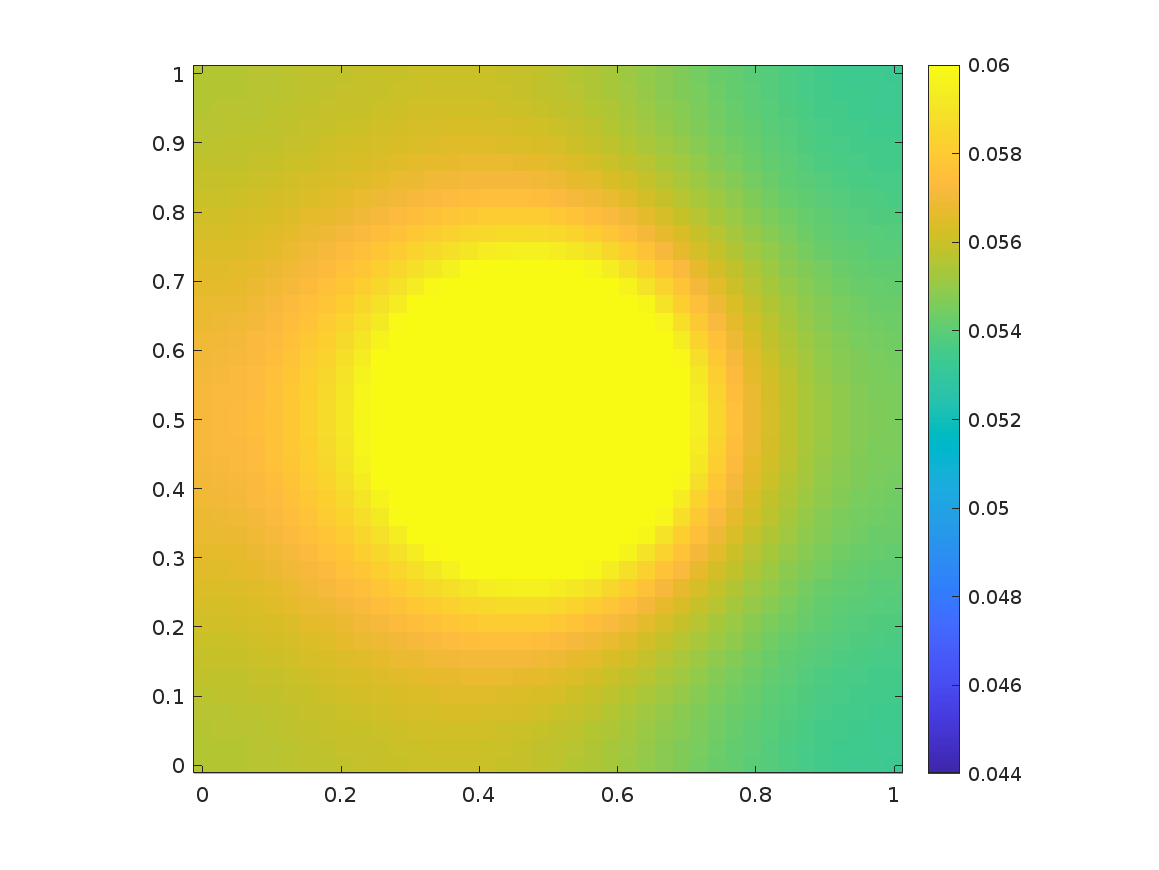}
  \caption{Solutions of concentration when $H=1/40$ for Case 3 in Example 2. First column: multiscale solution $C_1$ at $t=0.02$, $0.1$, $0.5$, $1$, $2$. Second column: reference averaged solution in $\Omega_1$ at the corresponding time instants.}
  \label{fig:Example2_Case3_U1}
\end{figure}

\begin{figure}
  \centering
  \includegraphics[width=5.5cm]{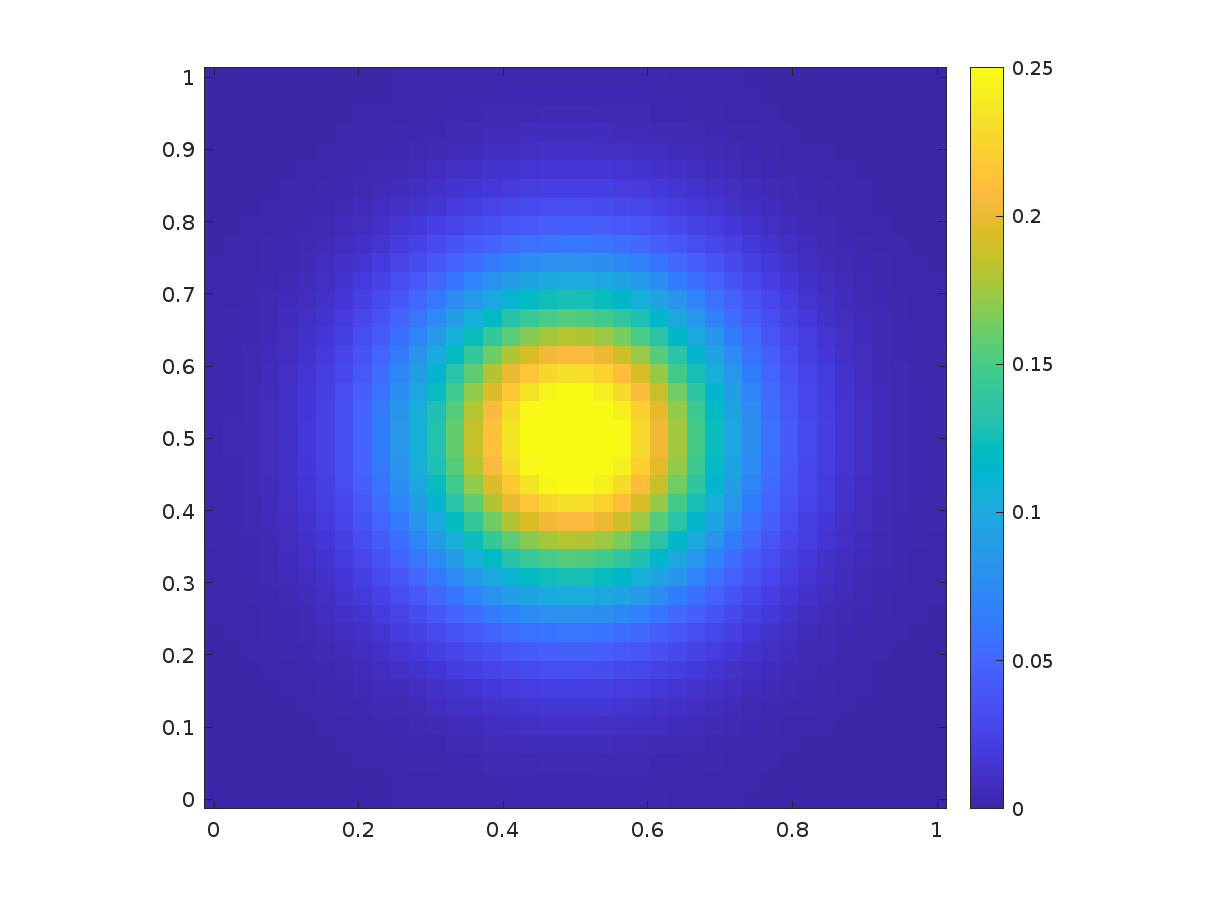}
  \quad
  \includegraphics[width=5.5cm]{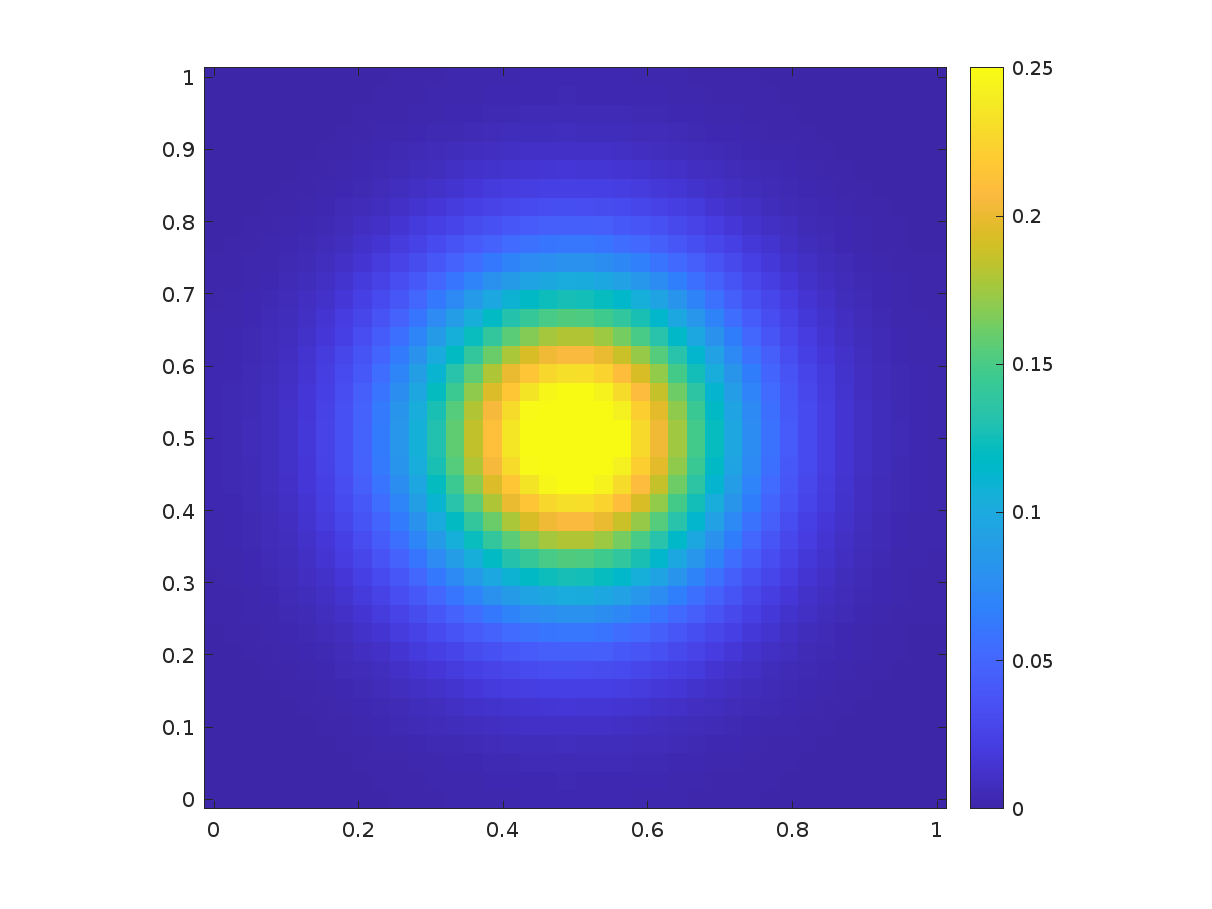}
  \quad
  \includegraphics[width=5.5cm]{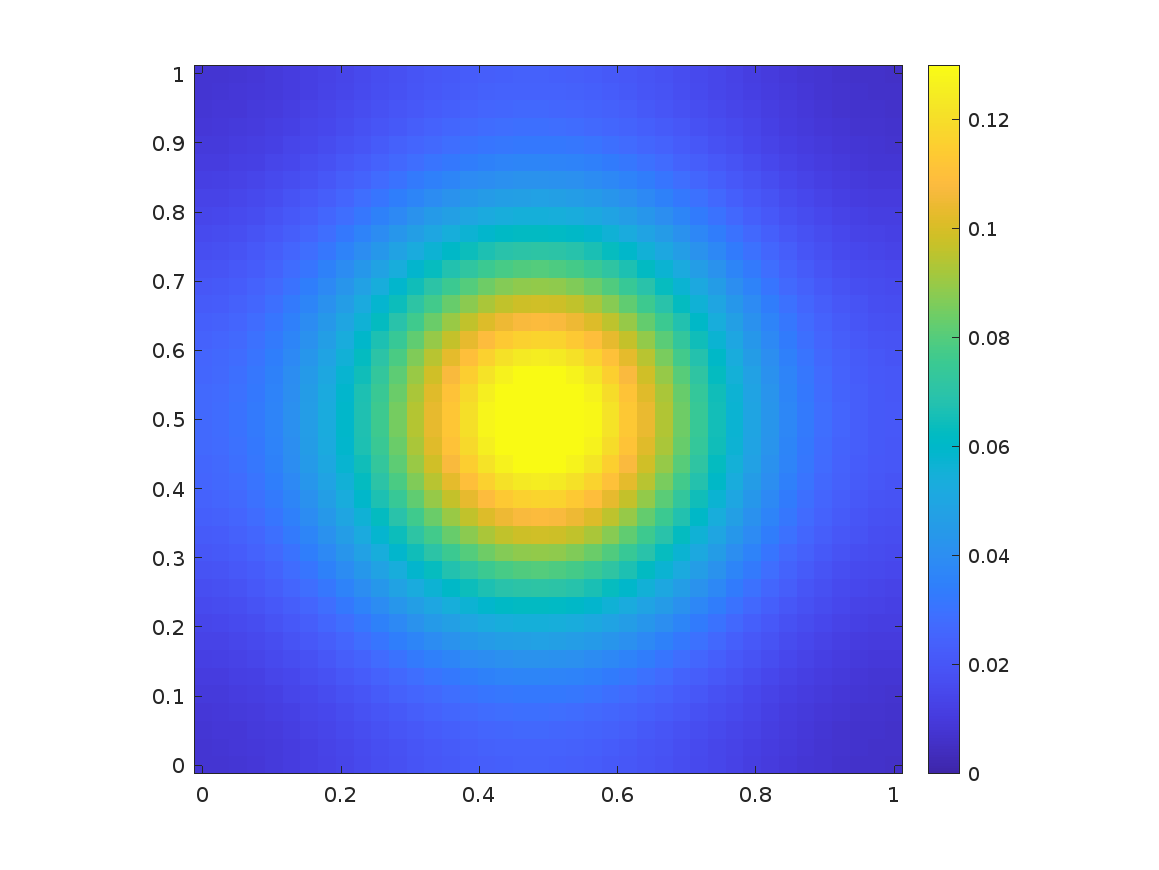}
  \quad
  \includegraphics[width=5.5cm]{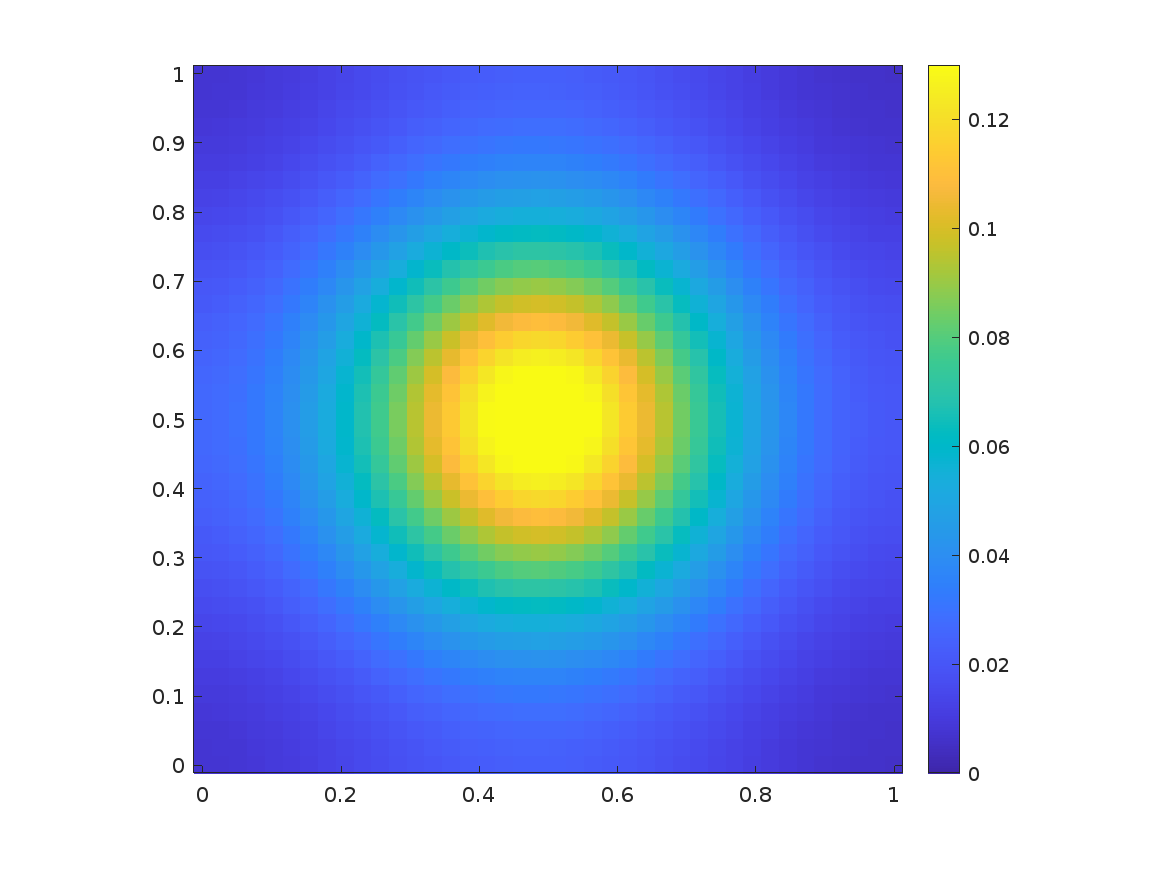}
  \quad
  \includegraphics[width=5.5cm]{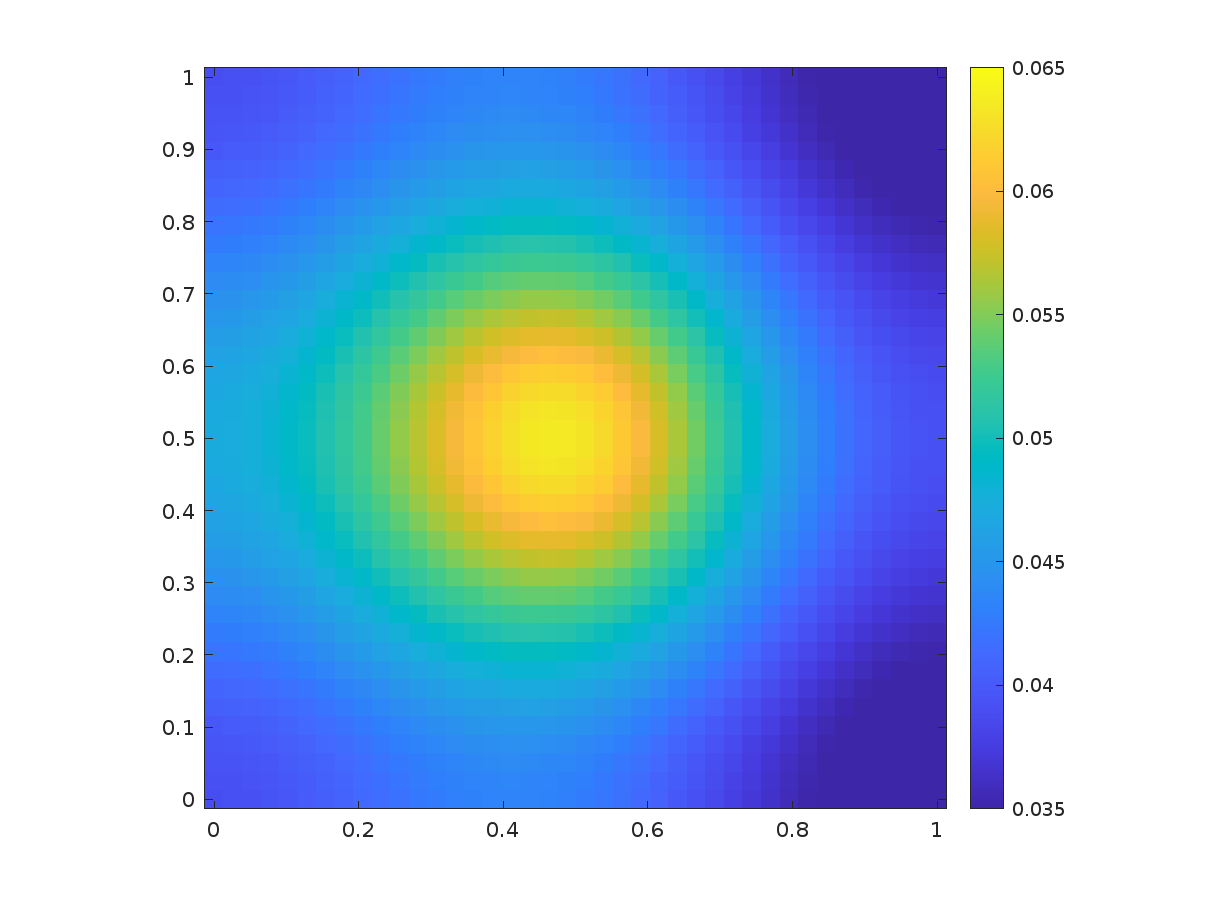}
  \quad
  \includegraphics[width=5.5cm]{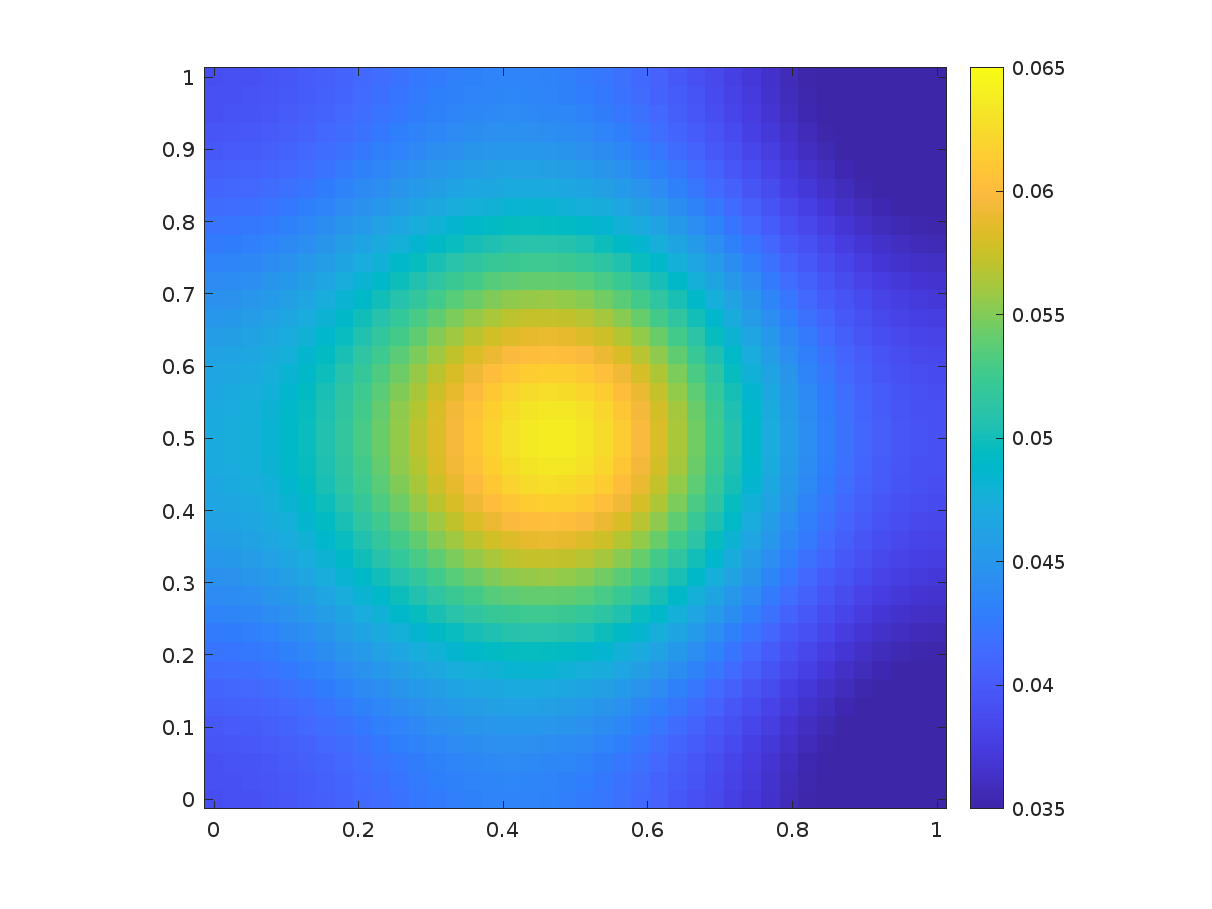}
  \quad
  \includegraphics[width=5.5cm]{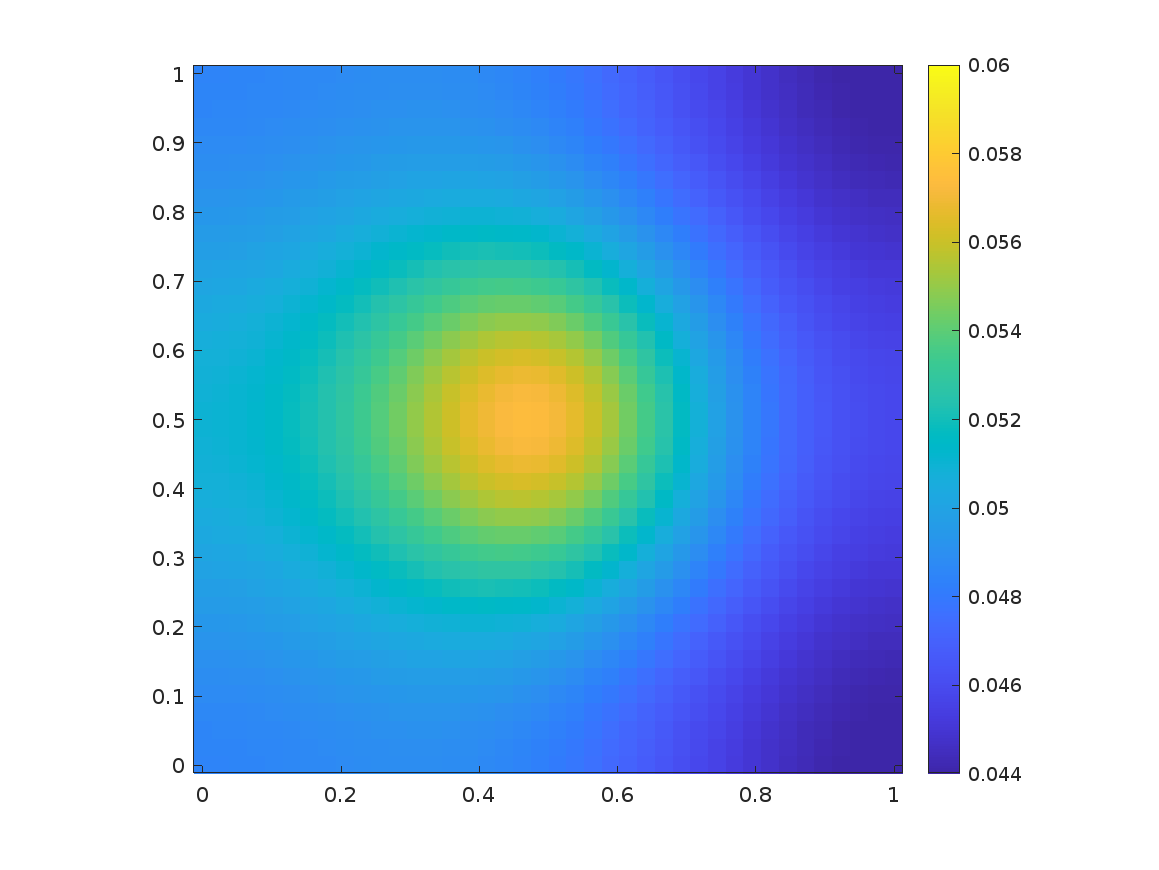}
  \quad
  \includegraphics[width=5.5cm]{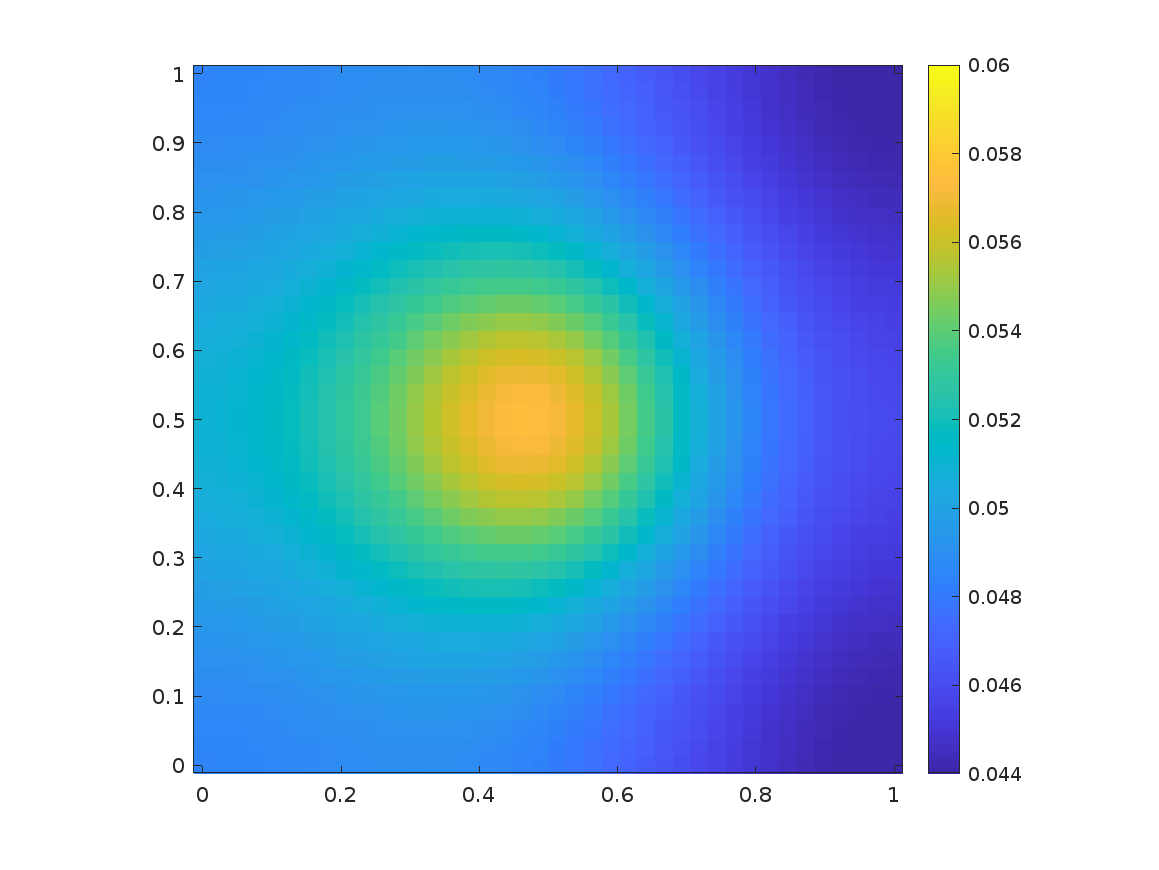}
  \quad
  \includegraphics[width=5.5cm]{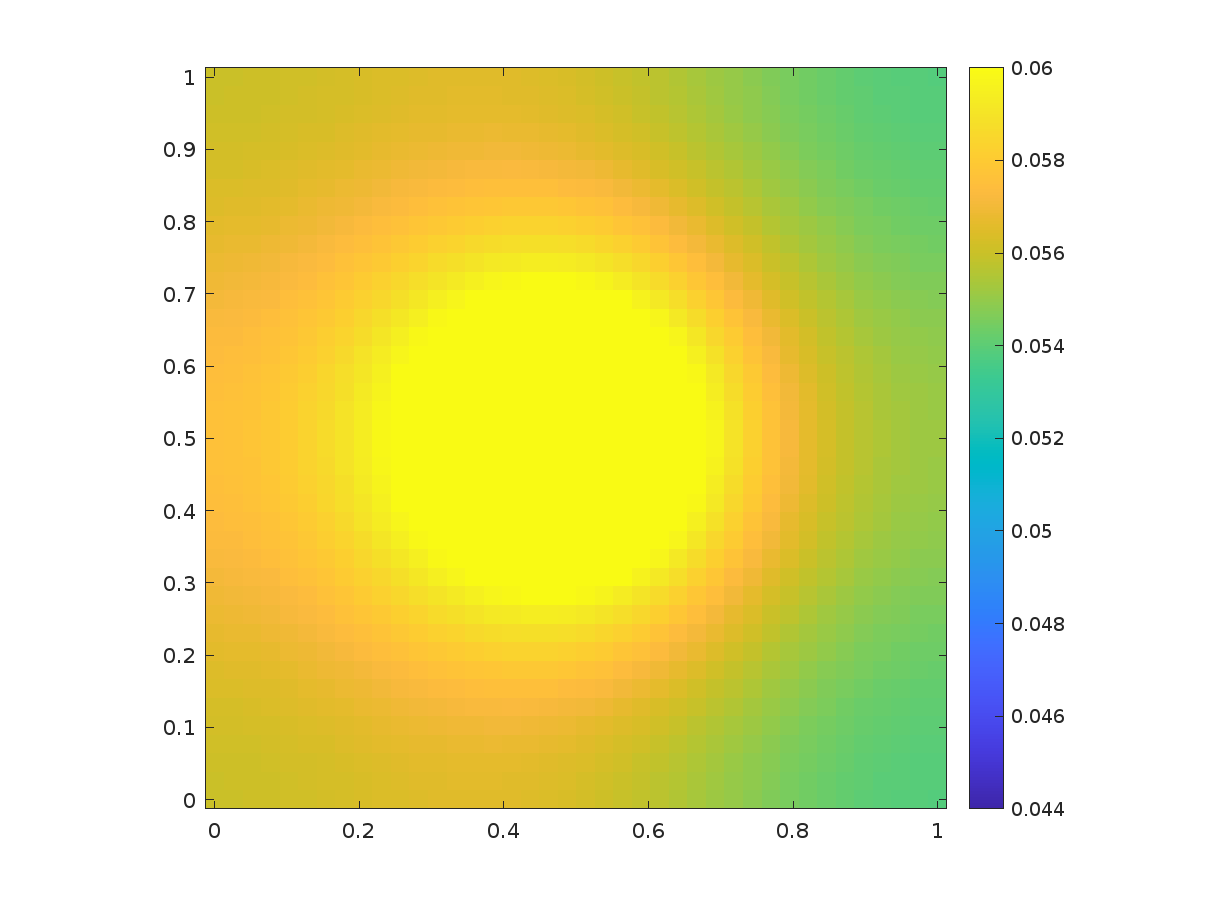}
  \quad
  \includegraphics[width=5.5cm]{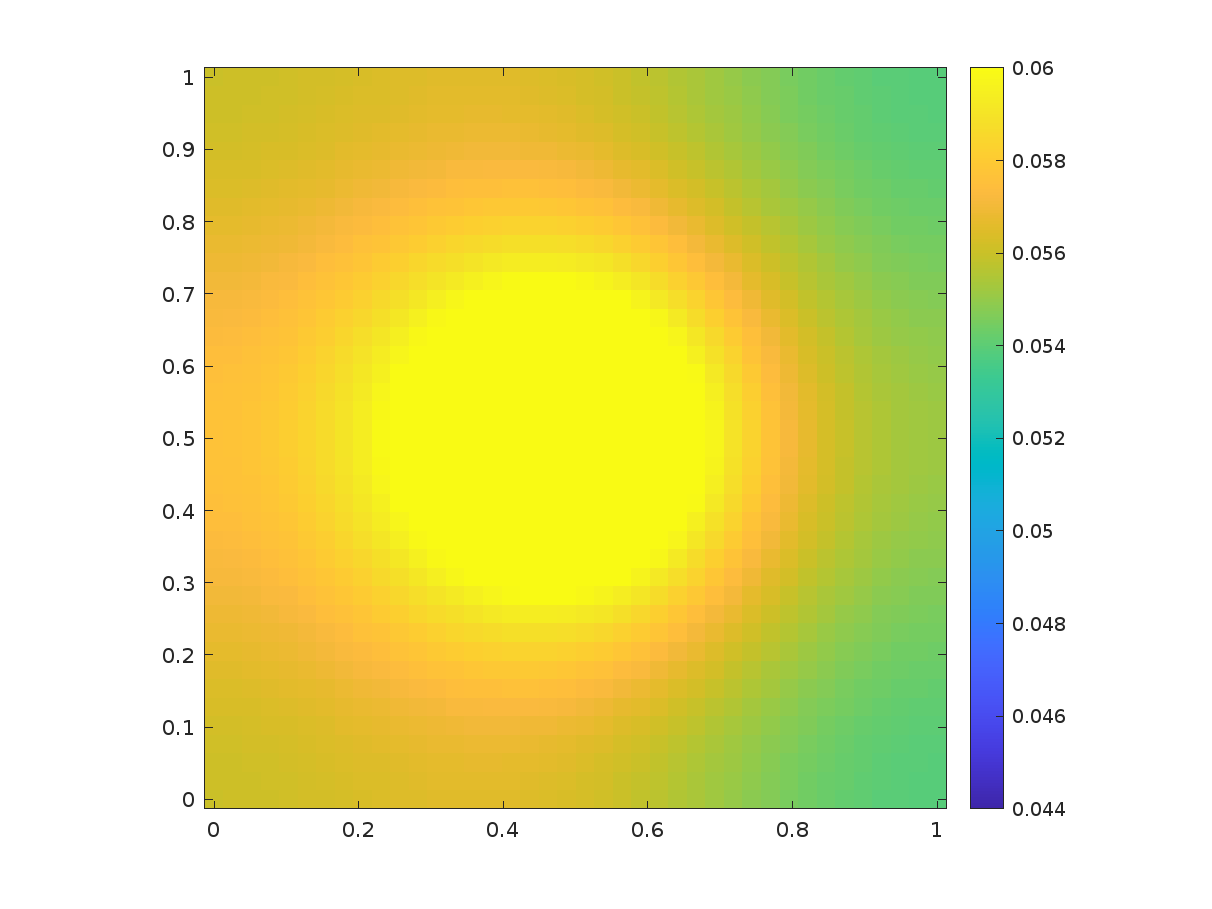}
  \caption{Solutions of concentration when $H=1/40$ for Case 3 in Example 2. First column: multiscale solution $C_2$ at $t=0.02$, $0.1$, $0.5$, $1$, $2$. Second column: reference averaged solution in $\Omega_2$ at the corresponding time instants.}
  \label{fig:Example2_Case3_U2}
\end{figure}

Again, we apply our algorithm to different boundary conditions corresponding to Cases 1, 2, and 3 in Example 1 and present the errors in Tables \ref{tab:Example2_Case1_error}, \ref{tab:Example2_Case2_error}, and \ref{tab:Example2_Case3_error}, respectively. Numerical solutions for Case 3 when $H=1/40$ are depicted in Figures \ref{fig:Example2_Case3_U1} and \ref{fig:Example2_Case3_U2}, where the color bars have been adjusted to display the results more clearly. One can observe that the transport of the concentration for Case 3 in Example 2 is isotropic and more rapid than in Example 1, which is a consequence of the symmetry and porosity of the circular field. We also note that the numerical results for all three cases are in tight agreement with reference solutions.

\begin{table}
\caption{Relative $L^2$ errors at $t=0.02$, $0.1$, $0.5$, $1$, $2$ for Case 1 in Example 2. Left: $H=1/20$ and $l=6$. Right: $H=1/40$ and $l=8$.}
\centering
\begin{tabular}{ccc}
   \toprule
   $t$      & $e^{(1)}(t)$    &   $e^{(2)}(t)$ \\
   \midrule
   $0.02$ & $4.04\%$ & $1.90\%$ \\
   $0.1$ & $2.07\%$ & $1.75\%$ \\
   $0.5$ & $1.21\%$ & $0.91\%$ \\
   $1.0$ & $0.87\%$ & $0.79\%$ \\
   $2.0$ & $0.60\%$ & $0.45\%$ \\
   \bottomrule 
\end{tabular}
\qquad
\begin{tabular}{ccc}
   \toprule
   $t$      & $e^{(1)}(t)$    &   $e^{(2)}(t)$ \\
   \midrule
   $0.02$ & $1.74\%$ & $0.44\%$ \\
   $0.1$ & $0.19\%$ & $0.73\%$ \\
   $0.5$ & $0.41\%$ & $0.21\%$ \\
   $1.0$ & $0.28\%$ & $0.22\%$ \\
   $2.0$ & $0.20\%$ & $0.13\%$ \\
   \bottomrule 
\end{tabular}
\label{tab:Example2_Case1_error}
\end{table}

\begin{table}
\caption{Relative $L^2$ errors at $t=0.02$, $0.1$, $0.5$, $1$, $2$ for Case 2 in Example 2. Left: $H=1/20$ and $l=6$. Right: $H=1/40$ and $l=8$.}
\centering
\begin{tabular}{ccc}
   \toprule
   $t$      & $e^{(1)}(t)$    &   $e^{(2)}(t)$ \\
   \midrule
   $0.02$ & $4.04\%$ & $1.91\%$ \\
   $0.1$ & $2.07\%$ & $1.76\%$ \\
   $0.5$ & $1.20\%$ & $0.90\%$ \\
   $1.0$ & $0.86\%$ & $0.78\%$ \\
   $2.0$ & $0.59\%$ & $0.45\%$ \\
   \bottomrule 
\end{tabular}
\qquad
\begin{tabular}{ccc}
   \toprule
   $t$      & $e^{(1)}(t)$    &   $e^{(2)}(t)$ \\
   \midrule
   $0.02$ & $1.74\%$ & $0.45\%$ \\
   $0.1$ & $0.20\%$ & $0.73\%$ \\
   $0.5$ & $0.41\%$ & $0.21\%$ \\
   $1.0$ & $0.28\%$ & $0.22\%$ \\
   $2.0$ & $0.21\%$ & $0.13\%$ \\
   \bottomrule 
\end{tabular}
\label{tab:Example2_Case2_error}
\end{table}

\begin{table}
\caption{Relative $L^2$ errors at $t=0.02$, $0.1$, $0.5$, $1$, $2$ for Case 3 in Example 2. Left: $H=1/20$ and $l=6$. Right: $H=1/40$ and $l=8$.}
\centering
\begin{tabular}{ccc}
   \toprule
   $t$      & $e^{(1)}(t)$    &   $e^{(2)}(t)$ \\
   \midrule
   $0.02$ & $4.04\%$ & $1.90\%$ \\
   $0.1$ & $2.06\%$ & $1.80\%$ \\
   $0.5$ & $0.65\%$ & $0.38\%$ \\
   $1.0$ & $0.22\%$ & $0.20\%$ \\
   $2.0$ & $0.17\%$ & $0.16\%$ \\
   \bottomrule 
\end{tabular}
\qquad
\begin{tabular}{ccc}
   \toprule
   $t$      & $e^{(1)}(t)$    &   $e^{(2)}(t)$ \\
   \midrule
   $0.02$ & $1.74\%$ & $0.48\%$ \\
   $0.1$ & $0.21\%$ & $0.78\%$ \\
   $0.5$ & $0.22\%$ & $0.08\%$ \\
   $1.0$ & $0.07\%$ & $0.05\%$ \\
   $2.0$ & $0.04\%$ & $0.04\%$ \\
   \bottomrule 
\end{tabular}
\label{tab:Example2_Case3_error}
\end{table}

\FloatBarrier
\section{Conclusion}
\label{sec:conclusion}

In this work, we have applied the multicontinuum homogenization method to develop a new multicontinuum model for the coupled flow and transport equations. Coupled constraint cell problems were carefully formulated to obtain the localizable multiscale basis functions, which can capture the heterogeneity and high contrast properties. Macroscopic variables were introduced to represent the local averages of solutions in each continuum. By performing the multicontinuum expansions to both the flow and transport equations, we have achieved a macroscopic system consisting of homogenized elliptic equations and convection-diffusion-reaction equations. 

We have conducted a number of numerical experiments to verify the obtained multicontinuum model. We have considered high-contrast layered and circular heterogeneous media with homogeneous and inhomogeneous Dirichlet boundary conditions and mixed Dirichlet-Neumann boundary conditions. Numerical results have shown that the proposed model allows us to approximate the reference averaged solutions for all considered cases accurately. The computed errors demonstrate the convergence of the proposed method on the coarse grid size.

\bibliographystyle{unsrt}
\bibliography{ref}

\end{document}